\documentclass{amsproc}
\usepackage{amsmath}
\usepackage{amssymb}
\usepackage{epsf}
\usepackage[scanall]{psfrag}
\usepackage{graphicx}

\newtheorem{theorem}{Theorem}[section]
\newtheorem{lemma}[theorem]{Lemma}
\theoremstyle{definition}
\newtheorem{defn}[theorem]{Definition}

\newtheorem{proposition}[theorem]{Proposition}
\newtheorem{corollary}[theorem]{Corollary}

\theoremstyle{remark}
\newtheorem{rem}[theorem]{Remark}

\numberwithin{equation}{section}

\numberwithin{equation}{section}

\renewcommand{\a}{\alpha}

\renewcommand{\d}{\delta}

\renewcommand{\L}{\Lambda}

\renewcommand{\O}{\Omega}
\renewcommand{\o}{\omega}

\renewcommand{\div}{\mathop{\mathrm{div}}\nolimits}

\date{}

\numberwithin{equation}{section}
\newtheorem{remark}{Remark}

\newenvironment{notation}{\medskip\noindent{\bf Notation\ }}{\medskip}
\newcommand{\ba}{\boldsymbol\alpha}
\newcommand{\bg}{\boldsymbol\gamma}
\newcommand{\sh}{Schr\"odinger\ }
\newcommand{\ka}{K\"ahler\ }
\newcommand{\RR}{\mathbb{R}}

\newcommand{\tr}{\mathrm{Tr}}
\newcommand{\G}{\Gamma}
\newcommand{\cF}{\mathcal{F}}
\newcommand{\cH}{\mathcal{H}}
\newcommand{\cO}{\mathcal{O}}
\newcommand{\cG}{\mathcal{G}}
\newcommand{\s}{\sigma}

\renewcommand{\d}{\delta}

\newcommand{\D}{\Delta}
\renewcommand{\O}{\Omega}
\newcommand{\p}{\partial}

\newcommand{\NN}{\mathbb{N}}

\newcommand{\e}{\varepsilon}
\newcommand{\dist}{\,\mathrm{dist}\,}

\newcommand{\nint}{\int \kern-1.13em {\begin{turn}{-20}$\bigm/$%
\end{turn}}\!}
\newcommand\notint{{{\,\int \kern-1.01em \raise1pt\hbox{{\begin{turn}{-30}$/$%
\end{turn}\!\!}}}}}

\newcommand{\id}{\mathrm{Id}}

\renewcommand{\a}{\alpha}
\renewcommand{\div}{\,\mathrm{div}\,}

\renewcommand{\L}{\Lambda}

\renewcommand{\o}{\omega}

\newcounter{mnotecount}[section]

\newcommand{\rmnote}[1]{}

\begin{document}

\title{The Cauchy problem for Schr\"odinger flows into K\"ahler manifolds}
\author[C. Kenig]{Carlos Kenig}
\thanks{CK Partially supported by the NSF under Grant DMS-0456583}
\address{University of Chicago}

\author[T. Lamm]{Tobias Lamm}
\thanks{TL Partially supported by a PIMS Postdoctoral Fellowship}
\address{University of British Columbia}

\author[D. Pollack]{Daniel Pollack}
\thanks{DP Partially supported by the NSF under Grant DMS-0305048}
\address{University of Washington}

\author[G. Staffilani]{Gigliola Staffilani}
\thanks{GS Partially supported by the NSF under Grant DMS-0330731}
\address{MIT}

\author[T. Toro]{Tatiana Toro}
\thanks{TT Partially supported by the NSF under Grant DMS-0244834 and an ADVANCE TSP grant at the University of Washington}
\address{University of Washington}

\begin{abstract}We prove local well-posedness of the Schr\"odinger flow from $\RR^n$ into a
compact K\"ahler manifold $N$ with initial data in $H^{s+1}(\RR^n,N)$ for
$s\geq\left[\frac{n}{2}\right]+4$. 
\end{abstract}

\maketitle


\section{Introduction}

We consider maps
\[
u:\RR^n\to N
\]
where $N$ is a $k$-dimensional compact K\"ahler manifold with complex
structure $J$ and K\"ahler form $\omega$ (so that $\o$ is a nondegenerate,
skew-symmetric two-form). Thus $N$ is a complex manifold and $J$ is an
endomorphism of the tangent bundle whose square, at each point, is minus the
identity. $N$ has a Riemannian metric $g$ defined by
\[
g(\cdot,\cdot)=\o(\cdot,J\cdot).
\]
The condition that $N$ is K\"ahler is equivalent to assuming that $\nabla
J=0$ where $\nabla$ is the Levi-Civita covariant derivative with respect to
$g$.
The energy of a map $u$ is defined by
\[
E(u)=\frac{1}{2}\int_{\RR^n}|du|^2dx\]
where the energy density $|du|^2$ is simply the trace with respect to the
Euclidean metric of the pullback of the metric $g$ under $u$, $|du|^2=\tr\,
u^\ast(g)$. In local coordinates we have
\[
|du|^2(x)=\sum^n_{\a=1}g_{ij}(u(x))\frac{\p u^i}{\p x^\a}\frac{\p u^j}{\p
x^\a}.
\]
(We use the Einstein summation convention and sum over repeated indices.)

The $L^2$-gradient of $E(u)$ is given by minus the tension of the map,
$-\tau(u)$, $\tau(u)$ is a vector field on $N$ which can be expressed in
local coordinates as $\tau(u)=(\tau(u)^1,\ldots,\tau(u)^k)$ with
\begin{equation}\label{int-eqn1}
\tau(u)^i=\Delta u^i+\sum^n_{\a=1}\G^i_{jk}(u) \frac{\p u^j}{\p x^\a}
\frac{\p u^k}{\p x^\a}\qquad\mbox{for}\qquad i=1,\ldots, k
\end{equation}
where $\G^i_{jk}(u)$ are the Christoffel symbols of the metric $g$ at
$u(x)$. Critical points of the energy are Harmonic maps and are
characterized by the equation $\tau(u)=0$. The foundational result on the
existence of harmonic maps is due to Eells and Sampson \cite{ES} and is
achieved by studying the harmonic map flow
\[
\frac{\p u}{\p t} = \tau(u)
\]
which is simply the gradient flow for the energy functional on the space of
maps. Eells and Sampson proved the existence of harmonic maps as stationary points of 
this flow when the domain is a compact manifold and the target is a compact manifold
of non-positive curvature.  In our setting, 
the symplectic structure on $N$ induces a symplectic structure on the
space of maps. Let $X_s=H^s(\RR^n, N)$ be the Sobolev space of maps between
$\RR^n$ and $N$ as defined below. For $s\ge \frac{n}{2}+1$, $X_s$ is a
Banach manifold with a symplectic structure $\O$ induced from that of
$(N,\o)$ as follows. 
The tangent space to $X_s$ at a map $u$ is identified with sections of
the pull-back tangent bundle over $\RR^n$. We let $\G(V)$ denote the space
of sections of the bundle $V$, for example $du\in \G(T^s(\RR^n)\otimes
u^{-1}(TN))$. For $\s,\mu\in\G(u^{-1}(TN))=T_uX_s$ we define
\[
\O(\s,\mu)=\int_{\RR^n}\o(\s,\mu)dx.
\]
In this setting we are interested in the Hamiltonian flow for the 
energy functional
$E(\cdot)$ on $(X_s,\O)$. This is the \sh flow which takes
the form
\begin{equation}\label{int-eqn2}
\frac{\p u}{\p t}=J(u)\tau(u).
\end{equation}
This natural geometric motivation for the flow (\ref{int-eqn2}) was 
elucidated in \cite{DW1998}.

A key aspect of our approach to understanding the flow (\ref{int-eqn2}) is to
isometrically embed $N$ in some Euclidean space $\RR^p$ and study
``ambient'' flows of maps from $\RR^n$ to $\RR^p$ which are related to
(\ref{int-eqn2}). This is also central to the Eells-Sampson treatment of
the harmonic map flow. Toward this end we use the Nash embedding theorem to
assume that we have an isometric embedding
\[
w:(N,g)\to(\RR^p,\d).
\]
Using this we can now define $H^s(\RR^n;N)$, the $L^2$-based Sobolev spaces
of maps from $\RR^n$ to $N$ as follows. Note that since the domain is
noncompact some care must be taken even when $s=0$.

\begin{defn}\label{int-defn1}
For $s\ge 1$ let
\begin{eqnarray}
H^s(\RR^n; N)&=&\{u:\RR^n\to\RR^p: u(x)\in N\mbox{ a.e. and }\exists\ 
y_u\in N\mbox{ such that} \\
&& v-w(y_u)\in L^{2}(\RR^n;\RR^p), \, \partial v\in H^{s-1}(\RR^n;\RR^p),\nonumber\\
&&\mbox{ where }v=w\circ u\}.\nonumber
\end{eqnarray}
\end{defn}
With this definition in mind we can state our main result.
\begin{theorem}\label{int-thm1}
Given $\beta\geq 0$, the initial value problem
\begin{equation}\label{int-eqn3}
\left\{\begin{array}{rcl}
\frac{\p u}{\p t} & = & J(u)\tau(u)+\beta\tau(u) \\
u(0) & = & u_0 \end{array}\right.
\end{equation}
for the generalized \sh flow 
has a solution whenever the initial data $u_0\in
H^{s+1}(\RR^n,N)$ for $s\ge\left[\frac{n}{2}\right]+4$.
Moreover (\ref{int-eqn3}) 
is locally well posed in 
$H^{s+1}(\RR^n,N)$ for $s\ge\left[\frac{n}{2}\right]+4$.
\end{theorem}

The question of the local and global well-posedness of equation (\ref{int-eqn3}) with
data in Sobolev spaces has been previously studied by many authors (see
\cite{DW1998, DW2001, D, PWW2000, PWW2001, SSB, CSU, NSU, NSU-errat, TU, WW, McG,  
McG1, K, K-N, K-K}).
A common feature in most existence results for smooth solutions of 
Schr\"odinger maps is that they are obtained by using  the  energy method. This method 
consists in finding an appropriate  regularizing equation which approximates
the \sh flow, and for which smooth solutions exist.  One then proves that the regularizing equations 
satisfy a priori bounds in certain Sobolev norms, independent of the approximation, 
and that they converge to a solution of the original equation. The differences in the distinct results and
proofs lie in the type of  regularization used. 

Ding and Wang \cite{DW2001}  established a similar result
to Theorem 1.2 for $s\ge \left[\frac{n}{2}\right]+3$. Their work proceeds by
direct study of equation (\ref{int-eqn3}) with $\beta>0$, with a passage to the limit
for $\beta=0$. Thus the regularizing equation they use is obtained by adding
the second order dissipative term $\beta \tau(u)$.   In this paper we analyze equation (\ref{int-eqn3}) 
by adding a fourth order dissipative term (note that we allow the case $\beta=0$ from the
start).  This term arises naturally in the geometric setting as the first variation of the 
$L^2$-norm of the tension.
We believe that our regularization of (\ref{int-eqn3}) by a fourth order
equation, which is geometric in nature, is of intrinsic and independent
interest. H.\ McGahagan \cite{McG, McG1} in her doctoral dissertation
also proved a version of Theorem 1.2. Her work proceeds by a different
regularization, this time hyperbolic, implemented by adding a term of the
form $-\epsilon\frac{\partial^2 u }{\partial t^2}$ which  transforms the equation into one whose
solutions are wave maps. 

We note that while our existence proof in Theorem 1.2 is different from the ones in \cite{DW2001} and \cite{McG, McG1}, our proof of uniqueness is the same, using parallel transport. In fact, in Appendix A we extend the uniqueness argument in \cite{McG1}, carried out there in the case when $\beta=0$ to the case $\beta \ge 0$, which gives the uniqueness statement in Theorem 1.2. 

Our proof of Theorem 1.2 actually only shows that the mapping $u_0\mapsto
u\in C([0,T], H^{s'+1}(\RR^n,N))$, with $s'<s$, is continuous. However,
one can show, by combining the parallel transport argument with the standard Bona-Smith regularization procedure
(\cite{B-S,I-I,Ke}) that the statement in Theorem 1.2 also holds.

Equations of the type (\ref{int-eqn3}), but with $N$ being Euclidean space are
generally known as derivative Schr\"odinger equations and have been the
object of extensive study recently (see for instance
\cite{KPV,H-O,C,HH,KPV-2,KPV-3,KPRV1,KPRV2}). The results in these
investigations however do not apply directly to (\ref{int-eqn3}) for two reasons. The
first one is the constraint imposed by the target being the manifold $N$. The second one
is that in these works one needs to have data $u_0$ in a weighted Sobolev
space, a condition that we would like to avoid in the study of (\ref{int-eqn3}).

It turns out that for special choices of the target $N$, the equations (\ref{int-eqn3})
are related to various theories in mechanics and physics. They are examples
of gauge theories which are abelian in the case of Riemann surfaces
(K\"ahler manifolds of dimension 1 such as the 2-sphere $S^2$ or
hyperbolic 2-space $H^2$). In the case of the 2-sphere $S^2$, Schr\"odinger
maps arise naturally from the Landau-Lifschitz equations (a $U(1)$-gauge
theory) governing the static as well as dynamical properties of
magnetization \cite{LL,PT}. More precisely, for maps $s:\RR\times\RR^n\to
S^2\hookrightarrow\RR^3$, equation (\ref{int-eqn3}) takes the form
\begin{equation}\label{lan-lip}
\partial_ts=s\times\Delta s,\qquad |s|=1
\end{equation}
which is the Landau-Lifschitz equation at zero dissipation, when only the
exchange field is retained \cite{K-P,PT}. When $n=2$ this equation is also
known as the two-dimensional classical continuous isotropic Heisenberg spin
model (2d-CCIHS); i.e. the long  wave length limit of the isotropic
Heisenberg ferromagnet (\cite{K-P,PT,SSB}). It also occurs as a continuous
limit of a cubic lattice of classical spins evolving in the magnetic field
created by their close neighbors \cite{SSB}. The paper \cite{SSB} contains,
in fact, for the cases $n=1,2$, $N=S^2$ the first local well-posedness
results for equation (\ref{int-eqn3}) or (\ref{lan-lip}) that we are aware of. In \cite{CSU},
Chang-Shatah-Uhlenbeck showed that, when $n=1$, $(\ref{lan-lip})$ is globally (in
time) well-posed for data in the energy space $H^1(\RR^1; S^2)$. When
$n=2$, for either radially symmetric or $S^1$-equivariant maps, 
they show that small energy implies global existence. For global existence results see also \cite{rrs}. In
\cite{NSU,NSU-errat}, the authors show that, when $n=2$, the problem is
locally well-posed in the space $H^{2+\e}(\RR^2;S^2)$, while the existence
was extended to the space $H^{3/2+\e}(\RR^2;S^2)$  in \cite{K} and \cite{K-N},
and the uniqueness to the space $H^{7/4+\epsilon}(\RR^2; S^2)$ in \cite{K-K}.

More recently, in \cite{Be1,IK1}, a direct method, in the case of small data, using fixed point arguments in suitable spaces was introduced. The first global well-posedness result for \eqref{lan-lip} in critical spaces (precisely, global well-posedness for small data in the critical Besov spaces in dimensions $n\ge 3$) was proved, independently in \cite{IK2} and \cite{Be2}. This was later improved to global regularity for small data in the critical Sobolev spaces in dimensions $n\ge 4$ in \cite{BeIK}. Finally, in \cite{BeIKT}, the global well-posedness of \eqref{lan-lip}, for ``small data'' in the critical Sobolev space $H^{\frac{n}{2}}(\RR^n,S^2)$, $n\ge 2$, was proved. 

\begin{rem}
A first version of this paper was posted on arxiv:0511701, in November 2005, by C.K., D.P., G.S. and T.T. The paper was withdrawn in May 2007. The reason for this was that Jesse Holzer, at the time a graduate student of Alex Ionescu at the University of Wisconsin, Madison, discovered an error in the original argument. The source of the error was in the construction of the ambient flow equations, which were introduced in the first version of the paper. The construction of the ambient flow equations resulted in equations of quasilinear type in the leading order term, which could not be solved by the Duhamel principle as was pointed out by Holzer. 

C.K., D.P., G.S. and T.T. are indebted to Jesse Holzer for pointing out this error and to T.L., who constructed new ambient flow equations, whose leading order terms are $\e\Delta^2$, which can then be dealt with directly by Duhamel's principle. The price one pays for this change in the ambient flow equation, is that it turns out to be more difficult to show that if the initial value takes values in the manifold, so does the whole flow. All of this is carried out in Lemma 2.5 and Lemma 2.10.
\end{rem}

\begin{notation}
    We will use $C, c$ to denote various constant, usually depending
    only on $s, \, n$ and the manifold $N$.  In case a constant
    depends upon other quantities, we will try to make  that explicit.
    We use $A\lesssim B$ to denote an estimate of the form $A\leq CB$. 
    \end{notation}


\section{A fourth order parabolic regularization}

The method we employ in order to establish short-time existence to
(\ref{int-eqn3}) is in part inspired by the work of Ding and Wang
\cite{DW1998}. We seek to approximate equation (\ref{int-eqn3}) by a family
(parametrized by $0<\e<1$) of parabolic equations. We establish short time
existence for these systems and use energy methods to show that the time of
existence is independent of $\e$ and obtain $\e$ independent bounds which
allow us to pass to the limit as $\e\to 0$ and thus obtain a solution to
(\ref{int-eqn3}). The regularization we use differs substantially from that
of Ding and Wang because we wish to view the right hand side of 
(\ref{int-eqn3}) as a lower order term (in the regularization) 
so that we can use Duhamel's principle and a contraction
mapping argument to establish and study the existence of our derived
parabolic system.

The energy method we employ ultimately depends on establishing $\e$
independent $L^2$-estimates for the tension, $\tau(u)$ and its derivatives.
This suggests that we regularize (\ref{int-eqn2}) by $\e$ times the gradient
flow for the functional
\[
G(u)=\frac{1}{2}\int_{\RR^n}|\tau(u)|^2dx.
\]

\subsection{Geometric Preliminaries}

We perform many of our computations in the appropriate pull-back tensor
bundles over $\RR^n$. We begin by recalling alternative formulations of the
tension $\tau(u)$ in this setting (see \cite{ES}). First note that $du$
is a closed 1-form with values in $u^{-1}(TN)$. The tension is simply minus
the divergence of the differential of $u$
\[
\tau(u)=-\d du\in \G(u^{-1}(TN))
\]
where $\d$ denotes the divergence operator with respect to the metric $g$.
In particular, this shows that a map $u$ is harmonic if and only if its
differential is a harmonic 1-form. Let $\nabla$ denote the covariant
derivative on $T^\ast(\RR^n)\otimes u^{-1}(TN)$ defined with respect to the
Levi-Civita connection of the Euclidean metric on $\RR^n$ (i.e. the ordinary
directional derivative) and the Riemannian metric $g$ on $N$. For
$\a=1,\ldots, n$ we let $\nabla_\a u\in\G(u^{-1}(TN))$ be the vector field
given by
\begin{equation}\label{fo-eqn4}
\nabla_\a u=\p_\a u=\frac{\p u^i}{\p x^\a} \frac{\p}{\p u^i}
\end{equation}
where $(u^1,\ldots,u^k)$ are coordinates about $u(x)\in N$. In particular
\[
du=\frac{\p u^i}{\p x^\a} dx^\a \otimes \frac{\p}{\p u^i} = (\nabla_\a
u)^idx^\a\otimes\frac{\p}{\p u^i}.
\]
The second fundamental form of the map $u$ is defined to be the covariant
derivative of $du$, $\nabla du\in\G((T^2\RR^n)\otimes u^{-1}(TN))$. In local
coordinates we have for $i=1,\ldots, k$ and $\a,\beta\in 1,\ldots n$,
\begin{eqnarray}\label{fo-eqn5tilde}
(\nabla du)^i_{\a\beta} & = & \nabla_\a\nabla_\beta u^i \\
& = & \frac{\p^2 u^i}{\p x^\a\p x^\beta} + \G^i_{jk}(u) \frac{\p u^j}{\p
x^\a} \frac{\p u^k}{\p x^\beta}.\nonumber
\end{eqnarray}
Note that here the subscript $\a$ actually denotes covariant differentiation
with respect to the vector field $\nabla_\a u$ as defined in (\ref{fo-eqn4})
and we have $\nabla_\a\nabla_\beta u=\nabla_\beta\nabla_\a u$. The tension
is simply the trace of $\nabla du$ with respect to the Euclidean metric,
$\d=\d_{\a\beta}$
\begin{eqnarray}\label{fo-eqn6}
\tau(u)^i & = & \nabla_\a\nabla_\a u^i \\
& = & \frac{\p^2 u^i}{\p x^\a\p x^\a} + \G^i_{jk}(u) \frac{\p u^j}{\p x^\a}
\frac{\p u^k}{\p x^\a} \nonumber
\end{eqnarray}
from which we recover (\ref{int-eqn1}).

\subsection{The gradient flow for $G(u)$}

For a given vector field $\xi\in\G(u^{-1}(TN))$, we construct a variation of
$u:\RR^n\to N$ with initial velocity $\xi$ as follows. Define the map
\[
U:\RR^n\times\RR\to N
\]
by setting
\[
U(x,s)=\exp_{u(x)}s\xi(x)
\]
where $\exp_{u(x)}:T_{u(x)}N\to N$ denotes the exponential map. Set
$u_s(x)=U(x,s)$ and now let $\nabla$ denote the natural covariant derivative
on $T^\ast(\RR^n\times \RR)\otimes U^{-1}(TN)$. Then
\begin{eqnarray*}
\left.\frac{d}{ds} G(u_s)\right|_{s=0} & = & \left.\frac{1}{2} \frac{d}{ds}
\int_{\RR^n} |\tau(u_s)|^2 dx\right|_{s=0} \\
& = & \left.\frac{1}{2}\int_{\RR^n}\frac{\p}{\p s}
\langle\tau(u_s),\tau(u_s)\rangle dx\right|_{s=0} \\
& = & \left.\int_{\RR^n}\langle\nabla_s\tau(u_s), \tau(u_s)\rangle
dx\right|_{s=0} 
\end{eqnarray*}
where the inner products are taken with respect to $g$ and we have used the
metric compatibility of $\nabla$. Let $R=R(\cdot,\cdot)\cdot$ denote the
Riemann curvature endomorphism of $\nabla$. Using (\ref{fo-eqn6}) and the
definition of $R$ we see that
\begin{eqnarray*}
\nabla_s\tau(u_s) & = & \nabla_s\nabla_\a\nabla_\a u_s \\
& = & \nabla_\a\nabla_s\nabla_\a u_s-R(\nabla_\a u_s, \nabla_su_s)\nabla_\a
u_s \\
& = & \nabla_\a\nabla_\a\nabla_s u_s - R(\nabla_\a u_s, \nabla_s u_s)\nabla_\a
u_s.
\end{eqnarray*}
Therefore
\begin{eqnarray*}
\left.\frac{d}{ds}G(u_s)\right|_{s=0} & = &
\left.\int_{\RR^n}\langle\nabla_\a\nabla_\a\nabla_su_s - R(\nabla_\a u_s,
\nabla_su_s) \nabla_\a u_s, \tau(u_s)\rangle dx\right|_{s=0} \\
& = & \int_{\RR^n}\langle\nabla_\a\nabla_\a\xi,\tau(u)\rangle dx -
\int_{\RR^n}\langle R(\nabla_\a u,\xi)\nabla_\a u,\tau(u)\rangle dx.
\end{eqnarray*}
By the symmetries of the curvature we have
\[
\int_{\RR^n}\langle R(\nabla_\a u,\xi)\nabla_\a u, \tau(u)\rangle dx =
\int_{\RR^n}\langle R(\nabla_\a u,\tau(u))\nabla_\a u,\xi\rangle dx
\]
and provided that $\tau(u)$ and $\nabla_\a \tau(u)$, for $\a=1,\ldots, n$,
are in $L^2$ (and likewise for $v$) we may integrate by parts to obtain
\begin{equation}\label{fo-eqn7}
\left.\frac{d}{ds}G(u_s)\right|_{s=0} = \int_{\RR^n}
\langle\xi,\nabla_\a\nabla_\a\tau(u)\rangle dx - \int_{\RR^n}\langle
R(\nabla_\a u,\tau(u))\nabla_\a u, \xi\rangle dx.
\end{equation}

\begin{proposition}\label{fo-prop1}
The Euler-Lagrange equation for $G$ acting on $H^{s+1}(\RR^n,N)$, for
$s\ge 3$ is
\begin{equation}\label{fo-eqn8}
F(u)\equiv\nabla_\a\nabla_\a\tau(u)-R(\nabla_\a u,\tau(u))\nabla_\a u=0.
\end{equation}
\end{proposition}

The parabolic regularization of (\ref{int-eqn3}) which we now proceed to
study is
\begin{equation}\label{fo-eqn9}
\left\{\begin{array}{rcl}
\frac{\p u}{\p t} & = & -\e F(u)+J(u)\tau(u)+\beta\tau(u) \\
u(0) & = & u_0 \end{array}\right.
\end{equation}

\subsection{The ambient flow equations}

Rather than attempting to study the parabolic equations (\ref{fo-eqn9})
directly we will focus on the induced ``ambient flow equations'' for
$v=w\circ u$ where $w:(N,g)\to\RR^p$ is a fixed isometric embedding. We fix
a $\d>0$, chosen sufficiently small so that on the $\d$-tubular neighborhood
$w(N)_\delta\subset\RR^p$, the nearest point projection map
\[
\Pi:w(N)_\d\to w(N)
\]
is a smooth map (cf.\ \cite{Si} \S2.12.3). For a point $Q\in w(N)_\d$ set
\[
\rho(Q)=Q-\Pi(Q)\in\RR^p
\]
so that $|\rho(Q)|=\dist(Q,w(N))$, and viewing $\rho$ and $\Pi$ as maps from
$w(N)_\d$ into itself we have 
\begin{equation}\label{fo-eqn10}
\left.\Pi+\rho=\id\right|_{w(N)_\d}.
\end{equation}
Note that then the differentials of the maps satisfy
\begin{equation}\label{fo-eqn11}
d\Pi+d\rho=\id
\end{equation}
as a linear map from $\RR^p$ to itself. For any map $v:\RR^n\to w(N)_\d$ we
set
\[
T(v)=\Delta v-\Pi_{ab}(v) v^a_\a v^b_\a
\]
where $\Pi_{ab}(v)$, $1\le a,b\le p$ are the components of the Hessian of
$\Pi$ at $v(\cdot)$. At a point $y\in N$ the Hessian of $\Pi$ is minus the
second fundamental form of $N$ at $y$. So if $v=w\circ u$, with $u:\RR\to
N$, then $T(v)$ is simply the tangential component of the Laplacian of $v$
which corresponds to the tension of the map $u$, i.e.
\[
dw(\tau(u))=(\D v)^T=d\Pi(\D v)=T(v).
\]
Therefore, in direct analogy with the functional $G(\cdot)$, we now consider
\begin{eqnarray*}
\cG(v) & = & \frac{1}{2}\int_{\RR^n}|T(v)|^2dx \\
& = & \frac{1}{2}\int_{\RR^n} |\D v-\Pi_{ab}(v) v^a_\a v^b_\a|^2 dx
\end{eqnarray*}
Our point here (and hence the seemingly odd notation) is that we wish to
consider $T(v)$ for arbitrary maps into $w(N)_\delta$ whose image does not
necessarily lie on $N$.

\begin{defn}\label{fo-defn2}
For $v:\RR^n\to w(N)_\delta$, let $\cF(v)$ denote the Euler-Lagrange operator
of $\cG(v)$ with respect to unconstrained variations. A simple computation
shows that its components are given by
\begin{eqnarray*}
(\cF(v))^c&=&(\D T(v))^c\\
&&-\sum^n_{\a,\beta=1}\left(T(v)^e\Pi^e_{abc}(v) v^a_\a
v^b_\beta - (T(v)^e\Pi^e_{ac}(v) v^a_\a)_\beta - (T(v)^e\Pi^e_{cb}(v)
v^b_\beta)_\a\right)\\
&=& \Delta^2 v^c-(\tilde{\cF}(v))^c,
\end{eqnarray*}
where 
\begin{align*}
(\tilde{\cF}(v))^c=&\sum^n_{\alpha,\beta=1} \big(\Delta (\Pi^c_{ab}(v)v_\alpha^a v_\beta^b)+T(v)^e
\Pi^e_{abc}(v) v^a_\alpha v^b_\beta - (T(v)^e\Pi^e_{ac}(v) v^a_\alpha)_\beta \\
&- (T(v)^e\Pi^e_{cb}(v) v^b_\beta)_\alpha\big)
\end{align*}
denotes the lower order terms. Note that the subscripts here refer to coordinate differentiation in $\RR^n$ (Greek indices) or $\RR^p$ (Roman indices).
\end{defn}

For $v=w\circ u$, we wish to consider compactly supported tangential
variations of $\cG(v)$. Such variations correspond to (compactly supported)
vector fields $\phi$ on $w(N)_\delta$ which satisfy $d\rho(\phi)=0$.

\begin{proposition}\label{fo-prop2}
If $u:\RR^n\to N$ and $v=w\circ u$ then
for all $\phi\in\G(Tw(N)_\d)$ with compact support such that $d\rho(\phi)=0$
we have
\[
\left.\frac{d}{ds}\cG(v+s\phi)\right|_{s=0}=\int_{\RR^n} \langle
\cF(v),\phi\rangle dx =\int_{\RR^n} \langle
d\Pi(\cF(v)),\phi\rangle dx.
\]
Recall that $d\Pi=\rm{Id} -d\rho$.
\end{proposition}

\begin{defn}\label{fo-defn3}
If $v=w\circ u$, then the ambient form of the \sh vector field
$J(u)\tau(u)$, is given by the vector field $f_v$ with
\begin{equation}\label{fo-eqn12}
f_v=dw_{|_{w^{-1}\Pi(v(x))}}[J(w^{-1}\Pi(v))(dw)^{-1}_{|_{\Pi(v(x))}}
(d\Pi_{|_{v(x)}}(\D v))].
\end{equation}
Note that $f_v$ is defined for maps $v:\RR^n\to w(N)_\d$ whose image does
not necessarily lie on $N$.
\end{defn}
Next we have the following Lemma.
\begin{lemma}\label{normal}
If $u:\RR^n \rightarrow N$ and $v=w\circ u$ then we have
\begin{align*}
d\rho(v)(\cF(v))=& \Delta(\Pi_{ab}(v) \nabla_\alpha v^a \nabla_\alpha v^b)+\text{div}(\Pi_{ab}(v) \Delta v^a \nabla v^b)+\Pi_{ab}(v) \nabla_\alpha \Delta v^a \nabla_\alpha v^b\\
&-d\rho(v)(\tilde{\cF}(v))\\
=:& \cH(v).
\end{align*}
\end{lemma}
\begin{proof}
This follows from the facts that
\begin{align*}
d\rho(v)(\cF(v))=d\rho(v)(\Delta^2 v)-d\rho(v)(\tilde{\cF}(v)),
\end{align*}
\begin{align*}
\Delta \text{div}\big(d\rho(v)(\nabla v)\big)=& \Delta \big(d\rho(v)(\Delta v)\big) +\Delta\big(\rho_{ab}(v) \nabla v^a \nabla v^b\big)\\
=&\text{div}\big(d\rho(v)(\nabla \Delta v)\big)+\text{div}\big(\rho_{ab}(v) \Delta v^a \nabla v^b\big)\\
&+\Delta\big(\rho_{ab}(v) \nabla v^a \nabla v^b\big)\\
=& d\rho(v)(\Delta^2 v)+\rho_{ab} (v)\nabla \Delta v^a \nabla v^b\\
&+\text{div}\big(\rho_{ab}(v) \Delta v^a \nabla v^b\big)+\Delta\big(\rho_{ab}(v) \nabla v^a \nabla v^b\big)
\end{align*}
and
\begin{align*}
d\rho(v)(\nabla v)=0.
\end{align*}
Finally we note that
\begin{align*}
\rho_{ab}(v)=-\Pi_{ab}(v).
\end{align*}
\end{proof}
\begin{rem}
Note that $\cH(v)$ only contains derivatives of $v$ up to third order. Moreover this term is well-defined for every $v:\RR^n \rightarrow w(N)_\delta$.
\end{rem}
The regularized ambient equations are given by
\begin{equation}\label{fo-eqn13}
\left\{\begin{array}{rcl}
\frac{\p v}{\p t} & = & -\e(\cF(v)-\cH(v))+f_v+\beta Tv \\
v(0) & = & v_0 \end{array}\right.
\end{equation}
The basic relationship between the regularized geometric flows
(\ref{fo-eqn9}) and the regularized ambient flows (\ref{fo-eqn13}) is
provided by the following Lemma (cf.\ \S7 of \cite{ES}).

\begin{lemma}\label{fo-lem1}
Fix $\e\in[0,1]$. Given $u_0\in H^{s+1}(\RR^n,N)$ with
$s\ge3$, $w:N\to\RR^p$ an isometric embedding, and
$T_\e>0$, a flow $u:\RR^n\times[0,T_\e]\to N$ satisfies (\ref{fo-eqn9}) if
and only if the flow $v=w\circ u:\RR^n\times[0, T_\e]\to\RR^p$ satisfies
(\ref{fo-eqn13}) with $v_0=w\circ u_0$.
\end{lemma}

\begin{proof}
First note that since $w$ is an isometry we have
\begin{equation}\label{fo-eqn14}
|T(v)|^2=|\tau(u)|^2
\end{equation}
and therefore $\cG(v)=G(u)$. Given $\xi\in \G(u^{-1}(TN))$ a smooth compactly
supported vector field set $\phi=dw(\xi)\in \G(u^{-1}(T\RR^p))$. As before
we consider the variation of $u$ given by $u_s(x)=\exp_{u(x)}s\xi$. We then
have
\[
w\circ u_s=v+s\phi+\cO(s^2)
\]
so that
\[
G(u_s)=\cG(v+s\phi)+\cO(s^2).
\]
Therefore
\[
\int_{\RR^n}\langle F(u),\xi\rangle dx =
\int_{\RR^n}\langle\cF(v),\phi\rangle dx.
\]
Observe that
\begin{equation}\label{fo-eqn15}
\int_{\RR^n}\left\langle\frac{\p u}{\p t},\xi\right\rangle dx = \int_{\RR^n}
\left\langle dw\left(\frac{\p u}{\p t}\right), dw(\xi)\right\rangle dx =
\int_{\RR^n} \left\langle\frac{\p v}{\p t},\phi\right\rangle dx.
\end{equation}
Since $d\rho(\phi)=0$ and $\cH(v)=d\rho(\cF(v))$ we also have
\begin{equation}\label{fo-eqn16}
-\e\int_{\RR^n}\langle F(u),\xi\rangle dx =
-\e\int_{\RR^n}\langle\cF(v),\phi\rangle
dx=-\e\int_{\RR^n}\langle\cF(v)-\cH(v),\phi\rangle dx.
\end{equation}
Note that
\begin{eqnarray*}
\int_{\RR^n}\langle J(u)\tau(u),\xi\rangle dx & = & \int_{\RR^n}\langle
dw(J(u)\tau(u)),dw(\xi)\rangle dx \\
& = & \int_{\RR^n}\langle dw[J(w^{-1}(\Pi(v)))(dw)^{-1}(T(v))],
dw(\xi)\rangle dx \\
& = & \int_{\RR^n}\langle f_v,\phi\rangle dx
\end{eqnarray*}
and 
\begin{eqnarray*}
\int_{\RR^n}\langle \tau(u),\xi\rangle dx & = & \int_{\RR^n}\langle
dw(\tau(u)),dw(\xi)\rangle dx \\
& = & \int_{\RR^n}\langle Tv,\phi\rangle dx.
\end{eqnarray*}
This together with
(\ref{fo-eqn15}) and (\ref{fo-eqn16}) implies that the flows correspond as
claimed.
\end{proof}

We end this section by exhibiting in a more practical form the structure of
the parabolic operator appearing in the regularized ambient flow equations
(\ref{fo-eqn13}).

\begin{defn}\label{fo-defn4}
For $v:\RR^n\to\RR^p$, and $j\in\NN$ we let $\p^jv$ denote an arbitrary
$j^th$-order partial derivative of $v$
\[
\p^jv = \frac{\p^jv}{\p x^{\a_1}\cdots\p x^{\a_r}}\qquad\mbox{with}\qquad
\a_1+\cdots\a_r=j
\]
and let
\[
\p^{j_1}v\ast\cdots\ast\p^{j_l}v
\]
denote terms which are a sum of products of terms of the form $\p^{j_1}v,
\ldots, \p^{j_l}v$.
\end{defn}

\begin{proposition}\label{fo-prop3}
Let $v:\RR^n\to w(N)_\delta\subset\RR^p$, then
\begin{gather*}
-\e(\cF(v)-\cH(v)) + f_v + \beta Tv \\
 = -\e\D^2v-\e\sum^4_{l=2}\sum_{j_1+\cdots+j_l=4} A_{(j_1\cdots
j_l)}(v)\p^{j_1}v\ast\cdots\ast\p^{j_l}v + B_0(v)\p^2v+B_1(v)\p v\ast\p v
\end{gather*}
where each $j_s\ge 1$ and each of $A_{(j_1\cdots j_l)}(v)$, $B_0(v)$ and
$B_1(v)$ are bounded smooth functions of $v$.
\end{proposition}

\begin{proof}
This follows from the explicit expressions for $\cF(v)$, $\cH(v)$, $f_v$ and $Tv$.
\end{proof}

In the following Lemma (which is a suitable modification of Theorem $7$C of \cite{ES}) we show that if $v:\RR^n \times [0,T]\rightarrow w(N)_\delta$ is a solution of \eqref{fo-eqn13} and if $v_0:\RR^n \rightarrow w(N)$, then $v: \RR^n \times [0,T]\rightarrow w(N)$.

\begin{lemma}\label{fo-lem2.3A}
Fix $\e \ge0$ and $\beta\ge 0$. Let $v:\RR^n\times [t_0,t_1]\to w(N)_\delta$ be a solution of \eqref{fo-eqn13} with $v(x,t_0)  =  v_0(x) \in w(N)$, where $v_0\in H^{s+1}(\RR^n, w(N))$ with $s\ge [\frac{n}{2}] +4$. Then 
$v(x,t)\in w(N)$ for all $x\in\RR^n$ and all $t\in [t_0,t_1]$.
\end{lemma}

Note that in this case by Lemma \ref{fo-lem1} 
$u(x,t)=w^{-1}\circ v(x,t)$ solves
\begin{equation}
\left\{\begin{array}{rcl}
\frac{\p u}{\p t} & = & -\e F(u)+J(u)\tau(u)+\beta\tau(u)\nonumber \\
u(x,t_0) & = & u_0(x)=w^{-1}\circ v_0(x). \end{array}\right.\nonumber
\end{equation}

\begin{proof}
We start by calculating
\begin{align*}
\partial_t \rho(v)=&d_a \rho(v) \partial_t v^a,\\
\Delta \rho(v)=& d_a \rho(v) \Delta v^a- \Pi_{ab}(v)\nabla v^a \nabla v^b\ \ \ \text{and} \\
\Delta^2 \rho(v)=& \Delta \Big(d_a \rho(v) \Delta v^a\Big)-\Delta\Big( \Pi_{ab}(v)\nabla v^a \nabla v^b\Big)\\
=& \text{div}\Big(d_a \rho(v) \nabla \Delta v^a\Big)-\text{div}\Big( \Pi_{ab}(v)\Delta v^a \nabla v^b\Big)-\Delta\Big( \Pi_{ab}(v)\nabla v^a \nabla v^b\Big)\\
=& d_a \rho(v) \Delta^2 v^a- \Pi_{ab}(v)\nabla\Delta v^a \nabla v^b-\text{div}\Big( \Pi_{ab}(v)\Delta v^a \nabla v^b\Big)\\
&-\Delta\Big( \Pi_{ab}(v)\nabla v^a \nabla v^b\Big).
\end{align*}
Here we used again that $\Pi_{ab}(v)=-\rho_{ab}(v)$. Now if $v$ is a solution of \eqref{fo-eqn13} we get that
\begin{align*}
\partial_t v=&-\e \Delta^2 v+ \e \tilde{\cF}(v)-\e d\rho(v)\tilde{\cF}(v)+f_v+\beta \Delta v-\beta \Pi_{ab}(v)\nabla v^a \nabla v^b\\
&+\e\Big(\Delta(\Pi_{ab}(v) \nabla v^a \nabla v^b)+\text{div}(\Pi_{ab}(v) \Delta v^a \nabla v^b)+\Pi_{ab}(v) \nabla \Delta v^a \nabla v^b\Big)\\
=&\beta \Delta v-\e \Delta^2 v+\e d\Pi(v)\tilde{\cF}(v)+f_v-\beta \Pi_{ab}(v)\nabla v^a \nabla v^b\\
&+\e\Big(\Delta(\Pi_{ab}(v) \nabla v^a \nabla v^b)+\text{div}(\Pi_{ab} (v)\Delta v^a \nabla v^b)+\Pi_{ab}(v) \nabla \Delta v^a \nabla v^b\Big).
\end{align*}
Combining these two calculations yields
\begin{align*}
(\partial_t-&\beta \Delta+\e \Delta^2)\rho(v)= \e d\rho(v)\big( d\Pi(v) \tilde{\cF}(v)\big)+d\rho(v)f_v-\beta d\rho(v)\big(\Pi_{ab}(v)\nabla v^a \nabla v^b\big)\\
&+\e d\rho(v)\Big(\Delta(\Pi_{ab}(v) \nabla v^a \nabla v^b)+\text{div}(\Pi_{ab}(v) \Delta v^a \nabla v^b)+\Pi_{ab}(v) \nabla \Delta v^a \nabla v^b\Big)\\
&+\beta\Pi_{ab}(v)\nabla v^a \nabla v^b-\e \Big(\Delta(\Pi_{ab}(v) \nabla v^a \nabla v^b)+\text{div}(\Pi_{ab}(v) \Delta v^a \nabla v^b)\\
&+\Pi_{ab}(v) \nabla \Delta v^a \nabla v^b\Big)\\
=& d\rho(v)f_v+\beta d\Pi(v)\big(\Pi_{ab}(v)\nabla v^a \nabla v^b\big)  -\e d\Pi(v)\Big(\Delta(\Pi_{ab}(v) \nabla v^a \nabla v^b)\\
&+\text{div}(\Pi_{ab}(v) \Delta v^a \nabla v^b)+\Pi_{ab}(v) \nabla \Delta v^a \nabla v^b\Big)+\e d\rho(v)\big( d\Pi(v) \tilde{\cF}(v)\big).
\end{align*}
Multiplying this equation with $\rho(v)$ and using the facts that (note that $f_v \in T_{\Pi(v)}w(N)$)
\begin{align*}
\rho(v)\cdot d\Pi(v)(\phi)=&0\ \ \ \forall \phi\in w(N)_\delta,\\
\rho(v)\cdot d\rho(v)f_v=\rho(v)\cdot(f_v-d\Pi(v) f_v)=&0,\ \ \ \ \text{and}\\
\rho(v)\cdot d\rho(v)\big(d\Pi(v)(\phi)\big)=\rho(v)\cdot d\Pi(v)\big(d\rho(v)(\phi)\big)=&0 \ \ \ \forall \phi\in w(N)_\delta,
\end{align*}
gives for all $t\in (t_0,t_1)$
\begin{align*}
\frac12 \partial_t |\rho(v)|^2=& \langle \rho(v), \beta \Delta \rho(v)-\e \Delta^2 \rho(v) \rangle.
\end{align*}
Integrating this equation over $\RR^n$ and using integration by parts, we have for all $t\in (t_0,t_1)$
\begin{align*}
 \partial_t \int_{\RR^n} |\rho(v)|^2 =& -2\int_{\RR^n} (\beta |\nabla \rho(v)|^2+\e|\Delta \rho(v)|^2)\\
\le& 0.
\end{align*}
Since $\rho(v_0)=0$ this implies that $\rho(v)=0$ for all $t\in [t_0,t_1]$ and hence finishes the proof of the Lemma.
\end{proof}

\section{The Duhamel solution to the ambient flow equations}

In this section we introduce a fixed point method that solves the initial
value problem (\ref{fo-eqn13}) in the Sobolev space $H^{s+1}(\RR^n, \RR^p)$, for
$s\geq\frac{n}{2}+4$. To simplify the notation, using Proposition
\ref{fo-prop3}, we rewrite (\ref{fo-eqn13}) as
\begin{equation}\label{du-eqn17}
\left\{\begin{array}{rcl}
\frac{\p v}{\p t} & = & -\e\D^2v+N(v) \\
v(x,0) & = & v_0,
\end{array}
\right.
\end{equation}
where
\begin{eqnarray}\label{du-eqn18}
N(v) & = & -\e\sum^4_{l=2}\sum_{j_1+\cdots+j_l=4} A_{(j_1,\ldots,
j_l)}(v)\p^{j_1}v\ast\cdots\ast\p^{j_l}v \\
&&\qquad +B_0(v)\p^2 v+B_1(v)\p v\ast\p v.\nonumber
\end{eqnarray}
We now state the well-posedness theorem for (\ref{du-eqn17}). 
For any fixed $v_0$, define the
spaces
\[
L^{2}_\d=\{v:\RR^n\rightarrow \RR^p |\,\   \|v-v_0\|_{L^2}<\d\}.
\]
and
\[
L^{2,\infty}_\d=\{v:\RR^n\rightarrow \RR^p |\,\   \|v-v_0\|_{L^2}, \|v-v_0\|_{L^\infty}<\d\}.
\]
We then have the following theorem:
\begin{theorem}\label{du-thm2}
Assume $\d>0$, $\e>0$, and $\gamma \in\RR^p$ are fixed. Then for any 
$(v_0-\gamma)\in H^{s+1}(\RR^n, \RR^p)$, $ s\geq \frac{n}{2}+4$ there
exist $T_\e=T(\e,\d,\|\p v_0\|_{H^s}, \|v_0-\gamma\|_{L^2})$ and a unique solution $v=v_\e$ for
(\ref{du-eqn17}) such that $v\in C([0,T_\e],H^{s+1}\cap L^{2,\infty}_\d)$. 
\end{theorem}

To prove the theorem we rewrite (\ref{du-eqn17}) as an integral equation
using the Duhamel principle:
\begin{equation}\label{du-eqn19}
v(x,t)=S_\e(t)(v_0-\gamma)+\int^t_0S_\e(t-t')N(v)(x,t')dt'+\gamma,
\end{equation}
where for $f\in H^{s+1}(\RR^n, \RR^p)$
\begin{equation}\label{du-eqn20}
S_\e(t)f(x) = \int_{\RR^n} e^{(i\langle x,\xi\rangle-\e|\xi|^4t)}
\widehat{f}(\xi)d\xi
\end{equation}
is the solution of the linear and homogeneous initial value problem
associated to (\ref{du-eqn17}). The main idea is to consider the operator
\begin{equation}\label{du-eqn21}
Lv(x,t)=S_\e(t)(v_0-\gamma)+\int^t_0S_\e(t-t')N(v)(x,t')dt' +\gamma
\end{equation}
and prove that for a certain $T_\e$ the operator $L$ is a contraction from a suitable ball in $C([0,T_\e],H^s\cap L^{2,\infty}_\d)$ into itself.

To estimate $L$ we need to study the smoothing properties of the linear
solution $S_\e(t)v_0$. Because the order of derivatives that appears in
$N(v)$ is 3, in order to be able to estimate the nonlinear part of $L$ in
$H^{s+1}$, we should prove that the operator $S_\e(t)$ provides a smoothing
effect also of order 3. We have in fact the following lemma:

\begin{lemma}\label{du-lem3}
Define the operator $D^s,$ $s\in\RR$ as the multiplier operator such that
$\widehat{D^sf}(\xi)=|\xi|^s\widehat f$. Then for any $t>0$ and $i=1,2,3$, 
\begin{equation}\label{du-eqn22}
\|S_\e(t)f\|_{L^2} \lesssim\|f\|_{L^2},
\end{equation}
\begin{equation}\label{du-eqn23}
\|D^sS_\e(t)f\|_{L^2} \lesssim
t^{-\frac{i}{4}}\e^{-\frac{i}{4}}\|D^{s-i}f\|_{L^2}.
\end{equation}
\end{lemma}

\begin{proof}
The proof follows from Plancherel theorem and the two estimates
\begin{eqnarray*}
\biggl|e^{-\e|\xi|^4t}\biggr| & \lesssim & 1 \\
|\xi|^s\biggl|e^{-\e|\xi|^4t}\biggr| & \lesssim &
|\xi|^{s-i}t^{-\frac{i}{4}}\e^{-\frac{i}{4}}.
\end{eqnarray*}
\end{proof}

\noindent
To state the next lemma, where we show how for small intervals of time the evolution 
$S_\e(t)(v_0-\gamma)$ stays close to $v_0-\gamma$, we need to
introduce the space $\dot H^s$. This space  denotes the homogeneous Sobolev space defined 
as the set of all functions $f$ such that $D^s f\in L^2$.

\begin{lemma}\label{du-lem4}
Let $\s\in (0,1)$, $ s>\frac{n}{2}+4\s$ and assume that
$f\in H^{4\s}\cap \dot H^{s}$. Then
\begin{equation}\label{du-eqn24}
\|S_\e(t) f-f\|_{L^\infty} \le \e^\s t^\s[\|f\|_{\dot H^s}
+\|f\|_{\dot H^{4\s}}],\ \hbox{and}\ 
\|S_\e(t) f-f\|_{L^2} \le \e^\s t^\s\|f\|_{\dot H^{4\s}}.
\end{equation}
\end{lemma}

\begin{proof}
By the mean value theorem
\[
\biggl|e^{-\e|\xi|^4t}-1\biggr|\lesssim |\xi|^4t\e,
\]
which combined with the trivial bound
\[
\biggl|e^{-\e|\xi|^4t}-1\biggr|\le 2
\]
gives, for any $\s\in[0,1]$
\begin{equation}\label{du-eqn25}
\biggl|e^{-\e|\xi|^4t}-1\biggr|\lesssim(|\xi|^4t\e)^\s.
\end{equation}
We now write
\begin{eqnarray}\label{tt1}
\ \ \ \ \ |(S_\e(t)-1)f(x)| &=& \biggl|\int_{\RR^n} e^{i\langle
x,\xi\rangle}[e^{-\e|\xi|^4t}-1]\widehat{f}(\xi)d\xi\biggr| \\
&\lesssim& (t\e)^\s\int_{\RR^n} |\widehat{f}|(\xi)|\xi|^{4\s}\nonumber\\
&=& (t\e)^\s\left[\int_{|\xi|\leq 1} |\widehat{f}|(\xi)|\xi|^{4\s}
+\int_{|\xi|\geq 1} |\widehat{f}|(\xi)
|\xi|^s\frac{1}{|\xi|^{s-4\s}}\right],\nonumber
\end{eqnarray}
and this concludes the argument since $s>  \frac{n}{2} +4\s$
guarantees the summability after the application of Cauchy-Schwartz.
Note also that 
\begin{eqnarray}\label{tt2}
\|(S_\e(t)-1)f\|_{L^2}& \le& \| (e^{-|\xi|^4\e t}-1)\widehat f\|_{L^2}\\
&\lesssim & (t\e)^\sigma 
\left(\int( |\xi|^{4\s})^{2}|\widehat f|^2\right)^{\frac{1}{2}}\nonumber\\
&\lesssim & (t\e)^\sigma \|f\|_{\dot H^{4\s}}
\lesssim (t\e)^\sigma \|\p f\|_{H^{4\s-1}}
.\nonumber
\end{eqnarray}

\end{proof}

\noindent
We are now ready to prove Theorem \ref{du-thm2}.

\begin{proof}
For $T_\e,r>0$ and $s\geq \frac{n}{2}+4$ consider the ball
\[
B_r=\{ \p v\in H^s: \|\p(v-v_0)\|_{L^\infty_{T_\e},H^s} \le r\}\cap L^{2,\infty}_\d.
\]
We want to prove that  for the appropriate 
$T_\e$ and $r$, the operator $L$ maps
$B_r$ to itself and is a contraction. We start with the estimate of the
linear part of $L$. By (\ref{du-eqn22}) we have
\begin{equation}\label{du-eqn26}
\|\p(S_\e(t)(v_0-\gamma)-(v_0-\gamma))\|_{H^s} \lesssim \|(1+D^s)S_\e(t)\p v_0\|_{L^2} +
\|\p v_0\|_{H^s}\lesssim \|\p v_0\|_{H^s}.
\end{equation}
To estimate the nonlinear term we use (\ref{du-eqn22}) 
(\ref{du-eqn23}), and interpolation :
\begin{eqnarray}\label{du-eqn27}
&&\biggl\| \partial \biggl(\int^t_0S_\e(t-t')N(v)(x,t')dt'\biggr)\biggr\|_{H^s}\\
&& = \left\|\int^t_0S_\e(t-t')\p N(v)(x,t')dt'\right\|_{L^2} + 
\left\|\int^t_0D^sS_\e(t-t')\p N(v)(x,t')dt'\right\|_{L^2}\nonumber \\
&&\lesssim  \int_0^t\|\p N(v)\|_{L_x^2}(t')\, dt'+\int_0^t(t')^{-3/4}\e^{-3/4}\|D^{s-3}\p N(v)\|_{L_x^2}(t')\, dt'\nonumber\\
&&\lesssim \int_0^t (1+(t')^{-3/4}\e^{-3/4})\|\p v\|^m_{H_x^s}(t')\,
dt',\nonumber
\end{eqnarray}
where $m$ is the order of the nonlinearity\footnote{In our case
actually one can compute that $m=4$.} $N(v)$.
Note that to control $\p N(v)$ and $D^{s-3}\p N(v)$ in the previous inequality
we are never in the position of estimating $v$ in $L^2$.
By (\ref{du-eqn21}), (\ref{du-eqn23}), (\ref{du-eqn26}) and (\ref{du-eqn27}), we obtain the estimate
\begin{equation}\label{du-eqn28}
\|\p(Lv-v_0)\|_{H_x^s}(t)\le C_0\|\p v_0\|_{H^s} +
C_1\int_0^t (1+(t')^{-3/4}\e^{-3/4})\|\p v\|^m_{H_x^s}(t')\, dt'.
\end{equation}
Thus 
\begin{equation}\label{tt3}
\|\p(Lv- v_0)\|_{L^{\infty}_{T_\e}H_x^s}\le C_0\|\p v_0\|_{H^s} +
C_1\e^{-3/4}T_{\e}^{1/4}\|\p v\|^m_{L^{\infty}_{T_\e}H_x^s}.
\end{equation}

We still need to check that $Lv$ is continuous in time and that $Lv\in
L^{\infty,2}_\d$. The continuity follows directly from the continuity of the
operator $S_\e(t)$. To prove the $L^\infty$ and $L^2$ estimates one uses
(\ref{du-eqn24}) with $\sigma=1/4$ applied to $f=v_0-\gamma$, the 
Sobolev inequality and estimates similar to the ones 
used to obtain (\ref{du-eqn28}). One gets 
\begin{eqnarray}\label{du-eqn29}
\|Lv-v_0\|_{L^\infty_{T_\e} L^2_x} & + & \|Lv-v_0\|_{L^\infty_{T_\e}
L^\infty_x}\le C_1\e^{1/4} T^{1/4}_\e\|\p v_0\|_{H^s} \\
& + &
C_1\e^{-\frac{3}{4}}T^{1/4}_\e\|\p v\|^m_{L^{\infty}_{T_\e}H_x^s}.\nonumber
\end{eqnarray}
We now take $r=3C_0\|v_0\|_{H^s}$ and
\begin{equation}\label{du-eqn30}
T_\e\le\min(\d^4C_1^{-4}\e^3 6^{-4m} C^{-4m}_0\|v_0\|^{-4m+4}_{H^s},
\d^{4} 2^{-4}
C^{-4}_1\e^{-1}\|v_0\|^{-4}_{H^s})
\end{equation}
so that (\ref{du-eqn28}) and (\ref{du-eqn29}) guarantees that $L$ maps
$B_r$ into itself. Note that (\ref{du-eqn27}) yields for $v,w\in B_r$
\begin{eqnarray}\label{tt4}
&& \kern-.7in\biggl\|\int_0^t S_\e(t-t')\p[N(v)-N(w)](x,t')\, 
\ dt'\biggr\|_{H^s} \\
& \lesssim &
\e^{-3/4}T_\e^{1/4}\|\p[N(v)-N(w)]\|_{L^\infty_{T_\e} H_x^{s-3}}.\nonumber
\end{eqnarray}
Therefore 
\begin{eqnarray}\label{tt5}
&&\|\p(Lv-Lw)\|_{L^\infty_{T_\e} H_x^{s}} \\
&& \qquad\lesssim \e^{-3/4}T_\e^{1/4}\|\p[N(v)-N(w)]\|_{L^\infty_{T_\e} H_x^{s-3}}
\nonumber\\
&&\qquad \lesssim 
\e^{-3/4}T_\e^{1/4}C(\d)(\|\p v\|^{m-1}_{L^\infty_{T_\e} H_x^{s}} +\|\p
w\|^{m-1}_{L^\infty_{T_\e} H_x^{s}})
\|\p(v-w)\|_{L^\infty_{T_\e} H_x^{s}}\nonumber\\
&&\qquad \lesssim
\e^{-3/4}T_\e^{1/4}C(\d,\|\p v_0\|_{H_x^{s}})
\|\p(v-w)\|_{L^\infty_{T_\e} H_x^{s}}.\nonumber
\end{eqnarray}
Similarly one shows that 
\begin{equation}\label{tt5A}
\|Lv-Lw\|_{L^\infty_{T_\e} H_x^{s}} 
\lesssim
\e^{-3/4}T_\e^{1/4}C(\d,\|\p v_0\|_{H_x^{s}})
\|\p(v-w)\|_{L^\infty_{T_\e} H_x^{s}}.\nonumber
\end{equation}
By shrinking $T_\e$ further by an absolute constant if necessary,  
from (\ref{tt5}) and (\ref{tt5A}) we obtain
\begin{equation}\label{tt6}
\|Lv-Lw\|_{L^\infty_{T_\e} H_x^{s+1}} \le
\frac{1}{2}
\|v-w\|_{L^\infty_{T_\e} H_x^{s+1}}.
\end{equation}
The contraction mapping theorem ensures that there exists a unique 
function $v=v_\e$ in $L^2_\delta\cap 
\{ \p v\in H^s: \|\p(v-v_0)\|_{L^\infty_{T_\e},H^s} \le r\}$
which solves the integral equation (\ref{du-eqn19}) in
the time interval [0, $T_\e$] defined in (\ref{du-eqn30}).  
Moreover $v\in B_r$ by our choice of $T_\epsilon$.
The uniqueness in the whole space 
$H_x^{s+1}\cap L^{2,\infty}_\d$ follows by similar 
and by now classical arguments.
\end{proof}

\section{Analytic preliminaries}

In this section we state and present the detailed proof of an interpolation
inequality for Sobolev sections on vector bundles which appears in
\cite{DW2001} (see Theorem 2.1). This inequality was first proved for
functions on $\RR^n$ by Gagliardo and Nirenberg, and for functions on
Riemannian manifolds by Aubin \cite{A}. The justification for presenting a
complete proof is that this estimate plays a crucial role in the energy
estimates and therefore in the proof of the results this paper. The precise
dependence of the constants involved in this inequality is vital to our
argument and we feel compelled to emphasize it.

Let $\Pi:E\to\RR^n$ be a Riemannian vector bundle over $\RR^n$. We have the
bundle $\L^PT^\ast\RR^n\otimes E\to\RR^n$ over $\RR^n$ which is a tensor
product of the bundle $E$ and the induced $P$-form bundle over $\RR^n$, with
$P=1,2,\ldots, n$. We define $T(\L^PT^\ast\RR^n\otimes E)$ as the set of all
smooth sections of $\L^PT^\ast\RR^n\otimes E\to\RR^n$. There exists an
induced metric on $\L^PT^\ast\RR^n\otimes E\to\RR^n$ from the metric on
$T^\ast\RR^n$ and $E$ such that for any $s_1,
s_2\in\G(\L^PT^\ast\RR^n\otimes E)$
\begin{equation}\label{an-eqn1}
\langle s_1,s_2\rangle = \sum_{i_1\le\cdots\le i_p}\langle
s_1(e_{i_1},\ldots e_{i_p}),s_2(e_1,\ldots e_{i_p})\rangle
\end{equation}
where $\{e_i\}$ is an orthonormal local frame for $T\RR^n$. We define the
inner product on $\G(\Lambda^PT^\ast\RR^n\otimes E)$ as follows
\begin{equation}\label{an-eqn2}
(s_1,s_2) = \int_{\RR^n}\langle s_1,s_2\rangle(x)dx.
\end{equation}
The Sobolev space $L^2(\RR^n,\Lambda^PT^\ast\RR^n\otimes E)$ is the
completion of $\G(\L^PT^\ast\RR^n\otimes E)$ with respect to the above inner
product. To define the bundle-valued Sobolev space
$H^{k,r}(\RR^n,\L^PT^\ast\RR^n\otimes E)$ consider $\nabla$ the covariant
derivative induced by the metric on $E$, then take the completion of smooth
sections of $E$ in the norm
\begin{equation}\label{an-eqn3}
\|s\|_{H^{k,r}}=\|s\|_{k,r}=\left(\sum^k_{i=0}\int_{\RR^n}|\nabla^is|^rdx\right)^{\frac{1}{r}}
\end{equation}
where 
\begin{equation}\label{an-eqn4}
|\nabla^is|^2=\langle\underbrace{\nabla\cdots\nabla}_{i-\mbox{\footnotesize{times}}}s,
\underbrace{\nabla\cdots\nabla}_{i-\mbox{\footnotesize{times}}}s\rangle.
\end{equation}
If $r=2$, $H^{k,r}=H^k$.

\begin{proposition}\label{an-prop1}
Let $s\in C^\infty_c(E)$ where $E$ is a finite dimensional $C^\infty$ vector
bundle over $\RR^n$. Then given $q, r\in[1,\infty]$ and integers $0\le j\le
k$ we have that
\begin{equation}\label{an-eqn5}
\|\nabla^js\|_{L^p}\le C\|\nabla^ks\|^a_{L^q}\|s\|^{1-a}_{L^r}
\end{equation}
with $p\in[2,\infty)$, $a\in \left(\frac{j}{k}, 1\right]$ and satisfying
\begin{equation}\label{an-eqn6}
\frac{1}{p}=\frac{j}{n}+\frac{1}{r}+a\left(\frac{1}{q}-\frac{1}{r}-\frac{k}{n}\right).
\end{equation}
If $r=n/k-1\ne 1$ then (\ref{an-eqn5}) does not hold for $a=1$. The constant 
$C$ that appears in (\ref{an-eqn5}) only depends on $n, k, j, q, r$ and $a$.
\end{proposition}

\begin{proof}
If $f$ is a real valued smooth function with compact support on $E$ then
 Theorem 3.70 in \cite{Au} ensures that
(\ref{an-eqn5}) holds. 
\medskip

\noindent{\bf Case 1}: Let $j=0$ and $k=1$. Then for $f=|s|$ we have by
(\ref{an-eqn5}) that 
\begin{equation}\label{an-eqn7}
\|s\|_{L^p}\le C\|\nabla|s|\|^a_{L^q}\|s\|^{1-a}_{L^r}.
\end{equation}
Kato's inequality ensures that $|\nabla|s||\le |\nabla s|$ which using
(\ref{an-eqn7}) yields
\begin{equation}\label{an-eqn8}
\|s\|_{L^p}\le C\|\nabla s\|^a_{L^q}\|s\|^{1-a}_{L^r},
\end{equation}
which proves (\ref{an-eqn5}) for $j=0$ and $k=1$. In general if
$f=|\nabla^js|$ Kato's inequality ensures that $|\nabla|\nabla^js||\le
|\nabla^{j+1}s|$ which yields using (\ref{an-eqn8})
\begin{eqnarray}\label{an-eqn9}
\|\nabla^js\|_{L^p} & \le & C\|\nabla|\nabla^js|\,\|^a_{L^q}
\|\nabla^js\|^{1-a}_{L^r} \\
& \le & C\|\nabla^{j+1}s\|^a_{L^q}\|\nabla^js\|^{1-a}_{L^r} \nonumber
\end{eqnarray}
where $a\in (0,1)$ and
\begin{equation}\label{an-eqn10}
\frac{1}{p} = \frac{1}{r}+a\left(\frac{1}{q} - \frac{1}{r} -
\frac{1}{n}\right).
\end{equation}
Note that so far the condition $p\ge 2$ has not played a role.
\medskip

\noindent{\bf Case 2}: Let $j=1$, $k=2$ and $\frac{1}{2}\le a\le 1$. If
$a=1$ (\ref{an-eqn9}) yields
\begin{equation}\label{an-eqn11}
\|\nabla s\|_{L^p}\le C\|\nabla^2s\|_{L^q}
\end{equation}
with
\begin{equation}\label{an-eqn12}
\frac{1}{p}=\frac{1}{q}-\frac{1}{n}.
\end{equation}
If $a=\frac{1}{2}$, assume $p\ge 2$ then
\begin{gather}\label{an-eqn13}
\div\langle|\nabla s|^{p-2}\nabla s, s\rangle  =  |\nabla s|^p + 
|\nabla s|^{p-2}\langle \nabla_\a\nabla_\a s,s\rangle + \\
+ (p-2)|\nabla s|^{p-4}\langle \nabla_{\beta}s,
\nabla_\a\nabla_{\beta}s\rangle\langle\nabla_\a s,s\rangle. \nonumber
\end{gather}
Since
\begin{equation}\label{an-eqn14}
\int_{\RR^n}\div\langle|\nabla s|^{p-2}\nabla s, s\rangle=0
\end{equation}
then (\ref{an-eqn13}) gives
\begin{equation}\label{an-eqn15}
\int_{\RR^n}|\nabla s|^p\le (n+p-2)\int_{\RR^n}|\nabla
s|^{p-2}|\nabla^2s|\,|s|.
\end{equation}
Given our choice of $j=1$, $k=2$ and $a=\frac{1}{2}$ we have
$\frac{1}{q}+\frac{1}{r}=\frac{2}{p}$, i.e.
$\frac{1}{q}+\frac{1}{r}+\frac{p-2}{p}=1$.

Thus H\"older's inequality yields
\begin{equation}\label{an-eqn16}
\|\nabla s\|^p_{L^p}\le (n+p-2)\|\nabla^2s\|_{L^q}\|s\|_{L^r}\|\nabla
s\|^{p-2}_{L^p},
\end{equation}
thus
\begin{equation}\label{an-eqn17}
\|\nabla s\|_{L^p} \le \sqrt{n+p-2} \|\nabla^2
s\|^{\frac{1}{2}}_{L^q}\|s\|^{\frac{1}{2}}_{L^r}
\end{equation}
with
\begin{equation}\label{an-eqn18}
\frac{1}{p} = \frac{1}{2}\left(\frac{1}{r} + \frac{1}{n}\right).
\end{equation}

For $a\in\left(\frac{1}{2},1\right)$ we consider two cases: $q<n$, and $q\ge
n$. For $q<n$ using the convexity of $\log\|f\|^p_{L^p}$ as a function of
$p$ we have
\begin{equation}\label{an-eqn19}
\|\nabla s\|_{L^p} \le \|\nabla s\|^\a_{L^t} \|\nabla s\|^{1-\a}_{L^\s}
\mbox{ with } \a=\frac{p^{-1}-\s^{-1}}{t^{-1}-\s^{-1}}\in (0,1)
\end{equation}
where $t<p<\s$ are such that
\begin{equation}\label{an-eqn20}
\frac{2}{t} = \frac{1}{q} + \frac{1}{r} \mbox{ and } \frac{1}{\s} =
\frac{1}{q} - \frac{1}{n}.
\end{equation}

Using (\ref{an-eqn11}) and (\ref{an-eqn17}) we have that 
\begin{equation}\label{an-eqn21}
\|\nabla s\|_{L^\s}\le C\|\nabla^2s\|_{L^q}
\end{equation}
and
\begin{equation}\label{an-eqn22}
\|\nabla s\|_{L^t} \le C\|\nabla^2s\|^{\frac{1}{2}}_{L^q}
\|s\|_{L^r}^{\frac{1}{2}}.
\end{equation}
Combining (\ref{an-eqn19}), (\ref{an-eqn21}) and (\ref{an-eqn22}) we obtain
\begin{equation}\label{an-eqn23}
\|\nabla s\|_{L^p}\le C\|\nabla^2s\|^{1-\frac{\a}{2}}_{L^q}
\|s\|^{\frac{\a}{2}}_{L^r}
\end{equation}
where
\begin{equation}\label{an-eqn24}
\frac{1}{p} = \frac{1}{n} + \frac{1}{r} +
\left(1-\frac{\a}{2}\right)\,\left(\frac{1}{q} - \frac{1}{r} -
\frac{2}{n}\right),
\end{equation}
which proves the case $a\in\left(\frac{1}{2},1\right)$ and $q<n$.

For $q\ge n$, $t>0$ and $b\in(0,1)$ such that 
\begin{equation}\label{an-eqn25}
\frac{1}{p} = \frac{1}{t} + b\left(\frac{1}{q} - \frac{1}{t} -
\frac{1}{n}\right)
\end{equation}
we have by (\ref{an-eqn9})
\begin{equation}\label{an-eqn26}
\|\nabla s\|_{L^p} \le C\|\nabla^2s\|^b_{L^q} \|\nabla s\|^{1-b}_{L^t}.
\end{equation}
Choosing $t>0$ so that
\begin{equation}\label{an-eqn27}
\frac{2}{t} = \frac{1}{q}+\frac{1}{r}
\end{equation}
we have by (\ref{an-eqn17})
\begin{equation}\label{an-eqn28}
\|\nabla s\|_{L^t} \le C\|\nabla^2s\|^{\frac{1}{2}}_{L^q}
\|s\|^{\frac{1}{2}}_{L^r}.
\end{equation}
Combining (\ref{an-eqn26}) and (\ref{an-eqn28}) we obtain
\begin{equation}\label{an-eqn29}
\|\nabla s\|_{L^p} \le
C\|\nabla^2s\|^{\frac{b+1}{2}}_{L^q}\|s\|^{\frac{1-b}{2}}_{L^r}
\end{equation}
with
\begin{equation}\label{an-eqn30}
\frac{1}{p} = \frac{1}{n}+\frac{1}{r} + \left(\frac{b+1}{2}\right)\,
\left(\frac{1}{q} - \frac{1}{r} - \frac{2}{n}\right)
\end{equation}
by (\ref{an-eqn25}) and (\ref{an-eqn27}). This concludes the proof of Case
2.
\medskip

\noindent{\bf Case 3}: Let $j=0$ and $k=2$. From (\ref{an-eqn8}) we have
\begin{equation}\label{an-eqn31}
\|s\|_{L^p} \le C\|\nabla s\|^{a_1}_{L^{q_1}}\|s\|^{1-a_1}_{L^r}
\end{equation}
with $a_1\in(0,1)$ and
\begin{equation}\label{an-eqn32}
\frac{1}{p}=\frac{1}{r} + a_1 \left(\frac{1}{q_1} - \frac{1}{r} -
\frac{1}{n}\right).
\end{equation}

Choosing $q_1$ so that
\begin{equation}\label{an-eqn33}
\frac{1}{q_1} = \frac{1}{r} + \frac{1}{n} + a_2 \left( \frac{1}{q} -
\frac{1}{r} - \frac{2}{n}\right)
\end{equation}
then $a_2\in\left(\frac{1}{2},1\right)$ and
\begin{equation}\label{an-eqn34}
\|\nabla s\|_{L^{q_1}}\le C\|\nabla^2s\|^{a_2}_{L^q} \|s\|^{1-a_2}_{L^r}.
\end{equation}

Combining (\ref{an-eqn31}) and (\ref{an-eqn34}) we have that
\begin{equation}\label{an-eqn35}
\|s\|_{L^p} \le C\|\nabla^2s\|^{a_1a_2}_{L^q} \|s\|^{1-a_1a_2}_{L^r}
\end{equation}
with 
\begin{equation}\label{an-eqn36}
\frac{1}{p} = \frac{1}{r} + a_1a_2\left(\frac{1}{q} - \frac{1}{r} -
\frac{2}{n}\right)
\end{equation}
from (\ref{an-eqn32}) and (\ref{an-eqn33}).
\medskip

\noindent{\bf Case 4}: We now proceed by induction on $k$. Assume that for
$k\ge 2$ and $j<k$ we have proved (\ref{an-eqn5}). Let $j<k<k+1$. By
(\ref{an-eqn9}) we have
\begin{equation}\label{an-eqn37}
\|\nabla ^ks\|_{L^{q_1}} \le
C\|\nabla^{k+1}s\|^{a_2}_{L^{q_2}}\|\nabla^ks\|^{1-a_2}_{L^{r_2}}
\end{equation}
with
\begin{equation}\label{an-eqn38}
\frac{1}{q_1} = \frac{1}{r_2} + a_2
\left(\frac{1}{q_2}-\frac{1}{r_2}-\frac{1}{n}\right).
\end{equation}
By the induction hypothesis, applied to $\nabla^{k-1}s$, we also have
\begin{equation}\label{an-eqn39}
\|\nabla^ks\|_{L^{r_2}} \le C\|\nabla^{k+1}s\|^{a_3}_{L^{q_3}}
\|\nabla^{k-1}s\|^{1-a_3}_{L^{r_3}}
\end{equation}
with
\begin{equation}\label{an-eqn40}
\frac{1}{r_2} = \frac{1}{r_3} + \frac{1}{n} + a_3\left(\frac{1}{q_3} -
\frac{1}{r_3} - \frac{2}{n}\right)
\end{equation}
and
\begin{equation}\label{an-eqn41}
\|\nabla^{k-1}s\|_{L^{r_0}} \le C\|\nabla^ks\|^{a_4}_{L^{q_4}}
\|s\|^{1-a_4}_{L^{r_4}}
\end{equation}
with
\begin{equation}\label{an-eqn42}
\frac{1}{r_3} = \frac{1}{r_4} + \frac{k-1}{n} + a_4 \left(\frac{1}{q_4} -
\frac{1}{r_4} - \frac{k}{n}\right).
\end{equation}

Letting $q_4=r_2$, $q_3=q$, $r_4=r$, $r_2=p$ we obtain
\begin{equation}\label{an-eqn43}
\|\nabla^ks\|_{L^p} \le C\|\nabla^{k+1}s\|^a_{L^q}\|s\|^{1-a}_{L^r}\mbox{
with }a=\frac{a_3}{1-a_4+a_3a_4} \in \left[\frac{k}{k+1},1\right]
\end{equation}
and
\begin{equation}\label{an-eqn44}
\frac{1}{p} = \frac{1}{r} + \frac{k}{n} + a\left( \frac{1}{q} - \frac{1}{r}
+ \frac{k+1}{n}\right).
\end{equation}

By hypothesis for $j<k$ and using (\ref{an-eqn43}) we have
\begin{eqnarray}\label{an-eqn45}
\|\nabla^js\|_{L^p} & \le & C\|\nabla^ks\|^{a_1}_{L^{q_1}}
\|s\|^{1-a_1}_{L^r} \\
& \le & C\|\nabla^{k+1}s\|^{a_0a_1}_{L^q}\|s\|^{1-a_0a_1}_{L^r} \nonumber 
\end{eqnarray}
with $a_0a_1\in\left[\frac{j}{k+1}, 1\right]$ and
\begin{equation}\label{an-eqn46}
\frac{1}{p} = \frac{1}{r} + \frac{j}{n} + a_1a_0 \left(\frac{1}{q} -
\frac{1}{r} - \frac{k+1}{n}\right)
\end{equation}
which finishes the proof of the proposition.
\end{proof}

\begin{corollary}\label{an-cor1}
Let $u\in C^\infty(\RR^n,N)$ be constant outside a compact set. Then for
$k\ge 1$, $q,r\in[1,\infty)$ and $0\le j\le k-1$ we have
\begin{equation}\label{an-eqn47}
\|\nabla^{j+1}u\|_{L^p} \le C\|\nabla^ku\|^a_{L^q} \|\nabla u\|^{1-a}_{L^r}
\end{equation}
with
\begin{equation}\label{an-eqn48}
\frac{1}{p}=\frac{j}{n}+\frac{1}{r}+a\left(\frac{1}{q}-\frac{1}{r}-\frac{k-1}{n}\right).
\end{equation}
If $r=\frac{n}{k-1-j}\ne 1$ then (\ref{an-eqn47}) does not hold for $a=1$.
The constant $C$ that appears in (\ref{an-eqn48}) only depends on $n,k, j,
q, r$ and $a$.
\end{corollary}

\begin{proof}
Apply (\ref{an-eqn5}) to $s=\nabla u$ a section of the bundle
$u^\ast(TN)\otimes T^\ast\RR^n$. Since $\nabla u$ is not necessarily
compactly supported a standard approximation argument might be needed to
complete the proof.
\end{proof}

In the second part of this section we establish the equivalence of the
Sobolev norms defined in either the intrinsic, geometric setting or in the
ambient, Euclidean setting. These results hold when we are above the range
in which these spaces have suitable multiplication properties. Since we are
working with the gradients of the maps we must consider the $H^s$ spaces
with $s>\frac{n}{2}+1$.

We begin by assuming that we have chosen coordinate systems on $(N,g)$ so
that the eigenvalues of $g$ are bounded above and below by a fixed constant
$C>1$, i.e. we assume that 
\[
C^{-1}|\xi|^2\le g_{ij}\xi_i\xi_j\le C|\xi|^2\qquad\mbox{for all}\qquad
\xi\in\RR^k.
\]
We denote these coordinates by either $(y^1,\ldots, y^k)$ or $(u^1,\ldots,
u^k)$. As before $(x^1,\ldots, x^n)$ denotes Euclidean coordinates on
$\RR^n$.

For $v:\RR^n\to\RR^p$ we let
\[
\p_\a v=\frac{\p v^a}{\p x^\a}e_a
\]
where $\{e_1,\ldots,e_p\}$ is an orthonormal basis for $\RR^p$. Recall that
if $X\in\G(u^{-1}(TN))$ then
\[
(\nabla_\a X)^j=\frac{\p X^j}{\p x^\a} + \G^j_{ik} X^i\frac{\p u^k}{\p x^\a}
\]
and $\nabla_\a u=\p_\a u\in\G(u^{-1}(TN))$ denotes the vector field along
$u$ defined in \eqref{fo-eqn4}. We use the following notation for higher order
derivatives.

\begin{defn}\label{an-defn5}
Let $\a=(\a_1,\ldots,\a_{l+1})$ denote a multi-index of length
$l+1$ $(|\a|=l+1)$ with each $\a_s\in\{1,\ldots,n\}$. We let
$\nabla^{l+1}u\in\G(u^{-1}(TN))$ denote any covariant derivative of $u$ of
order $l+1$ e.g.
\[
\nabla^{l+1}u=\sum_{\a=(\a_1,\ldots,\a_{l+1})}\nabla_{\a_1}\cdots\nabla_{\a_{l+1}}u
\]
Similarly
\[
\p^{l+1}v=\sum_{\a=(\a_1,\ldots,\a_{l+1})}
\frac{\p^{l+1}v^a}{\p x^{\a_1}\cdots\p
x^{\a_{l+1}}}e_a
\]
and
\[
\p^{l+1}u=\sum_{\a=(\a_1,\ldots,\a_{l+1})}
\frac{\p^{l+1}u^a}{\p x^{\a_1}\cdots\p x^{\a_{l+1}}} \frac{\p}{\p
y^a}.
\]
\end{defn}

\begin{remark}\label{an-rem1}
Note that our use of the multi-index notation differs from the usual one. 
\end{remark}

Recall that for
\begin{eqnarray*}
v & : & \RR^n\to\RR^p \\
u & : & \RR^n\to N
\end{eqnarray*}
the Sobolev norms of $\p v$ and $\nabla u$ for $k\in\NN$ are defined by
\begin{eqnarray*}
\|\p v\|_{H^k} & = & \sum^{k}_{l=0} \|\p^{l+1} v\|_{L^2(\RR^n)} \\
\|\nabla u\|_{H^k} & = & \sum^k_{l=0}\|\nabla^{l+1}u\|_{L^2(\RR^n)} \\
& = &
\sum^{k}_{l=0}\left(\int_{\RR^n}g_{ij}(\nabla^{l+1}u)^i(\nabla^{l+1}u)^j\right)^{\frac{1}{2}}
\end{eqnarray*}
where here
\[
\|\nabla^{l+1}u\|_{L^2(\RR^n)} =
\sum_{|\a|=l+1}\|\nabla_{\a_1}\cdots\nabla_{\a_{l+1}}u\|_{L^2(\RR^n)}
\]
and the sum is taken over all distinct multi-indices of length $l+1$. The $L^2$ norm of each of these is computed with respect to
the metric $g$ as indicated. We use the obvious analogous definition for
$\|\p^{l+1}v\|_{L^2(\RR^n)}$.

Note that by definition $u\in H^k(\RR^n,N)$ if $\exists\, y_u\in N$ such
that for $v=w\circ u$
\[
\|v-w(y_u)\|_{L^2}+\|\p v\|_{H^{k-1}}<\infty.
\]
Our immediate goal is to show that for $k>\frac{n}{2}+1$ if $v=w\circ u$
then
\[
\|\p v\|_{H^k}<\infty\quad\mbox{if and only if}\quad \|\nabla
u\|_{H^k}<\infty.
\]

\begin{lemma}\label{an-lem5}
For each $k\ge 0$ we have
\begin{gather}
\nabla^{k+1}u=\p^{k+1}u+\sum^{k+1}_{l=2} \sum_{j_1+\cdots+j_l=k+1}
G_{(j_1,\ldots, j_l)}(u)\p^{j_1} u\ast\cdots\ast\p^{j_l}u \label{an-eqn31a} \\
\p^{k+1}u=\nabla^{k+1}u+\sum^{k+1}_{l=2} \sum_{j_1+\cdots+j_l=k+1}
E_{(j_1,\ldots, j_l)}(u)\nabla^{j_1}u\ast\cdots\ast\nabla^{j_l}u
\label{an-eqn32a} \\
\p^{k+1}v=\frac{\p w^a}{\p y^j}\p^{k+1}u^je_a+\sum^{k+1}_{l=2}
\sum_{j_1+\cdots+j_l=k+1} F_{(j_1,\ldots,
j_l)}(u)\p^{j_1}u\ast\cdots\ast\p^{j_l}u \label{an-eqn33a}
\end{gather}
where each subscript $j_s\ge 1$ and
\begin{eqnarray*}
G_{(j_1,\ldots, j_l)}(u) & = & G^j_{(j_1,\ldots,j_l)}(u) \frac{\p}{\p y^j}
\\
E_{(j_1,\ldots, j_l)}(u) & = & E^j_{(j_1,\ldots,j_l)}(u) \frac{\p}{\p y^j}
\\
F_{(j_1,\ldots, j_l)}(u) & = & F^a_{(j_1,\ldots,j_l)}(u)e_a
\end{eqnarray*}
and each $G$, $E$ and $F$ are smooth, bounded functions of $u$. 
\end{lemma}

The notation $a_{j_1}\ast\cdots\ast a_{j_l}$ corresponds to a product  of
the $a_{j_{k}}$'s.

\begin{remark}\label{an-rem2}
Throughout this section whenever expressions similar to the right hand side
of (\ref{an-eqn31a}), (\ref{an-eqn32a}) or (\ref{an-eqn33a}) occur, a key
point is to note that all the subscripts $j_s\ge 1$, for $s\in \{1, \ldots,
l\}$. This is always to be understood even if it is not explicitly stated.
\end{remark}

\begin{proof}
We establish each of these by induction, beginning with (\ref{an-eqn31a}).
Note that for $k=0$, $\nabla u=\p u$. For $k=1$
\begin{eqnarray*}
\nabla_{\a_2}\nabla_{\a_1}u & = & \nabla_{\a_2}\left(\frac{\p u^j}{\p
x^{\a_1}} \frac{\p}{\p y^j}\right) \\
& = & \frac{\p^2u^j}{\p x^{\a_2}\p x^{\a_1}} \frac{\p}{\p y^j} + \G^j_{ik}
\frac{\p u^i}{\p x^{\a_1}} \frac{\p u^k}{\p x^{\a_2}} \frac{\p}{\p y^j}.
\end{eqnarray*}
Then
\begin{eqnarray*}
\nabla^{k+2}u & = & \nabla_{\a_{k+2}}\nabla^{k+1}u \\
& = & \nabla_{\a_{k+2}}(\nabla_{\a_{k+1}}\cdots\nabla_{\a_1}u) \\
& = & \nabla_{\a_{k+2}}\left(\frac{\p^{k+1}u^j}{\p x^{\a_{k+1}}\cdots\p
x^{\a_1}} \frac{\p}{\p y^j}\right) \\
&& + \nabla_{\a_{k+2}}\left(\sum^{k+1}_{l=2} \sum_{j_1+\cdots+j_l=k+1}
G^j_{(j_1,\ldots,j_l)}(u)\p^{j_1}u\ast\cdots\ast\p^{j_l}u\frac{\p}{\p
y^j}\right) \\
& = & \p^{k+2}u + \G^j_{il}\frac{\p^{k+1}u^i}{\p x^{\a_{k+1}}\cdots\p
x^{\a_1}} \frac{\p u^l}{\p x^{\a_{k+2}}} \frac{\p}{\p y^j} \\
&& + \sum^{k+1}_{l=2} \sum_{j_1+\cdots+j_l=k+1}
\left(G^j_{(j_1,\ldots,j_l)}(u)\right)'\p^{j_1}u\ast\cdots\ast\p^{j_l}u
\frac{\p}{\p y^j} \\
&& + \sum^{k+2}_{l=2} \sum_{j_1+\cdots+j_l=k+2}
G^i_{(j_1,\ldots,j_l)}(u)\p^{j_1}u\ast\cdots\ast\p^{j_l}u\frac{\p}{\p y^j}
\\
&& + \sum^{k+2}_{l=2} \sum_{j_1+\cdots+j_l=k+2}
G^i_{(j_1,\ldots, j_l)}(u)\p^{j_1} u\ast\cdots\ast\p^{j_l}u\G^j_{il} \frac{\p
u^l}{\p x^{\a_{k+2}}} \frac{\p}{\p y^j}.
\end{eqnarray*}
Therefore
\[
\nabla^{k+2}u=\p^{k+2}u+\sum^{k+2}_{l=2} \sum_{j_1+\cdots+j_l=k+2}
G_{(j_1,\ldots,j_l)} (u) \p^{j_1}u\ast\cdots\ast\p^{j_l}u
\]
which completes the proof of (\ref{an-eqn31a}). The proof of
(\ref{an-eqn32a}) proceeds in a similar fashion and is left to the reader.

To prove (\ref{an-eqn33a}) we recall that $v=w\circ u$ and thus
\begin{equation}\label{an-eqn34a}
\frac{\p v^a}{\p x^{\a_1}} = \frac{\p w^a}{\p y^j} \frac{\p u^j}{\p
x^{\a_1}}
\end{equation}
(which is the case $k=0$). When $k=1$ we differentiate this to obtain
\[
\frac{\p^2v^a}{\p x^{\a_2}\p x^{\a_1}} = \frac{\p w^a}{\p y^j}
\frac{\p^2u^j}{\p x^{\a_2}\p x^{\a_1}} + \frac{\p^2w^a}{\p y^i\p y^j}
\frac{\p u^i}{\p x^{\a_2}} \frac{\p u^j}{\p x^{\a_1}}.
\]
Assume now that (\ref{an-eqn33a}) holds for some $k\ge 1$. Then
\begin{eqnarray*}
\p^{k+2}v & = & \p_{\a_{k+2}}(\p^{k+1}v) \\
& = & \frac{\p w^a}{\p y^j} \p^{k+2}u^je_a + \frac{\p^2w^a}{\p y^i\p y^j}
\p^{k+1}u^j \frac{\p u^i}{\p x^{\a_{k+2}}} e_a \\
&& + \sum^{k+2}_{l=2} \sum_{j_1+\cdots+j_l=k+2}
F_{(j_1,\ldots,j_l)}(u)\p^{j_1}u\ast\cdots\ast\p^{j_l}u.
\end{eqnarray*}
This implies (\ref{an-eqn33a}) and completes the proof of the Lemma.
\end{proof}

Combining (\ref{an-eqn32a}) and (\ref{an-eqn33a}) in Lemma \ref{an-lem5} we
obtain the following.

\begin{lemma}\label{an-lem6}
For $v=w\circ u$ and $k\ge 0$ we have
\begin{equation}\label{an-eqn35a}
\p^{k+1}v^a = \frac{\p w^a}{\p y^j} (\nabla^{k+1}u)^j + \sum^{k+1}_{l=2}
\sum_{j_1+\cdots+j_l=k+1} H^a_{(j_1,\ldots, j_l)}(u) \nabla^{j_1}
u\ast\cdots\ast\nabla^{j_l} u
\end{equation}
where, as before, each subscript $j_s\ge 1$ and each $H^a$ is a smooth,
bounded function of $u$.
\end{lemma}

We now proceed to bound the pointwise norms in terms of each other.

\begin{lemma}\label{an-lem7}
For $v=w\circ u$ and $k\ge 0$ there is a constant $C>1$ depending only on
$n$ and $k$ such that
\begin{equation}\label{an-eqn36a}
|\p^{k+1}v|^2\le C|\nabla^{k+1}u|^2+C\sum^{k+1}_{l=2}
\sum_{j_1+\cdots+j_l=k+1} |\nabla^{j_1}u|^2 \cdots |\nabla^{j_l}u|^2
\end{equation}
and
\begin{equation}\label{an-eqn37a}
|\nabla^{k+1}u|^2\le C|\p^{k+1}v|^2 + C\sum^{k+1}_{l=2}
\sum_{j_1+\cdots+j_l=k+1} |\p^{j_1}v|^2 \cdots |\p^{j_l}v|^2
\end{equation}
\end{lemma}

\begin{proof}
Using (\ref{an-eqn35a}) we have
\begin{align*}
\sum^p_{a=1}|\p^{k+1}v^a&|^2  =  \sum^p_{a=1} \sum_{|\a|=k+1}
|\p_{\a_{k+1}} \cdots\p_{\a_1}v^a|^2 \\
=&  \sum^p_{a=1} \frac{\p w^a}{\p y^j} (\nabla^{k+1} u)^j \frac{\p w^a}{\p
y^i}(\nabla^{k+1}u)^i \\
&+ \left|\sum^{k+1}_{l=2} \sum_{j_1+\cdots+j_l=k+1}
H^a_{(j_1,\ldots,j_l)}(u)\nabla^{j_1}u\ast\cdots\ast\nabla^{j_l}u\right|^2
\\
&+ 2 \frac{\p w^a}{\p y^j} (\nabla^{k+1}u)^j \sum^{k+1}_{l=2}
\sum_{j_1+\cdots+j_l=k+1} H^a_{(j_1,\ldots,j_l)} (u) \nabla^{j_1}u
\ast\cdots\ast\nabla^{j_l} u.
\end{align*}
Since $w:N\to\RR^p$ is an isometric embedding we note that
\begin{equation}\label{an-eqn38a}
g_{ij} = \sum^p_{a=1} \frac{\p w^a}{\p y^i} \frac{\p w^a}{\p y^j}.
\end{equation}
Therefore, using the fact that for any $l\ge 1$
\[
C^{-1} \sum^k_{i=1} |(\nabla^lu)^i|^2\le |\nabla^lu|^2 =
g_{ij}(\nabla^lu)^i(\nabla^lu)^j \le C\sum^k_{i=1} |(\nabla^lu)^i|^2
\]
we have
\[
\sum^p_{a=1} |\p^{k+1}v^a|^2 \le 2|\nabla^{k+1}u|^2+2C\sum^{k+1}_{l=2}
\sum_{j_1+\cdots+j_l=k+1} |\nabla^{j_1}u|^2\cdots|\nabla^{j_l}u|^2
\]
which establishes (\ref{an-eqn36a}). To prove (\ref{an-eqn37a}) we proceed by
induction. For $k=0$ we have
\[
\nabla_au=\frac{\p u^i}{\p x^\a} \frac{\p}{\p y^i}.
\]
Using (\ref{an-eqn34a}) and (\ref{an-eqn38a}) this implies
\[
|\nabla_\a u|^2 = g_{ij} \frac{\p u^i}{\p x^\a} \frac{\p u^j}{\p x^\a} =
\frac{\p v^a}{\p x^\a} \frac{\p va}{\p x^\a}.
\]
Therefore
\begin{equation}\label{an-eqn39a}
|\nabla u|^2=|\p v|^2.
\end{equation}
Note that for $k=1$, by (\ref{an-eqn35a}) we have
\[
\frac{\p w^a}{\p y^j} (\nabla^2u)^j = \p^2v^a-H^a(u)\nabla u\ast\nabla u.
\]
So that
\[
|\nabla^2 u|^2 \le 2|\p^2v|^2+C|\nabla u|^4
\]
or
\begin{equation}\label{an-eqn40a}
|\nabla^2u|^2\le 2|\p^2v|^2+C|\p v|^4.
\end{equation}
Assume now that (\ref{an-eqn37a}) holds for any $k\ge 1$. Again using
(\ref{an-eqn35a}) we then have
\begin{eqnarray*}
|\nabla^{k+2}u|^2 & \le & 2|\p^{k+2}v|^2 + C\sum^{k+2}_{l=2}
\sum_{j_1+\cdots+j_l=k+2} |\nabla^{j_1}u|^2 \cdots|\nabla^{j_l}u|^2 \\
& \le & 2|\p^{k+2}v|^2 + C\sum^{k+2}_{l=2} \sum_{j_1+\cdots+j_l=k+2}
|\p^{j_1}v|^2 \cdots |\p^{j_l}v|^2
\end{eqnarray*}
which completes the proof of Lemma \ref{an-lem7}.
\end{proof}

\begin{lemma}\label{an-lem8}
Assume that $k>\frac{n}{2}+1$. There exists a constant $C=C(N,k,n)$ such
that for $u\in C^\infty(\RR^n,N)$ constant outside a compact set of $\RR^n$
if $v=w\circ u$ then
\begin{eqnarray}
\|\nabla^{k+1}u\|_{L^2} & \le & C\sum^k_{l=1}\|\p v\|^l_{H^k} \label{an-eqn49}
\\
\|\p^{k+1}v\|_{L^2} & \le & C\sum^k_{l=1}\|\nabla u\|^l_{H^k}.
\label{an-eqn50}
\end{eqnarray}
\end{lemma}

\begin{proof}
By (\ref{an-eqn37a}) we have
\begin{equation}\label{an-eqn51}
\|\nabla^{k+1}u\|_{L^2} \le C\|\p^{k+1}v\|_{L^2} + C\sum^{k+1}_{l=2}
\sum_{j_1+\cdots+j_l=k+1} \left(\int_{\RR^n} |\p^{j_1}_v|^2 \cdots
|\p^{j_l}v|^2\right)^{\frac{1}{2}}.
\end{equation}
Let $2\le p_i\le\infty$, $i=1,\ldots l$ be such that
\begin{equation}\label{an-eqn52}
\frac{1}{p_1}+\cdots+\frac{1}{p_l}=\frac{1}{2}.
\end{equation}
Then by H\"older's inequality
\begin{equation}\label{an-eqn53}
\| |\partial^{j_1}v|\cdots|\partial^{j_l}v|\|_{L^2}\le
c\|\partial^{j_1}v\|_{L^{p_1}} \cdots \|\partial^{j_l}v\|_{L^{p_l}}.
\end{equation}
Since $k\ge\frac{n}{2}+1$ then
\begin{equation}\label{an-eqn54}
\frac{j_i-1}{k}<a_i=\frac{j_i-1}{k}+\frac{n}{2k^2}
\left(k-j_i+\frac{1}{l}\right)<1
\end{equation}
and
\begin{equation}\label{an-eqn55}
\frac{1}{2}\ge
\frac{1}{p_i}=\frac{j_i-1}{n}+\frac{1}{2}-\frac{ka_i}{n}>0\quad .
\end{equation}
Note that to ensure that $a_i<1$ in (\ref{an-eqn54}) we either need $n\le 3$
or we must have $\left(\frac{n}{2k}-1\right)\left(k-j_i+\frac{1}{l}\right)<1-\frac{1}{l}$.
Since $2\le l\le k+1$ and $1\le j_i \le k$ the previous inequality holds
provided $\left(\frac{n}{2k}-1\right)\left(k-\frac{1}{2}\right)<\frac{1}{2}$
which requires $k>\frac{n+\sqrt{n(n-4)}}{4}$. Thus to accommodate all values
of $n$ simultaneously, it is enough to choose $k\ge \frac{n}{2}+1$ and
$k\in\NN$.
Thus (\ref{an-eqn5}) in Proposition 4.1 yields
\begin{eqnarray}\label{an-eqn56}
\|\partial^{j_i}v\|_{L^{p_i}} & \le &
C\|\partial^{k+1}v\|^{a_i}_{L^2}\|\partial v\|^{1-a_i}_{L^2} \\
& \le & C\|\partial v\|_{H^k}. \nonumber
\end{eqnarray}
Therefore combining (\ref{an-eqn51}), (\ref{an-eqn53}) and (\ref{an-eqn56})
we have
\begin{equation}\label{an-eqn57}
\|\nabla^{k+1}u\|_{L^2}\le C\sum^k_{\rho=1}\|\partial v\|^\rho_{H^k}.
\end{equation}
A similar argument to the one above where Proposition \ref{an-prop1} is now
applied to $\nabla u$ rather than $\partial v$ yields (\ref{an-eqn50}).
\end{proof}

\begin{lemma}\label{an-lem8A}
There exists a constant $C=C(N,n)$ such that if $u\in C^\infty(\RR^n,N)$ is
constant outside a compact set of $\RR^n$ and $v=w\circ u$ then for $1\le
k\le\frac{n}{2}+1$
\begin{eqnarray}
\|\nabla^{k+1}u\|_{L^2} & \le & C\sum^k_{l=1}\|\partial
v\|^l_{H^{\left[\frac{n}{2}\right]+2}} \label{an-eqn49A} \\
\|\partial^{k+1}v\|_{L^2} & \le & C\sum^k_{l=1}\|\nabla
u\|^l_{H^{\left[\frac{n}{2}\right]+2}}. \label{an-eqn50A}
\end{eqnarray}
\end{lemma}

\begin{proof}
The proof is very similar to that of Lemma \ref{an-lem8}, where the $a_i$'s
and $p_i$'s in the interpolation are taken as follows
\begin{equation}\label{an-eqn51A}
\frac{j_i-1}{s_0}<a_i = \frac{j_i-1}{s_0}+\frac{n}{2ks_0}(k-j_i+l^{-1})<1,
\end{equation}
where $s_0=\left[\frac{n}{2}\right]+2$, and
\begin{equation}\label{an-eqn52A}
\frac{1}{2}\ge \frac{1}{p_i}=\frac{j_i-1}{n} + \frac{1}{n} -
\frac{s_0}{n}a_i>0.
\end{equation}
\end{proof}

\begin{rem}
Proposition \ref{an-prop1} holds for $s\in C^m_c(E)$ where $E$ is a finite
dimensional $C^m$ vector bundle over $\RR^n$ provided $k<m$. Similarly Lemma
\ref{an-lem8} holds for $u\in C^m(\RR^n,N)$ and $u$ constant outside a
compact set of $\RR^n$, provided once again that $m>k$. A simple
approximation theorem ensures that Lemma \ref{an-lem8} holds for $u\in
C^m(\RR^n,N)\cap H^k(\RR^n,N)$ with $m>k$.
\end{rem}

\begin{corollary}\label{an-cor4.1}
Assume that $k\ge\frac{n}{2}+4$. There exists a constant $C=C(N,k,n)$ such
that for $u\in C^{k+1}(\RR^n,N)\cap H^k(\RR^n,N)$
\begin{eqnarray}
\|\nabla u\|_{L^\infty} & \le & C\sum^{\left[\frac{n}{2}\right]+2}_{j=1}
\|\nabla u\|^j_{H^{\left[\frac{n}{2}\right]+2}} \label{an-eqn58} \\
\|\nabla^2 u\|_{L^\infty} & \le & C\sum^{2\left[\frac{n}{2}\right]+4}_{l=1}
\|\nabla u\|^l_{H^{\left[\frac{n}{2}\right]+3}} \label{an-eqn59} \\
\|\nabla^3 u\|_{L^\infty} & \le & C\sum^{3\left[\frac{n}{2}\right]+12}_{l=1}
\|\nabla u\|^l_{H^{\left[\frac{n}{2}\right]+4}} .\label{an-eqn59A}
\end{eqnarray}
\end{corollary}

\begin{proof}
Recall that $\|\nabla u\|=|\partial v|$ if $v=w\circ u$, and by Sobolev
embedding theorem
\begin{align}
\|\partial v\|_{L^\infty} & \le  c\|\partial
v\|_{H^{\left[\frac{n}{2}\right]+2}} \label{an-eqn60} \\
\|\partial^2v\|_{L^\infty} & \le c\|\partial
v\|_{H^{\left[\frac{n}{2}\right]+3}} \label{an-eqn61} \\
\intertext{and}
\|\partial^3v\|_{L^\infty} & \le c\|\partial
v\|_{H^{\left[\frac{n}{2}\right]+4}} \label{an-eqn61A}
\end{align}

Therefore combining (\ref{an-eqn50}), (\ref{an-eqn50A}) and (\ref{an-eqn60})
we have
\begin{eqnarray}\label{an-eqn62}
\|\nabla u\|_{L^\infty} & \le & C\sum^{\left[\frac{n}{2}\right]+2}_{j=0}
\|\partial^{j+1}v\|_{L^2} \\
& \le & C\left(\|\nabla u\|_{L^2} + \sum^{\left[\frac{n}{2}\right]+2}_{j=1}
\|\nabla u\|^j_{H^{\left[\frac{n}{2}\right]+2}}\right) \nonumber \\
& \le & C\sum^{\left[\frac{n}{2}\right]+2}_{j=1} \|\nabla
u\|^j_{H^{\left[\frac{n}{2}\right]+2}}. \nonumber
\end{eqnarray}
Note that (\ref{an-eqn37a}) ensures that 
\begin{equation}\label{an-eqn63}
|\nabla^2u|\le C|\partial^2v|+c|\partial v|^2.
\end{equation}

Combining (\ref{an-eqn60}), (\ref{an-eqn61}), (\ref{an-eqn50}) and
(\ref{an-eqn50A}) we have
\begin{eqnarray}\label{an-eqn64}
\|\nabla^2 u\|_{L^\infty} & \le & C\|\partial^2v\|_{L^\infty} + C\|\partial
v\|^2_{L^\infty} \\
& \le & C\|\partial v\|_{H^{\left[\frac{n}{2}\right]+3}} + C\|\partial
v\|^2_{H^{\left[\frac{n}{2}\right]+2}} \nonumber \\
& \le & C\sum^{\left[\frac{n}{2}\right]+3}_{j=0} \|\partial^{j+1}v\|_{L^2} +
C\sum^{\left[\frac{n}{2}\right]+2}_{j=0}\|\partial^{j+1}v\|^2_{L^2}
\nonumber \\
& \le & C \sum^{\left[2\frac{n}{2}\right]+4}_{l=1} \|\nabla
u\|^l_{H^{\left[\frac{n}{2}\right]+3}} \nonumber
\end{eqnarray}

Note that (\ref{an-eqn37a}) also ensures that
\begin{equation}\label{an-eqn65}
|\nabla^3u|\le C(|\partial^3v|+|\partial^2v|\,|\partial v|+|\partial v|^3).
\end{equation}
Combining (\ref{an-eqn60}), (\ref{an-eqn61}), (\ref{an-eqn61A}),
(\ref{an-eqn50}) and (\ref{an-eqn50A}) we have
\begin{equation}\label{an-eqn66}
\|\nabla^3u\|_{L^\infty} \le
C\sum^{3\left(\left[\frac{n}{2}\right]+4\right)}_{j=1} \|\nabla
u\|^j_{H^{\left[\frac{n}{2}\right]+4}}.
\end{equation}
\end{proof}

\section{$\e$-independent energy estimates}

Theorem \ref{du-thm2} ensures that the initial value problem
\[
\left\{\begin{array}{rcl}
\frac{\p v}{\p t} & = & -\e\D^2v+N(v) \\
v(0) & = & v_0
\end{array}\right.
\]
has a unique solution $v_\e\in C([0,T_\e],H^{s+1}\cap L^{2,\infty}_\d)$ provided
$v_0\in H^{s+1}(\RR^n,\RR^p)$ for $s>[\frac{n}{2}]+4$. To prove that
(\ref{int-eqn3}) has a solution we need to show that (\ref{fo-eqn13}) has a
solution for $\e=0$. To do this we need to show that each $v_\e$ extends to
a solution in $C([0,T],H^{s+1}\cap L^{2,\infty}_\d)$ where $T>0$ is independent of
$\e$. This is accomplished by proving $\e$-independent energy estimates for
the function $v_\e$. It turns out that
thanks to the geometric nature of this flow, if one assumes enough regularity (i.e. $s>[\frac{n}{2}]+4$), it is easier to prove
$\e$-independent energy estimates for the corresponding $u_\e$. Lemma
\ref{an-lem8} and Lemma \ref{an-lem8A}
 then allows us to translate these into estimates for $v_\e$.

Let $u_\e=u\in C([0,T_\e],H^{s+1}(\RR^n,N))$ with $s$ large enough\footnote{ We will see later that 
$s>[\frac{n}{2}]+4$ will be enough. In this paper we do not attempt to obtain the lowest possible exponent $s$.} be a
solution of
\begin{equation}\label{in-eqn5.1}
\left\{\begin{array}{rcl}
\p_tu & = & -\e\D\tau(u)+\e R(\nabla u, \tau(u))\nabla u +
J(u)\tau(u)+\beta\tau(u) \\
u(0) & = & u_0,\end{array} 
\right.
\end{equation}
where $\e\in(0,1]$, $\beta\ge 0$, $\D=\sum^n_{\a=1}\nabla_\a\nabla_\a$. Our
goal is to understand how $\|\nabla u\|_{H^k}(t)$ varies with time.

Let $l\in\NN$. We denote by $\ba$ the multi-index of length $l$,
$\ba=(\a_1\cdots\a_l)$, and
$\nabla_{\ba}u=\nabla_{\a_1\cdots}\nabla_{\a_l}u$. The following lemma and
corollaries establish some computational identities which are very useful.

\begin{lemma}\label{in-lemma12}
Let $u\in C^1([0,T],H^s(\RR^n,N))$, $s\in\NN$, $s>\frac{n}{2}+2$. Let $X\in
TN$ for $1\le l\le s$ and $|\ba|=l$. We have
\begin{gather}
\nabla_{\a_0}\nabla_{\ba}u = \nabla_{\ba}\nabla_{\a_0}u + \sum^{l-2}_{j=0}
\nabla_{\a_1}\cdots\nabla_{\a_j}[R(\nabla_{\a_0}u,
\nabla_{\a_{j+1}}u)\nabla_{\a_{j+2}}\cdots\nabla_{\a_l}u] \label{in-eqn1}\\ 
\nabla_t\nabla_{\ba}u=\nabla_{\ba}\nabla_tu + \sum^{l-2}_{j=0}
\nabla_{\a_1}\cdots\nabla_{\a_j} [R(\nabla_tu,
\nabla_{\a_{j+1}}u)\nabla_{\a_{j+2}}\cdots\nabla_{\a_l}u]\label{in-eqn1a}\\ 
\nabla_{\a_0}\nabla_{\ba}X = \nabla_{\ba} \nabla_{\a_0} X + \sum^{l-1}_{j=0}
\nabla_{\a_1}\cdots\nabla_{\a_j} [R(\nabla_{\a_0}u,
\nabla_{\a_{j+1}}u)\nabla_{\a_{j+2}}\cdots\nabla_{\a_l}X].\label{in-eqn2}
\end{gather}
\end{lemma}

\begin{proof}
The proof is done by induction on the length of the multi-index $\ba$, i.e.,
on $l$. We prove (\ref{in-eqn2}) and leave (\ref{in-eqn1}) and
(\ref{in-eqn1a}) to the reader, as the proofs are very similar. If
$l=1$
\begin{equation}\label{in-eqn3}
\nabla_{\a_0}\nabla_{\a_1}X = \nabla_{\a_1}\nabla_{\a_0}X+R(\nabla_{\a_0}u,
\nabla_{\a_1}u)X.
\end{equation}
Suppose (\ref{in-eqn2}) holds for $l\ge 1$ and consider
\begin{eqnarray}\label{in-eqn4}
&&\kern-.4in\nabla_{\a_0}\nabla_{\a_1}\cdots \nabla_{\a_{l+1}}X\\
& = & \nabla_{\a_1}
[\nabla_{\a_0} \nabla_{\a_2}\cdots \nabla_{\a_{l+1}}X] + R(\nabla_{\a_0}u,
\nabla_{\a_1}u)\nabla_{\a_2}\cdots\nabla_{\a_{l+1}}X \nonumber \\
& = & \nabla_{\a_1}\cdots \nabla_{\a+1}\nabla_{\a_0}X 
 +
R(\nabla_{\a_0}u, \nabla_{\a_1}u)\nabla_{\a_2}\cdots\nabla_{\a_{l+1}}X
\nonumber \\
&&\quad+ \nabla_{\a_1}
[\sum^{l}_{j=1}\nabla_{\a_2}\cdots\nabla_{\a_j} [R(\nabla_{\a_0}u,
\nabla_{\a_{j+1}}u)\nabla_{\a_{j+2}}\cdots\nabla_{\a_{l+1}}X] 
\nonumber \\
& = & \nabla_{\a_1}\cdots\nabla_{\a_{l+1}}\nabla_{\a_0}X\nonumber \\
&&\quad + \sum^l_{j=0}
\nabla_{\a_1}\cdots\nabla_{\a_j}[R(\nabla_{\a_0}u,
\nabla_{\a_{j+1}}u)\nabla_{\a_{j+2}}\cdots\nabla_{\a_{l+1}}X]\nonumber
\end{eqnarray}
\end{proof}

\begin{corollary}\label{in-cor1a}
Let $u\in C^1([0,T],H^s(\RR^n,N))$ $s\in \NN$, $s>\frac{n}{2}+2$, then for
$1\le l\le s$, $|\ba|=l$ we have
\begin{eqnarray}
\ \ \ \ \ \D\nabla_{\ba}u & = & \nabla_{\ba}\tau(u) +
\sum^{l-1}_{j=0}\nabla_{\a_1}\cdots\nabla_{\a_j}[R(\nabla_{\a_0}u,
\nabla_{\a_{j+1}}u)\nabla_{\a_{j+2}}\cdots\nabla_{\a_l}\nabla_{\a_0}u]
\label{in-eqn5}\\
&&\;+
\sum^{l-2}_{j=1}\nabla_{\a_0}\nabla_{\a_1}\cdots\nabla_{\a_j}[R(\nabla_{\a_0}u,
\nabla_{\a_{j+1}}u)\nabla_{\a_{j+2}}\cdots\nabla_{\a_l}u] \nonumber 
\end{eqnarray}
\begin{eqnarray}
\nabla_t\nabla_{\ba}\nabla_{\a_0}u & = & \nabla_{\ba}\nabla_{\a_0}\nabla_tu
\label{in-eqn6}\\
&&\;+ \sum^{l-1}_{j=0}
\nabla_{\a_1}\cdots\nabla_{\a_j}[R(\nabla_tu,\nabla_{\a_{j+1}}u)\nabla_{\a_{j+2}}\cdots\nabla_{\a_l}\nabla_{\a_0}u]
\nonumber
\end{eqnarray}
\begin{eqnarray}
\nabla_{\beta_0}\nabla_{\ba}\nabla_{\a_0}u & = &
\nabla_{\ba}\nabla_{\a_0}\nabla_{\beta_0}u \label{in-eqn7} \\
&&\;+ \sum^{l-1}_{j=0}
\nabla_{\a_1}\cdots\nabla_{\a_j}[R(\nabla_{\beta_0}u,
\nabla_{\a_{j+1}}u)\nabla_{\a_{j+2}}\cdots\nabla_{\a_l}\nabla_{\a_0}u]
\nonumber
\end{eqnarray}
\end{corollary}

\begin{proof}
The proof of (\ref{in-eqn5}) is an application of (\ref{in-eqn1}) and
(\ref{in-eqn2}). To prove (\ref{in-eqn6}) and (\ref{in-eqn7}) apply
(\ref{in-eqn2}) to $\nabla_{\a_0}u=X$ and note that $\nabla_t$ and
$\nabla_{\beta_0}$ behave the same way. Moreover recall that
$\nabla_{\a_0}\nabla_tu=\nabla_t\nabla_{\a_0}u$. 
\end{proof}

\begin{corollary}\label{in-cor2a}
Let $u\in C^1([0,T],H^s(\RR^n,N))$, $s\in\NN$ $s>\frac{n}{2}+2$. Let $X\in TN$
then for $l\ge 1$ and $|\ba|=l$ we have
\begin{eqnarray}\label{in-eqn8}
\D\nabla_{\ba} X&=& \nabla_{\ba}\D X \\
&& + \sum^{l-1}_{j=0}\nabla_{\a_1} \cdots
\nabla_{\a_j}[R(\nabla_{\a_0}u, \nabla_{\a_{j+1}}u) \nabla_{\a_{j+2}}\cdots
\nabla_{\a_l} \nabla_{\a_0}X] \nonumber\\
&& + \sum^{l-1}_{j=0}\nabla_{\a_0}\nabla_{\a_1}\cdots\nabla_{\a_j}
[R(\nabla_{\a_0}u, \nabla_{\a_{j+1}}u) \nabla_{\a_{j+2}} \nabla_{\a_l}X].
\nonumber
\end{eqnarray}
\end{corollary}

\begin{proof}
To prove (\ref{in-eqn8}) we apply (\ref{in-eqn2}) twice, first to $X$ then
$\nabla_{\a_0}X$.
\begin{eqnarray}\label{in-eqn9}
\D\nabla_{\ba}X & = & \nabla_{\a_0}\nabla_{\a_0} \nabla_{\ba}X \\
& = &
\nabla_{\a_0} [\nabla_{\ba}\nabla_{\a_0}X + \nabla_{\a_0}
[\sum^{l-1}_{j=0} \nabla_{\a_1}\cdots \nonumber\\
&& \cdots\nabla_{\a_j}(R(\nabla_{\a_0}u,
\nabla_{\a_{j+1}}u)\nabla_{\a_{j+2}}\cdots\nabla_{\a_l}X] \nonumber \\
& = & \nabla_{\ba} \nabla_{\a_0} \nabla_{\a_0}X + \sum^{l-1}_{j=0}
\nabla_{\a_1} \cdots  \nonumber\\
&&\cdots \nabla_{\a_j} (R(\nabla_{\a_0}u, \nabla_{\a_{j+1}}u)
\nabla_{\a_{j+2}}\nabla_{\a_l} \nabla_{\a_0}X) \nonumber \\
&+&\sum^{l-1}_{j=0} \nabla_{\a_0}\nabla_{\a_1}\cdots\nabla_{\a_j}
[R(\nabla_{\a_0}u, \nabla_{\a_{j+1}}u)\nabla_{\a_{j+2}} \cdots
\nabla_{\a_l}X]. \nonumber
\end{eqnarray}
\end{proof}

\begin{rem}
Note that in particular (\ref{in-eqn8}) applied to $X=\tau(u)$ yields
\begin{eqnarray}\label{in-eqn10}
\kern.5in\D\nabla_{\ba}\tau(u) & = & \nabla_{\ba}\D\tau(u) \\
&&\kern-.2in+\sum^{l-1}_{j=0}\nabla_{\a_1}\cdots\nabla_{\a_j}[R(\nabla_{\a_0}u,
\nabla_{\a_{j+1}}u)\nabla_{\a_{j+2}}\cdots\nabla_{\a_l}\nabla_{\a_0}\tau(u)]
\nonumber \\
&&\kern-.2in+\sum^{l-1}_{j=0}\nabla_{\a_0}\nabla_{\a_1}\cdots\nabla_{\a_j}[R(\nabla_{\a_0}u,
\nabla_{\a_{j+1}}u)\nabla_{\a_{j+2}}\cdots\nabla_{\a_l} \tau(u)]\nonumber
\end{eqnarray}
\end{rem}

\begin{lemma}\label{in-lem1}
Let $u\in C^1([0,T], H^s(\RR^n,N))$ with $s\in\NN$ and
$s\ge\left[\frac{n}{2}\right]+4$ be a solution of (\ref{in-eqn5.1}). Then for
$\left[\frac{n}{2}\right]+4\le l\le s$ and $l\in\NN$ we have
\begin{equation}\label{in-eqn5.1ast}
\frac{d}{dt}\|\nabla^lu\|^2_{L^2}  \le  C\|\nabla
u\|^2_{H^{l-1}}(1+\|\nabla u\|^{3n+2l+14}_{H^{l-1}}).
\end{equation}
For $l=\left[\frac{n}{2}\right]+2$ and $l=\left[\frac{n}{2}\right]+3$ we have
\begin{equation}\label{in-eqn63}
\frac{d}{dt}\|\nabla^lu\|^2_{L^2} \le C\left(1+\|\nabla
u\|^{3n+12}_{H^{\left[\frac{n}{2}\right]+4}}\right) \|\nabla
u\|^2_{H^{l-1}}\left(1+\|\nabla u\|^{2l+2}_{H^{l-1}}\right).
\end{equation}
\end{lemma}

\begin{proof}
We first compute the evolution 
\begin{eqnarray}\label{in-eqn5.2}
\kern.5in\frac{1}{2}\frac{d}{dt}\int_{\RR^n}|\nabla u|^2dx & = &
\sum^n_{\a_0=1}\int\langle \nabla_t\nabla_{\a_0}u, \nabla_{\a_0}u\rangle \\
&&\kern-.5in= \int\nabla_{\a_{0}}(\langle \nabla_tu, \nabla_{\a_0}u\rangle) -
\int\langle\nabla_tu, \tau(u)\rangle 
\nonumber \\
&&\kern-.5in= \e\int\langle\D\tau(u),\tau(u)\rangle-\e\int\langle R(\nabla
u,\tau(u))\nabla u,\tau(u)\rangle \nonumber\\
&&- \int\langle J(u)\tau(u),\tau(u)\rangle -
\beta\int|\tau(u)|^2 \nonumber \\
&&\kern-.5in= -\e\int|\nabla\tau(u)|^2-\beta\int|\tau(u)|^2 \nonumber\\
&&- \e\int\langle R(\nabla
u,\tau(u))\nabla u,\tau(u)\rangle,
\nonumber
\end{eqnarray}
where we have used the fact that for $f\in L^1(\RR^n,\RR^n)\int_{\RR^n}\div
f=0$ as well as integration by parts. Note that using integration by parts
and Cauchy-Schwarz we have
\begin{eqnarray}\label{in-eqn5.3}
\left|\int\langle R(\nabla u,\tau(u))\nabla u,\tau(u)\rangle\right| & \le & C\|\nabla
u\|^2_{L^\infty}\int|\tau(u)|^2 \\
& \le & C\|\nabla u\|^2_{L^\infty} \|\nabla u\|_{L^2}\|\nabla\tau(u)\|_{L^2}
\nonumber \\
& \le & \frac{1}{2}\|\nabla\tau(u)\|^2_{L^2}+C\|\nabla
u\|^4_{L^\infty}\|\nabla u\|^2_{L^2}. \nonumber
\end{eqnarray}
Combining (\ref{in-eqn5.2}) and (\ref{in-eqn5.3}) we have 
\begin{equation}\label{in-eqn5.4}
\frac{1}{2}\frac{d}{dt} \|\nabla u\|^2_{L^2}\le C\|\nabla u\|^4_{L^\infty}
\|\nabla u\|^2_{L^2}.
\end{equation}
For $1\le l\le s$ applying (\ref{in-eqn1a}) we have
\begin{eqnarray}\label{in-eqn5.5}
\kern.2in\frac{1}{2}\frac{d}{dt}\|\nabla^l u\|^2_{L^2} & = &
\sum_{|\ba|=l}\int\langle\nabla_t\nabla_{\ba} u, \nabla_{\ba} u\rangle \\
& = & \sum_{|\ba|=l}\int\langle\nabla_{\ba}\nabla_tu,
\nabla_{\ba}u\rangle\nonumber \\
&&
+ \sum_{|\ba|=l}\sum^{l-2}_{j=0}\int\langle\nabla_{\a_1}\cdots\nabla_{\a_j}
[R(\nabla_tu,
\nabla_{\a_{j+1}}u)\nabla_{\a_{j+2}}\cdots \nonumber \\
&&\cdots\nabla_{\a_l}u],
\nabla_{\ba}u\rangle.\nonumber
\end{eqnarray}
Consider each term separately
\begin{eqnarray}\label{in-eqn5.6}
\int\langle\nabla_{\ba}\nabla_tu,\nabla_{\ba}u\rangle & = &
-\e\int\langle\nabla_{\ba}\D\tau(u),\nabla_{\ba}u\rangle \\
&&+\e\int\langle\nabla_{\ba}(R(\nabla u,\tau(u))\nabla
u,\nabla_{\ba}u\rangle \nonumber \\
&&+\int\langle\nabla_{\ba}J(u)\tau(u),\nabla_{\ba}u\rangle +
\beta\int\langle\nabla_{\ba}\tau(u),\nabla_{\ba}u\rangle.\nonumber
\end{eqnarray}
Using (\ref{in-eqn10}) and (\ref{in-eqn5}) and integrating by parts we have
that
\begin{eqnarray}\label{in-eqn5.7}
&&\kern.2in\int\langle\nabla_{\ba}\D\tau(u),\nabla_{\ba}u\rangle \\
& = & \int\langle\nabla_{\ba}\tau(u),\D\nabla_{\ba}u\rangle \nonumber \\
&&-\sum^{l-1}_{j=0}\int\langle\nabla_{\a_1}\cdots\nabla_{\a_j}
[R(\nabla_{\a_0}u,\nabla_{\a_{j+1}}u)\nabla_{\a_{j+2}}\cdots\nabla_{\a_l}
\nabla_{\a_0}\tau(u)],\nabla_{\ba}u\rangle \nonumber \\
&&
-\sum^{l-1}_{j=0}\int\langle\nabla_{\a_0}\nabla_{\a_1}\cdots\nabla_{\a_j}
[R(\nabla_{\a_0}u,\nabla_{\a_{j+1}}u)\nabla_{\a_{j+2}}\cdots\nabla_{\a_l}
\tau(u)],\nabla_{\ba}u\rangle \nonumber \\
& = & \int|\nabla_{\ba}\tau(u)|^2 \nonumber\\
&&+ \sum^{l-1}_{j=0}\int\langle\nabla_{\a_1}\cdots\nabla_{\a_j}
[R(\nabla_{\a_0}u,\nabla_{\a_{j+1}}u)\nabla_{\a_{j+2}}\cdots\nabla_{\a_l}\nabla_{\a_0}u],
\nabla_{\ba}\tau(u)\rangle \nonumber \\
&&- \sum^{l-1}_{j=0}\int\langle\nabla_{\a_1}\cdots\nabla_{\a_j}
[R(\nabla_{\a_0}u,
\nabla_{\a_{j+1}}u)\nabla_{\a_{j+2}}\cdots\nabla_{\a_l}\nabla_{\a_0}
\tau(u)],\nabla_{\ba}u\rangle \nonumber \\
&&-\sum^{l-1}_{j=0}\int\langle\nabla_{\a_0}\nabla_{\a_1}\cdots\nabla_{\a_j}
[R(\nabla_{\a_0}u,\nabla_{\a_{j+1}}u)\nabla_{\a_{j+2}}\cdots\nabla_{\a_l}\tau(u)],
\nabla_{\ba}u\rangle \nonumber \\
&&+\sum^{l-2}_{j=1}\langle\nabla_{\a_0}\nabla_{\a_1}
\cdots\nabla_{\a_j}[R(\nabla_{\a_0}u,
\nabla_{\a_{j+1}}u)\nabla_{\a_{j+2}}\cdots\nabla_{\a_l}u],
\nabla_{\ba}\tau(u)\rangle \nonumber
\end{eqnarray}
(\ref{in-eqn5.7}) yields
\begin{eqnarray}\label{in-eqn5.8}
-\e\sum_{|\ba|=l}\int \langle\nabla_{\ba}\D\tau(u),\nabla_{\ba}u\rangle & \le &
-\e\int|\nabla^l\tau(u)|^2 \\
&& \kern-1.5in+
C\e\sum^{l+2}_{m=3}\quad\sum_{\mathop{j_1+\cdots+j_m=l+2}\limits_{\scriptscriptstyle{j_s\ge
1}}}
\int|\nabla^l\tau(u)|\,|\nabla^{j_1}u|\cdots|\nabla^{j_m}u|
\nonumber\\
&&\kern-1.5in+C\e\sum^{l+2}_{m=3}\quad\sum_{\mathop{j_1+\cdots+j_m=l+2}\limits_{\scriptscriptstyle{j_s\ge
1}}}
\int|\nabla^{j_1}\tau(u)|\,|\nabla^{j_2}u|\cdots|\nabla^{j_m}u|\,
|\nabla^lu|. \nonumber
\end{eqnarray}
Similarly
\begin{gather}\label{in-eqn5.9}
\sum_{|\ba|=l}\int\langle\nabla_{\ba}[R(\nabla u,\tau(u))\nabla u],
\nabla_{\ba}u\rangle\le  \\
\le C\sum^{l+2}_{m=3}\quad
\sum_{\mathop{j_1+\cdots+j_m=l+2}\limits_{\scriptscriptstyle{j_s\ge
1\;\mbox{\tiny if }s\ge 2}}}\int|\nabla^l
u|\,|\nabla^{j_1}\tau(u)|\,|\nabla^{j_2}u|\cdots|\nabla^{j_m}u| \nonumber
\end{gather}

We now look at the third term in (\ref{in-eqn5.6}) and recall that $\nabla
J=0$ and $\langle JX,X\rangle=0$ for $X\in TN$. Integrating by parts and
applying (\ref{in-eqn5}) we obtain for $\gamma=(\a_2\cdots\a_l)$
\begin{eqnarray}\label{in-eqn5.10}
&&\kern.5in\int\langle\nabla_{\ba}J(u)\tau(u), \nabla_{\ba}u\rangle \\
& = &
-\int\langle\nabla_{\bg}J(u)\tau(u),\nabla_{\a_1}\nabla_{\ba}u\rangle 
\nonumber \\
& = & -\int \langle\nabla_{\bg}J(u)\tau(u),\D\nabla_{\bg}u\rangle \nonumber \\
& = & -\int\langle J(u)\nabla_{\bg}\tau(u);\nabla_{\bg}\tau(u)\rangle
\nonumber \\
&&- \sum^{l-1}_{j=1}\int\langle
J(u)\nabla_{\bg}\tau(u),\nabla_{\a_2}\cdots\nabla_{\a_j}[R(\nabla_{\a_1}u,
\nabla_{\a_{j+1}}u)
\nabla_{\a_{j+2}}\cdots\nabla_{\a_l}\nabla_{\a_1}u]\rangle \nonumber \\
&&-\sum^{l-2}_{j=2}\int\langle
J(u)\nabla_{\bg}\tau(u),\nabla_{\a_1}\nabla_{\a_2}\cdots\nabla_{\a_j}
[R(\nabla_{\a_1}u, \nabla_{\a_{j+1}}u)\nabla_{\a_{j+2}}\cdots \nabla_{\a_l}
u]\rangle \nonumber \\
& = & -\sum^{l-1}_{j=1}\int\langle
J(u)\nabla_{\a_3}\cdots\nabla_{\a_l}\tau(u),\nonumber \\
&&\qquad \nabla_{\a_2}\nabla_{\a_2}\nabla_{\a_3}\cdots\nabla_{\a_j} 
 [R(\nabla_{\a_1}u, \nabla_{\a_{j+1}}u)\nabla_{\a_{j+2}} \cdots \nabla_{\a_l}
\nabla_{\a_1}u]\rangle \nonumber \\
&&-\sum^{l-2}_{j=2} \int\langle J(u)\nabla_{\a_3} \cdots
\nabla_{\a_l}\tau(u), \nabla_{\a_2}\nabla_{\a_1}\nabla_{\a_2} \cdots
\nonumber \\
&&\qquad \cdots \nabla_{\a_j} [R (\nabla_{\a_1}u, \nabla_{\a_{j+1}}u)\nabla_{\a_{j+2}} \cdots\nabla_{\a_l}
u]\rangle.\nonumber
\end{eqnarray}

Thus (\ref{in-eqn5.10}) yields
\begin{eqnarray}\label{in-eqn5.11}
&&\kern-.5in\sum_{|\ba|=l}\int\langle\nabla_{\ba}J(u)\tau(u),\nabla_{\ba}u\rangle \\
& \le &
C\sum^{l+2}_{m=3}\sum_{\mathop{j_1+\cdots+j_m=l+2}\limits_{\scriptscriptstyle{j_s\ge
1}}}\int|\nabla^{l-2}\tau(u)|\,|\nabla^{j_1}u|\cdots|\nabla^{j_m}u|.\nonumber
\end{eqnarray}
A very similar computation yields
\begin{eqnarray}\label{in-eqn5.12}
&&\kern.3in\int\langle\nabla_{\ba}\tau(u),\nabla_{\ba}u\rangle \\ 
&\le & -
\int|\nabla_{\bg}\tau(u)|^2 + C
\sum^{l+2}_{m=3}
\sum_{\mathop{j_1+\cdots+j_m=l+2}\limits_{\scriptscriptstyle{j_s\ge
1}}}\int|\nabla^{l-2}\tau(u)|\,|\nabla^{j_1}u|\cdots|\nabla^{j_m}u|.\nonumber
\end{eqnarray}
Combining (\ref{in-eqn5.6}), (\ref{in-eqn5.8}), (\ref{in-eqn5.9}),
(\ref{in-eqn5.11}) and (\ref{in-eqn5.12}) we obtain\newpage
\begin{eqnarray}\label{in-eqn5.13}
&&\sum_{|\ba|=l}\int\langle\nabla_{\ba}\nabla_tu, \nabla_{\ba}u\rangle \\
& \le &
-\e\int|\nabla^l\tau(u)|^2+C\e\sum^{l+2}_{m=3}
\sum_{j_1+\cdots+j_m=l+2}\int|\nabla^l\tau(u)|\,|\nabla^{j_1}u|\cdots|
\nabla^{j_m}u| \nonumber\\
&&+C\e\sum^{l+2}_{m=3}\sum_{j_1+\cdots+j_m=l+2}\int|\nabla^{j_1}\tau(u)|\,|
\nabla^{j_2}u|\cdots|\nabla^{j_m}u|\,|\nabla^l u|
\nonumber \\
&&+ C\e
\sum^{l+2}_{m=3}\sum_{\mathop{j_1+\cdots+j_m=l+2}\limits_{\scriptscriptstyle{j_s\ge
1\mbox{\tiny{ if }} s\ge 2}}}
\int|\nabla^lu|\,|\nabla^{j_1}\tau(u)|\,|\nabla^{j_2}u|\cdots|\nabla^{j_m}u|
\nonumber \\
&& - \beta\int|\nabla^{l-1}\tau(u)|^2+C\sum^{l+2}_{m=3}
\sum_{j_1+\cdots+j_m=l+2}
\int|\nabla^{l-2}\tau(u)|\,|\nabla^{j_1}u|\cdots|\nabla^{j_m}u|. \nonumber
\end{eqnarray}
We now look at the second term in (\ref{in-eqn5.5}). Using equation
(\ref{in-eqn5.1}) we obtain
\begin{eqnarray}\label{in-eqn5.14}
&&\sum_{|\ba|=l}\sum^{l-2}_{j=0}\int\langle\nabla_{\a_1}\cdots\nabla_{\a_j}
[R (\nabla_tu, \nabla_{\a_{j+1}}u)\nabla_{\a_{j+2}}\cdots\nabla_{\a_l}u],
\nabla_{\ba}u\rangle \\
& \le & C \sum^{l+1}_{m=3}
\sum_{\mathop{j_1+\cdots+j_m=l}\limits_{\scriptscriptstyle{j_s\ge 1\mbox{
\tiny{ if }} s\ge 2}}}
\int|\nabla^lu|\,|\nabla^{j_1}\nabla_tu|\,|\nabla^{j_2}u|\cdots|\nabla^{j_m}u|
\nonumber \\
& \le &
C\e\sum^{l+1}_{m=3}\sum_{\mathop{j_1+\cdots+j_m=l}\limits_{\scriptscriptstyle{j_s\ge
1\mbox{ 
\tiny{ if }} s\ge 2}}}
\int|\nabla^lu|\,|\nabla^{j_1}\D\tau(u)|\,|\nabla^{j_2}u|\cdots|\nabla^{j_m}u|
\nonumber \\
&& + C\e\sum^{l+3}_{m=6}
\sum_{\mathop{j_1+\cdots+j_m=l+2}\limits_{\scriptscriptstyle{j_s\ge 1\mbox{ 
\tiny{ if }} s\ge 2}}}
\int|\nabla^lu|\,|\nabla^{j_1}\tau(u)|\,|\nabla^{j_2}u|\cdots|\nabla^{j_m}u|
\nonumber \\
&& +
C\sum^l_{m=3}\sum_{\mathop{j_1+\cdots+j_m=l}\limits_{\scriptscriptstyle{j_s\ge
1\mbox{ 
\tiny{ if }} s\ge 2}}}
\int|\nabla^lu|\,|\nabla^{j_1}\tau(u)|\,|\nabla^{j_2}u|
\cdots|\nabla^{j_m}u|. \nonumber
\end{eqnarray}

Combining (\ref{in-eqn5.5}), (\ref{in-eqn5.13}) and (\ref{in-eqn5.14}) we
obtain
\newpage
\begin{eqnarray}\label{in-eqn5.15}
&&\kern.4in\frac{1}{2}\frac{d}{dt}\|\nabla^lu\|^2_{L^2}\\
 & \le &\kern-.1in
-\e\int|\nabla^l\tau(u)|^2 \nonumber \\
&& + C\e\int|\nabla^l\tau(u)|\,|\nabla^lu|\,|\nabla u|^2 \nonumber\\
&& + C\e\int|\nabla^{l-1}\tau(u)|\,|\nabla^lu|\,\left(|\nabla u|^3 + |\nabla
u|\,|\nabla^2u|\right) \nonumber \\
&& + C\e\sum^{l+2}_{m=3}\;
\sum_{\mathop{j_1+\cdots+j_m=l+2}\limits_{\scriptscriptstyle{1\le
j_s\le l-1}}}\int|\nabla^l\tau(u)|\,|\nabla^{j_1}u|\cdots|\nabla^{j_m}u|
\nonumber \\
&&
C\e\sum^{l+2}_{m=3}\;\sum_{\mathop{j_1+\cdots+j_m=l+2}\limits_{{\scriptscriptstyle{j_s\ge
1\mbox{ \tiny{ if }} s\ge 2}}}} \int|\nabla^l
u|\,|\nabla^{j_1}\tau(u)|\cdots|\nabla^{j_m}u| \nonumber \\
&&
+C\e\sum^{l+2}_{m=3}\;\sum_{\mathop{j_1+\cdots+j_m=l+2}\limits_{%
{\scriptscriptstyle{j_s\ge 1}}}}
\int|\nabla^lu|\,|\nabla^{j_1}u|\cdots|\nabla^{j_m}u| \nonumber \\
&&
+C\e\sum^{l+3}_{m=5}\sum_{\mathop{j_1+\cdots+j_m=l+4}\limits_{%
\scriptscriptstyle{j_s\ge 1}}}
\int|\nabla^lu|\,|\nabla^{j_1}u|\cdots|\nabla^{j_m}u| \nonumber \\
&& +C\sum^{l+2}_{m=3}\;\sum_{j_1+\cdots+j_m=l+2}
\int|\nabla^lu|\,|\nabla^{j_1}u|\cdots|\nabla^{j_m}u|
-\beta\int|\nabla^{l-1}\tau(u)|^2 \nonumber \\
&\le &\kern-.1in -\e\int|\nabla^l\tau(u)|^2 +
C\e\int|\nabla^l\tau(u)|\,|\nabla^lu|\,|\nabla u|^2 \nonumber \\
&&+C\e\int|\nabla^{l-1}\tau(u)|\,|\nabla^lu|\,\left(|\nabla u|^3 +
|\nabla^2u|\,|\nabla u|\right)\nonumber\\
&&+ C\e\int|\nabla^lu|\,|\tau(u)|\,|\nabla^{l+1}u|\,|\nabla u| \nonumber \\
&&+C\e\int|\nabla^lu|^2\left(|\tau(u)|\,|\nabla u|^2 +
|\tau(u)|\,|\nabla^2u|\right)\nonumber\\
&&+ C\e\int|\nabla^lu|^2\left(|\nabla
u|^4 + |\nabla\tau(u)|\,|\nabla u|\right)\nonumber \\
&&+C\e\sum^{l+2}_{m=3}\;\sum_{\mathop{\scriptscriptstyle{j_1+\cdots+j_m=l+2}}%
\limits_{\scriptscriptstyle{1\le
js\le l-1}}} \int|\nabla^l\tau(u)|\,|\nabla^{j_1}u|\cdots|\nabla^{j_m}u|
\nonumber \\
&&+C\e\sum^{l+3}_{m=5}\;\sum_{\mathop{\scriptscriptstyle{j_1+\cdots+j_m=l+4}}\limits_{\scriptscriptstyle{1\le
j_s\le l-1}}} \int|\nabla^lu|\,|\nabla^{j_1}u|\cdots|\nabla^{j_m}u|
\nonumber \\
&&+C\sum^{l+2}_{m=3}\; \sum_{j_1+\cdots+j_m=l+2}
\int|\nabla^lu|\,|\nabla^{j_1}u|\cdots|\nabla^{j_m}u| \nonumber \\
&&+C\e\sum^{l+2}_{m=3}\;
\sum_{\mathop{j_1+\cdots+j_m=l+2}\limits_{\mathop{\scriptscriptstyle{1\le j_s
\le l-1,s\ge 2}}\limits_{\scriptscriptstyle{j_1\le l-2}}}} 
\int_l|\nabla^lu|\,|\nabla^{j_1}\tau(u)|\cdots|\nabla^{j_m}u|.
\nonumber
\end{eqnarray}
We now look at each term of (\ref{in-eqn5.15}) separately. Apply
Cauchy-Schwarz we have
\begin{eqnarray}\label{in-eqn5.16}
\kern.2in 
C\e\int|\nabla^l\tau(u)|\,|\nabla^lu|\,|\nabla u|^2 & \le & C\e\|\nabla
u\|^2_{L^\infty}(\int|\nabla^l\tau(u)|^2)^{\frac{1}{2}}
(\int|\nabla^lu|^2)^{\frac{1}{2}} \\
& \le & \frac{\e}{64} \int|\nabla^l\tau(u)|^2+C\|\nabla
u\|^4_{L^\infty}\int|\nabla^lu|^2. \nonumber
\end{eqnarray}

Using Cauchy-Schwarz and integration by parts we have
\begin{eqnarray}\label{in-eqn5.17}
&&\kern-.5in C\e\int|\nabla^{l-1}\tau(u)|\,|\nabla^lu|\,|\nabla u|^3 \\
& \le & C\e\|\nabla
u\|^3_{L^\infty}(\int|\nabla^{l-1}\tau(u)|^2)^{\frac{1}{2}}
(\int|\nabla^lu|^2)^{\frac{1}{2}}\nonumber \\
&\le & \frac{\e}{64} \int|\nabla^{l-1}\tau(u)|^2+C\|\nabla u\|^6_{L^\infty}
\int|\nabla^lu|^2 \nonumber \\
& \le & \frac{\e}{64}\int|\nabla^l\tau(u)|\,|\nabla^{l-2}\tau(u)|+C\|\nabla
u\|^6_{L^\infty}\int|\nabla^lu|^2 \nonumber \\
& \le & \frac{\e}{64}(\int|\nabla^l\tau(u)|^2)^{\frac{1}{2}}
(\int|\nabla^lu|^2)^{\frac{1}{2}}+C\|\nabla
u\|^6_{L^\infty}\int|\nabla^lu|^2 \nonumber \\
& \le & \frac{\e}{64} \int|\nabla^l\tau(u)|^2 + C(1+\|\nabla
u\|^6_{L^\infty}) \int|\nabla^lu|^2. \nonumber
\end{eqnarray}
\begin{eqnarray}\label{in-eqn5.29a}
&&\kern-.3in C\e\int |\nabla^{l-1}\tau(u)|\,|\nabla^lu|\,|\nabla^2u|\,|\nabla u|\\
 & \le &
C\e\|\nabla u\|_{L^\infty} \|\nabla^2u\|_{L^\infty} 
\left(\int|\nabla^{l-1}\tau(u)|^2\right)^{\frac{1}{2}}
\left(\int|\nabla^lu|^2\right)^{\frac{1}{2}} \nonumber\\
& \le & \frac{\e}{64}\int|\nabla^l\tau(u)|^2 
+ C\left(1+\|\nabla
u\|^6_{L^\infty} + \|\nabla^2 u\|^3_{L^\infty}\right)\,\int|\nabla^l
u|^2. \nonumber
\end{eqnarray}
Using Cauchy-Schwarz, integration by parts and (\ref{in-eqn5.4}) we have
\begin{eqnarray}\label{in-eqn5.18}
&&\kern-.5in C\e\int|\nabla^lu|\,|\tau(u)|\,|\nabla^{l+1}u|\,|\nabla u| \\
& \le &
C\e\|\tau(u)\|_{L^\infty}\|\nabla
u\|_{L^\infty}(\int|\nabla^lu|^2)^{\frac{1}{2}}(\int|\nabla^{l+1}u|^2)^{\frac{1}{2}}
\nonumber\\
& \le & \frac{\e}{64}\int|\nabla^{l+1}u|^2 +
C\|\tau(u)\|^2_{L^\infty}\|\nabla u\|^2_{L^\infty}\int|\nabla^lu|^2
\nonumber \\
& \le & -\frac{\e}{64} \int\langle\Delta\nabla^lu, \nabla^lu\rangle +
C\|\tau(u)\|^2_{L^\infty}\|\nabla u\|^2_{L^\infty} \int|\nabla^lu|^2
\nonumber \\
& \le &
\frac{\e}{64}\int|\nabla^l\tau(u)|\,|\nabla^lu|+C\|\tau(u)\|^2_{L^\infty}
\|\nabla u\|^2_{L^\infty}\int|\nabla^lu|^2 \nonumber \\
&&+C\e\int|\nabla^lu|\sum^{l+2}_{m=3} \;\sum_{j_1+\cdots+j_m=l+2}
|\nabla^{j_1}u|\cdots|\nabla^{j_m}u| \nonumber \\
& \le & \frac{\e}{64}\int|\nabla^l\tau(u)|^2 + C(\|\tau(u)\|^2_{L^\infty}
\|\nabla u\|^2_{L^\infty} +1) \int|\nabla^lu|^2 \nonumber \\
&&+C\e\sum^{l+2}_{m=3}
\;\sum_{j_1+\cdots+j_m=l+2}\int|\nabla^lu|\,|\nabla^{j_{1}¥}u|\cdots|\nabla^{j_m}u|.
\nonumber 
\end{eqnarray}

Combining (\ref{in-eqn5.15}), (\ref{in-eqn5.16}), (\ref{in-eqn5.17}),
(\ref{in-eqn5.29a}) and
(\ref{in-eqn5.18}) and using the fact that $ab\le
\frac{a^p}{p}+\frac{b^q}{q}$ if $\frac{1}{p}+\frac{1}{q}=1$ we have
\begin{eqnarray}\label{in-eqn5.19}
\kern.2in\frac{1}{2}&\frac{d}{dt}&\|\nabla^lu\|^2_{L^2}  \le 
\frac{-3\e}{4}\int|\nabla^l\tau(u)|^2 \\
&&+C(1+\|\nabla u\|^6_{L^\infty} + \|\nabla^2u\|^3_{L^\infty}+\|\nabla u\|_{L^\infty}\|\nabla^3 u\|_{L^\infty})
\int|\nabla^lu|^2 \nonumber \\
&&+C\e\sum^{l+2}_{m=3}\sum_{\mathop{\scriptscriptstyle{j_1+\cdots+j_m=l+2}}\limits_{\scriptscriptstyle{1\le
j_s\le l-1}}} \int|\nabla^l\tau(u)|\,|\nabla^{j_1}u|\cdots|\nabla^{j_m}u|
\nonumber \\
&& +
C\e\sum^{l+1}_{m=5}\sum_{\mathop{\scriptscriptstyle{j_1+\cdots+j_m=l+4}}\limits_{\scriptscriptstyle{1\le
j_s\le l-1}}} \int|\nabla^lu|\,|\nabla^{j_1}u|\cdots|\nabla^{j_m}u|
\nonumber \\
&& + C\e\sum^{l+2}_{m=3}
\sum_{\mathop{\scriptscriptstyle{j_1+\cdots+j_m=l+2}}%
\limits_{\mathop{\scriptscriptstyle{l-1\ge
j_s\ge 1\mbox{ {\tiny{if}} }s\ge 2}}\limits_{\scriptscriptstyle{j_1\le
l-2}}}} \int|\nabla^lu|\,|\nabla^{j_1}\tau(u)|\cdots|\nabla^{j_m}u|
\nonumber \\
&& +C\sum^{l+2}_{m=3}
\sum_{\mathop{j_1+\cdots+j_m=l+2}\limits_{\scriptscriptstyle{1\le j_s
\le l-1}}}
\int|\nabla^lu|\,|\nabla^{j_1}u|\cdots|\nabla^{j_m}u|.\nonumber
\end{eqnarray}
To finish the estimate we need to use the interpolation result that appears
in Proposition \ref{an-prop1}.
Consider $3\le m\le l+2$, $1\le j_s\le l-1$ and $j_1+\cdots+j_m=l+2$ then
by Cauchy-Schwarz we have
\begin{eqnarray}\label{in-eqn5.20}
&&\kern-.8in\int|\nabla^l\tau(u)|\,|\nabla^{j_1}u|\cdots|\nabla^{j_m}u| \\
&\le & \left(\int|\nabla^l\tau(u)|^2\right)^{\frac{1}{2}}\left(\int|
\nabla^{j_1}u|^2\cdots|\nabla^{j_m}u|^2\right)^{\frac{1}{2}}.\nonumber
\end{eqnarray}
Let $p_i\in[2,\infty]$ for $i=1,\ldots, m$ be such that 
\begin{equation}\label{in-eqn5.21}
\frac{1}{p_1}+\cdots+\frac{1}{p_m}=\frac{1}{2},
\end{equation}
by H\"older's inequality
\begin{equation}\label{in-eqn5.22}
\left(\int|\nabla^{j_1}u|^2\cdots|\nabla^{j_m}u|^2\right)^{\frac{1}{2}}\le\|\nabla^{j_1}u\|_{L^{p_1}}\cdots
\|\nabla^{j_m}u\|_{L^{p_m}}.
\end{equation}
Since $l>\left[\frac{n}{2}\right]+1$ for 
\begin{equation}\label{in-eqn5.23}
\frac{j_i-1}{l-1}<a_i = \frac{j_i-1}{l-1}+\frac{n}{2(l-1)^2}
\left(l-1-j_i+\frac{3}{m}\right)<1
\end{equation}
and when $m>3$ or $m=3$ and $j_i\ge 2$
\begin{equation}\label{in-eqn5.24}
\frac{1}{2}\ge \frac{1}{p_i}=\frac{j_i-1}{n}+\frac{1}{2} -
\frac{l-1}{n}a_i>0.
\end{equation}
Thus (\ref{an-eqn5}) yields
\begin{equation}\label{in-eqn5.25}
\|\nabla^{j_i}u\|_{L^{p_i}}  \le  C\|\nabla^lu\|^{a_i}_{L^2}\|\nabla
u\|^{1-a_i}_{L^2} 
 \le  C\|\nabla u\|_{H^{l-1}}.
\end{equation}

Combining (\ref{in-eqn5.20}), (\ref{in-eqn5.22}) and (\ref{in-eqn5.25}) we
have in the case $m>3$ that
\begin{equation}\label{in-eqn5.26}
\int|\nabla^l\tau(u)|\,|\nabla^{j_1}u|\cdots|\nabla^{j_m}u|  \le
c\left(\int|\nabla^l\tau(u)|^2\right)^{\frac{1}{2}}\|\nabla u\|^m_{H^{l-1}}.
\end{equation}
In the case when $m=3$, $j_1\ge j_2\ge j_3$, and $j_3=1$ we have
$j_1+j_2=l+1$ and
\begin{equation}\label{in-eqn5.27}
\left(\int|\nabla^{j_1}u|^2|\nabla^{j_2}u|^2|\nabla
u|^2\right)^{\frac{1}{2}} \le \|\nabla
u\|_{L^\infty}\left(\int|\nabla^{j_1}u|^2|\nabla^{j_2}u|^2\right)^{\frac{1}{2}}.
\end{equation}

If $j_2=1$ then (\ref{in-eqn5.20}) becomes 
\begin{equation}\label{in-eqn5.28}
\int|\nabla^l\tau(u)|\,|\nabla^lu|\,|\nabla u|^2 \le c\|\nabla
u\|^2_{L^\infty} \left(\int|\nabla^l\tau(u)|^2\right)^{\frac{1}{2}}\left(
\int|\nabla^lu|^2\right)^{\frac{1}{2}}.
\end{equation}
If $j_2>1$ then for $i=1,2$ let 
$l_0=\max\left\{\left[\frac{n}{2}\right]+4,l\right\}$. If
\begin{equation}\label{in-eqn5.29}
\frac{j_i-1}{l_0-1}\le a_i = \frac{j_i-1}{l_0-1}+\frac{n}{2(l_0-1)(l-1)} 
(l-j_i)<1
\end{equation}
and
\begin{equation}\label{in-eqn5.30}
\frac{1}{2}\ge \frac{1}{p_i}=\frac{j_i-1}{n} + \frac{1}{2} -
\frac{l_0-1}{n}a_i>0.
\end{equation}
H\"older's inequality and Proposition \ref{an-prop1} yield
\begin{eqnarray}\label{in-eqn5.31}
\left(\int|\nabla^{j_1}u|^2|\nabla^{j_2}u|^2\right)^{\frac{1}{2}} & \le &
\|\nabla^{j_1}u\|_{L^{p_1}}\|\nabla^{j_2}u\|_{L^{p_2}} \\
& \le & C\|\nabla^l_0u\|^{a_1}_{L^2}\|\nabla u\|^{1-a_1}_{L^2}
\|\nabla^l_0u\|^{a_2}_{L^2} \|\nabla u\|^{1-a_2}_{L^2} \nonumber \\
& \le & C\|\nabla u\|^2_{H^{l_0-1}}. \nonumber
\end{eqnarray}
Thus in this case (\ref{in-eqn5.20}) becomes combining (\ref{in-eqn5.27})
and (\ref{in-eqn5.31})
\begin{equation}\label{in-eqn5.32}
\int|\nabla^l\tau(u)|\,|\nabla^{j_1}u|\cdots|\nabla u|\le c\|\nabla
u\|_{L^\infty}\|\nabla
u\|^2_{H^{l_0-1}}\left(\int|\nabla^l\tau(u)|^2\right)^{\frac{1}{2}}.
\end{equation}

Combining (\ref{in-eqn5.26}), (\ref{in-eqn5.28}) and (\ref{in-eqn5.32}) we
can estimate the third term on the right hand side of (\ref{in-eqn5.19})
\begin{eqnarray}\label{in-eqn5.33}
&&\sum^{l+2}_{m=3}\sum_{\mathop{\scriptscriptstyle{j_1+\cdots+j_m=l+2}}\limits_{\scriptscriptstyle{1\le
j_s\le l-1}}}\int|\nabla^l\tau(u)|\,|\nabla^{j_1}u|\cdots|\nabla^{j_m}u|\\
&&\kern.2in\le
C\left(\int|\nabla^l\tau(u)|^2\right)^{\frac{1}{2}}\left(1+\|\nabla
u\|^2_{L^\infty}\right)\left(\sum^{l+2}_{m=1}\|\nabla u\|^m_{H^{l-1}}
+\|\nabla u\|^2_{H^{l_0-1}}\right).\nonumber
\end{eqnarray}
To estimate the fourth term in (\ref{in-eqn5.19}) consider $5\le m\le l+3$,
$1\le j_s\le l-1$ and $j_1+\cdots+j_m=l+4$ then by Cauchy-Schwarz we have
\begin{equation}\label{in-eqn5.34}
\int|\nabla^lu|\,|\nabla^{j_1}u|\cdots|\nabla^{j_m}u|\le
\left(\int|\nabla^lu|^2\right)^{\frac{1}{2}}\left(\int|\nabla^{j_1}u|^2\cdots|\nabla^{j_m}u|^2\right)^{\frac{1}{2}}.
\end{equation}

Let $2\le p_i\le \infty$ for $i=1,\ldots, m$ be such that 
\begin{equation}\label{in-eqn5.35}
\frac{1}{p_1}+\cdots+\frac{1}{p_m}=\frac{1}{2}
\end{equation}
by H\"older's inequality (\ref{in-eqn5.34}) becomes
\begin{equation}\label{in-eqn5.36}
\int|\nabla^lu|\,|\nabla^{j_1}u|\cdots|\nabla^{j_m}u|\le
\left(\int|\nabla^lu|^2\right)^{\frac{1}{2}}\|\nabla^{j_1}u\|_{L^{p_1}}
\cdots \|\nabla^{j_m}u\|_{L^{p_m}}
\end{equation}
since $l>\left[\frac{n}{2}\right]+1$ for
\begin{equation}\label{in-eqn5.37}
\frac{j_i-1}{l-1}\le a_i=\frac{j_i-1}{l-1} +
\frac{n}{2(l-1)^2}\left(l-1-j_i+\frac{5}{m}\right)<1
\end{equation}
and when $m>5$ or $m=5$ and $j_i\ge 2$
\begin{equation}\label{in-eqn5.38}
\frac{1}{2}\ge \frac{1}{p_i}=\frac{j_i-1}{n}+\frac{1}{2}-\frac{l-1}{n}a_i>0.
\end{equation}
Thus (\ref{an-eqn5}) yields
\begin{equation}\label{in-eqn5.39}
\int|\nabla^lu\|\nabla^{j_1}u|\cdots|\nabla^{j_m}u|\le
C\left(\int|\nabla^lu|^2\right)^{\frac{1}{2}}\|\nabla u\|^m_{H^{l-1}}.
\end{equation}

If $m=5$, $j_1\ge j_2\ge \cdots\ge j_5\ge 1$, and $j_5=1$ then
$j_1+j_2+j_3+j_4=l+3$, by Cauchy-Schwarz and H\"older's inequality
\begin{eqnarray}\label{in-eqn5.40}
&&\kern-.5in\int|\nabla^lu|\,|\nabla^{j_1}u|\cdots|\nabla^{j_4}u|\,|\nabla
u| \\
& \le &
\|\nabla u\|_{L^\infty} \left(\int|\nabla^lu|^2\right)^{\frac{1}{2}}
\left(\int|\nabla^{j_1}u|^2\cdots|\nabla^{j_4}u|^2\right)^{\frac{1}{2}} 
\nonumber \\
& \le & \|\nabla
u\|_{L^\infty}\left(\int|\nabla^lu|^2\right)^{\frac{1}{2}}\|\nabla^{j_1}u\|_{L^{p_1}}\cdots\|\nabla^{j_4}u\|_{L^{p_4}}
\nonumber
\end{eqnarray}
with $\frac{1}{p_1}+\frac{1}{p_2}+\frac{1}{p_3}+\frac{1}{p_4}=\frac{1}{2}$. For
\begin{equation}\label{in-eqn5.41}
\frac{j_i-1}{l-1} < a_i = \frac{j_i-1}{l-1}+\frac{n}{2(l-1)^2} (l-j_i)<1
\end{equation}
if $j_4>1$ we have
\begin{equation}\label{in-eqn5.42}
\frac{1}{2}\ge \frac{1}{p_i} = \frac{j_i-1}{n} + \frac{1}{2} - \frac{l-1}{n}
a_i>0
\end{equation}
and (\ref{in-eqn5.40}) becomes by Proposition \ref{an-prop1}
\begin{equation}\label{in-eqn5.41A}
\int|\nabla^lu|\,|\nabla^{j_1}u|\cdots|\nabla^{j_4}u|\,|\nabla u|\le
C\|\nabla u\|_{L^\infty} \left(\int|\nabla^lu|^2\right)^{\frac{1}{2}}
\|\nabla u\|^4_{H^{l-1}}.
\end{equation}

If $j_4=1$ and $j_3>1$ a similar argument yields 
\begin{equation}\label{in-eqn5.43}
\int|\nabla^lu|\,|\nabla^{j_1}u|\,|\nabla^{j_3}u|\,|\nabla u|^2 \le
C\|\nabla
u\|^2_{L^\infty}\left(\int|\nabla^lu|^2\right)^{\frac{1}{2}}\|\nabla
u\|^3_{H^{l-1}}.
\end{equation}
If $j_3=1$ then $j_1+j_2=l+1$ since $j_1\le l-1$ then $j_2>1$ and we have
\begin{equation}\label{in-eqn5.44}
\int|\nabla^lu|\,|\nabla^{j_1}u|\,|\nabla^{j_2}u|\,|\nabla u|^3 \le
C\|\nabla u\|^3_{L^\infty} \left(\int|\nabla^lu|^2\right)^{\frac{1}{2}}
\|\nabla u\|^2_{H^{l-1}}.
\end{equation}

Combining (\ref{in-eqn5.39}), (\ref{in-eqn5.41A}), (\ref{in-eqn5.43}) and
(\ref{in-eqn5.44}) we can estimate the fourth term on the right hand side of
(\ref{in-eqn5.19}) as follows
\begin{eqnarray}\label{in-eqn5.45}
&&\kern-.5in\sum^{l+3}_{m=5}\;\sum_{\mathop{\scriptscriptstyle{j_1+
\cdots+j_m=l+4}}\limits_{\scriptscriptstyle{1\le
j_s\le l-1}}} \int|\nabla^l u|\,|\nabla^{j_1}u|\cdots|\nabla^{j_m}u| \\
& \le &
C\left(\int|\nabla^lu|^2\right)^{\frac{1}{2}}(1+\|\nabla u\|^3_{L^\infty})
\sum^{l+3}_{m=2} \|\nabla u\|^m_{H^{l-1}} \nonumber\\
& \le & C\left(1+\|\nabla u\|^3_{L^\infty}\right) \sum^{l+4}_{m=3} \|\nabla
u\|^m_{H^{l-1}}. \nonumber
\end{eqnarray}

To estimate the fifth term in (\ref{in-eqn5.19}) consider $3\le m\le l+2$,
$j_1+\cdots+j_m=l+2$, $j_1\le l-2$, $1\le j_s\le l-1$ if $s\ge 2$.
Cauchy-Schwarz and H\"older's inequality ensure that for
$\frac{1}{p_1}+\cdots+\frac{1}{p_m}=\frac{1}{2}$
\begin{eqnarray}\label{in-eqn46}
&&\kern-.5in\int|\nabla^lu|\,|\nabla^{j_1}\tau(u)|\,|\nabla^{j_2}u|\cdots|\nabla^{j_m}u|
\\
&\le&
\left(\int|\nabla^lu|^2\right)^{\frac{1}{2}}\|\nabla^{j_1}\tau(u)\|_{L^{p_1}}
\cdots \|\nabla^{j_m}u\|_{L^{p_m}}.\nonumber
\end{eqnarray}
For $l_0=\max\left\{\left[\frac{n}{2}\right]+4,l\right\}>1$ for $i\ge 2$
\begin{align}
\frac{j_i-1}{l_0-1}<a_i&=\frac{j_i-1}{l_0-1}+\frac{n}{2(l_0-1)(l-1)}
\left(l-1-j_i
+\frac{3}{m}\right)<1 \label{in-eqn47} \\
\intertext{and}
\frac{j_1}{l_0-1}<a_1 &=\frac{j_1}{l_0-1}+\frac{n}{2(l-1)(l_0-1)}
\left(l-1-j_1+\frac{3}{m}\right)<1 \label{in-eqn48}
\end{align}
when $m>3$ or $m=3$ and $j_i\ge 2$ for $i\ge 2$
\begin{equation}\label{in-eqn49}
\frac{1}{2}\ge \frac{1}{p_i} =\frac{j_i-1}{n}+\frac{1}{2}-\frac{l_0-1}{n}a_i>0
\end{equation}
and $m>3$ or $m=3$ and $j_1\ge 2$
\begin{equation}\label{in-eqn50}
\frac{1}{2}\ge \frac{1}{p_1} = \frac{j_1}{n}+\frac{1}{2}-\frac{l_0-1}{n}a_1>0.
\end{equation}

In these cases (\ref{in-eqn46}) can be estimated by (\ref{an-eqn5}) as
follows
\begin{eqnarray}\label{in-eqn51}
&&\kern-.4in\int|\nabla^lu|\,|\nabla^{j_1}\tau(u)|\cdots|\nabla^{j_m}u| \\
&\le&
\left(\int|\nabla^lu|^2\right)^{\frac{1}{2}}
\|\nabla^{l_0-1}\tau(u)\|^{a_1}_{L^2}
\|\tau(u)\|^{1-a_1}_{L^2}\|\nabla u\|^{m-1}_{H^{l_0-1}}.\nonumber
\end{eqnarray}
If $m=3$ and $j_1\le 1$ then $j_2\ge 2$ and $j_3\ge 2$. Cauchy-Schwarz and
H\"older's inequality yield for $\frac{1}{p_2}+\frac{1}{p_3}=\frac{1}{2}$
\begin{eqnarray}\label{in-eqn52}
&&\kern-.4in\int|\nabla^lu|\,|\tau(u)|\,|\nabla^{j_2}u|\,|\nabla^{j_3}u| \\
&\le&
\|\tau(u)\|_{L^\infty}
\left(\int|\nabla^lu|^2\right)^{\frac{1}{2}}\|\nabla^{j_2}u\|_{L^{p_2}}\|\nabla^{j_3}u\|_{L^{p_3}}.\nonumber
\end{eqnarray}
and
\begin{equation}\label{in-eqn54}
\frac{1}{2}\ge \frac{1}{p_i} = \frac{j_i-1}{n} + \frac{1}{2} -
\frac{l-1}{n}a_i>0,
\end{equation}
Proposition \ref{an-prop1} ensures that
\begin{equation}\label{in-eqn5.53}
\int|\nabla^lu|\,|\tau(u)|\,|\nabla^{j_2}u|\,|\nabla^{j_3}u|\le
C\|\tau(u)\|_{L^\infty} \left(\int|\nabla^lu|^2\right)^{\frac{1}{2}}
\|\nabla u\|^{2}_{H^{l-1}}.
\end{equation}
Similarly
\begin{equation}\label{in-eqn5.54}
\int|\nabla^lu|\,|\nabla\tau(u)|\,|\nabla^{j_2}u|\,|\nabla^{j_3}u|
\le
C\|\nabla\tau(u)\|_{L^\infty}
\left(\int|\nabla^lu|^2\right)^{\frac{1}{2}}\|\nabla
u\|^2_{H^{l-1}}.
\end{equation}
If $m=3$, $j_1\ge 2$ and $j_2=1$, $j_3>1$ we have by Cauchy-Schwarz and
H\"older's inequality for $\frac{1}{p_1}+\frac{1}{p_3}=\frac{1}{2}$
\begin{eqnarray}\label{in-eqn55}
&&\kern-.4in\int|\nabla^lu|\,|\nabla^{j_1}\tau(u)|\,|\nabla
u|\,|\nabla^{j_3}u|\\
&\le &
C\|\nabla u\|_{L^\infty} \left(\int|\nabla^lu|^2\right)^{\frac{1}{2}}
\|\nabla^{j_1}\tau(u)\|_{L^{p_1}}\|\nabla^{j_3}u\|_{L^{p_3}}.\nonumber
\end{eqnarray}
For
\begin{equation}\label{in-eqn56}
a_1=\frac{j_1}{l-1}  + \frac{n}{2(l-1)^2} (l-j_1)<1\mbox{ and }\frac{1}{2}
\ge \frac{1}{p_1}=\frac{j_1}{n}+\frac{1}{2}-\frac{l-1}{n}a_1>0
\end{equation}
and
\begin{equation}\label{in-eqn57}
a_3=\frac{j_3-1}{l-1} + \frac{n}{2(l-1)^2} (l-j_3)<1 \mbox{ and }
\frac{1}{2}\ge \frac{1}{p_3} = \frac{j_3-1}{n} + \frac{1}{2} -
\frac{l-1}{n}a_3>0.
\end{equation}
Proposition \ref{an-prop1} ensures
\begin{eqnarray}\label{in-eqn58}
&&\kern-.5in
\int|\nabla^lu|\,|\nabla^{j_1}\tau(u)|\,|\nabla u|\,|\nabla^{j_3}u| \\
&\le&
C\|\nabla u\|_{L^\infty}\left(\int|\nabla^lu|^2\right)^{\frac{1}{2}}
\|\nabla^{l-1}\tau(u)\|^{a_1}_{L^2}\|\tau(u)\|^{1-a_1}_{L^2} \|\nabla
u\|_{H^{l-1}}.\nonumber
\end{eqnarray}
In the case $j_1=l$, $j_2=j_3=1$ see (\ref{in-eqn5.16}).

Combining (\ref{in-eqn51}), (\ref{in-eqn5.53}), (\ref{in-eqn5.54}) and
(\ref{in-eqn58}) we estimate the $5^{\mathrm{th}}$ term of
(\ref{in-eqn5.19}) as follows
\begin{eqnarray}\label{in-eqn59}
&&\kern-.1in\sum^{l+2}_{m=3} \,
\sum_{\mathop{\scriptscriptstyle{j_1+\cdots+j_m=l+2}}\limits_{\mathop{\scriptscriptstyle{1\le
j_s\le l-1\mbox{ {\tiny{if}} }s\ge 2}}\limits_{\scriptscriptstyle{j_1\le
l-2}}}} \int|\nabla^lu|\,|\nabla^{j_1}\tau(u)|\cdots|\nabla^{j_m}u| \\
& \le & C\left(\int|\nabla^lu|^2\right)^{\frac{1}{2}}\left(1+\|\nabla
u\|_{L^\infty}\right) \sum^{l+1}_{m=1}\|\tau(u)\|_{L^2}
\|\nabla u\|^m_{H^{l-1}}
\|\nabla^{l-1}\tau(u)\|_{L^2}\nonumber\\
& & +  C\left(\int|\nabla^lu|^2\right)^{\frac{1}{2}}\left(1+\|\nabla
u\|_{L^\infty}\right) \sum^{l+1}_{m=1}\|\tau(u)\|_{L^2}
\|\nabla u\|^m_{H^{l_0-1}}
\|\nabla^{l_0-1}\tau(u)\|_{L^2}  \nonumber \\
&&+ C(\|\tau(u)\|_{L^\infty} +
\|\nabla\tau(u)
\|_{L^\infty})\left(\int|\nabla^lu|^2\right)^{\frac{1}{2}}\|\nabla
u\|^2_{H^{l-1}}. \nonumber
\end{eqnarray}
Here we have used the fact that for
$a\in(0,1)\le r^as^{1-a}\le ar+(1-a)s$.

Finally we look at the last term of (\ref{in-eqn5.19}). Let $3\le m\le l+2$,
$j_1+\cdots+j_m=l+2$. Applying the same argument as the one used to obtain
(\ref{in-eqn5.33}) we conclude that 
\begin{eqnarray}\label{in-eqn60}
&&\sum^{l+2}_{m=3}
\,\sum_{j_1+\cdots+j_l=l+2}
\int|\nabla^lu|\,|\nabla^{j_1}u|\cdots|\nabla^{j_m}u|\\
&&\qquad \le
C(1+\|\nabla
u\|^2_{L^\infty})\left(\int|\nabla^lu|^2\right)^{\frac{1}{2}}
\left(\sum^{l+2}_{m=1}
\|\nabla u\|^m_{H^{l-1}}+\|\nabla u\|^2_{H^{l_0-1}}\right).\nonumber
\end{eqnarray}

Combining (\ref{in-eqn5.19}), (\ref{in-eqn5.33}), (\ref{in-eqn5.45}),
(\ref{in-eqn59}) and (\ref{in-eqn60}), using (\ref{an-eqn62}),
(\ref{an-eqn64}), (\ref{an-eqn66}) and the fact that
$l>\left[\frac{n}{2}\right]+1$ as well as $\e\in(0,1)$ and the fact that 
for
$a\in(0,1)\le r^as^{1-a}\le ar+(1-a)s$ we obtain for 
$l_0=l\ge \left[\frac{n}{2}\right]+4$

\begin{eqnarray}\label{in-eqn61}
&&\kern.3in\frac{1}{2}\,\frac{d}{dt}
\|\nabla^lu\|^2_{L^2} \\
 & \le & -\frac{3\e}{4}\int|\nabla^l\tau(u)|^2 +
C(1+\|\nabla u\|^6_{L^\infty} +
\|\nabla^2u\|^3_{L^\infty})\|\nabla^lu\|^2_{L^2} \nonumber\\
&&+C\e\left(\int|\nabla^l\tau(u)|^2\right)^{\frac{1}{2}}(1+\|\nabla
u\|^2_{L^\infty})\sum^{l+2}_{m=1}\|\nabla u\|^{m}_{H^{l-1}} \nonumber \\
&&+C\e(1+\|\nabla
u\|^3_{L^\infty})\left(\int|\nabla^lu|^2\right)^{\frac{1}{2}}\sum^{l+3}_{m=2}\|\nabla
u\|^m_{H^{l-1}} \nonumber \\
&&+C\e\left(\int|\nabla^lu|^2\right)^{\frac{1}{2}} (1+\|\nabla
u\|_{L^\infty})(\|\nabla^{l-1}\tau(u)\|_{L^2}+
\|\tau(u)\|_{L^2}) \sum^{l+1}_{m=1}\|\nabla u\|^m_{H^{l-1}} \nonumber
\\
&&+C(\|\tau(u)\|_{L^\infty} +
\|\nabla\tau(u)\|_{L^\infty})\left(\int|\nabla^lu|^2\right)^{\frac{1}{2}}
\|\nabla u\|^2_{H^{l-1}} \nonumber \\
&&+C(1+\|\nabla
u\|^2_{L^\infty})\left(\int|\nabla^lu|^2\right)^{\frac{1}{2}}
\sum^{l+2}_{m=1}\|\nabla u\|^m_{H^{l-1}} \nonumber \\
& \le & -\frac{3\e}{4}\int|\nabla^l\tau(u)|^2+C\left(\|\nabla
u\|_{H^{\left[\frac{n}{2}\right]+4}}^{6\left[\frac{n}{2}\right]+12}+1\right)
\sum^{l+4}_{m=2}\|\nabla u\|^m_{H^{l-1}} \nonumber \\
&&+C\e\|\nabla^lu\|_{L^2}\left(1+\|\nabla
u\|^{\left[\frac{n}{2}\right]+2}_{H^{\left[\frac{n}{2}\right]+4}}\right)
\sum^{l+2}_{m=1}\|\nabla u\|^m_{H^{l-1}}\|\nabla^{l-1}\tau(u)\|_{L^2}
\nonumber \\
& \le & -\frac{3\e}{4}\int|\nabla^l\tau(u)|^2+\frac{\e}{64}\|\nabla^{l-1}
\tau(u) \|^2_{L^2} \nonumber \\
&&+C\left(1+\|\nabla u\|^{3n+12}_{H^{\left[\frac{n}{2}\right]+4}}\right)
\sum^{2l+4}_{m=2} \|\nabla u\|^m_{H^{l-1}}. \nonumber 
\end{eqnarray}

Using the same trick as in (\ref{in-eqn5.17}) we obtain from
(\ref{in-eqn61}) for $l\ge\left[\frac{n}{2}\right]+4$
\begin{eqnarray}\label{in-eqn62}
\kern.3in\frac{d}{dt}\|\nabla^lu\|^2_{L^2} & \le &
-\frac{\e}{2}\int|\nabla^l\tau(u)|^2 + C\left(1+\|\nabla
u\|^{3n+12}_{H^{\left[\frac{n}{2}\right]+4}}\right) \sum^{2l+4}_{m=2}\|\nabla
u\|^m_{H^{l-1}} \\
& \le & C\left(1+\|\nabla
u\|^{3n+12}_{H^{\left[\frac{n}{2}\right]+4}}\right)\|\nabla
u\|^2_{H^{l-1}}\left(1+\|\nabla u \|^{2l+2}_{H^{l-1}}\right) \nonumber
\end{eqnarray}

In the case where $l=\left[\frac{n}{2}\right]+2$ or 
$l=\left[\frac{n}{2}\right]+3$ then (\ref{in-eqn61}) and (\ref{in-eqn62}) 
become

\begin{eqnarray}\label{in-eqn61A}
&&\kern.2in\frac{1}{2}\,\frac{d}{dt}
\|\nabla^lu\|^2_{L^2} \\
& \le & -\frac{3\e}{4}\int|\nabla^l\tau(u)|^2+\frac{\e}{64}\|\nabla^{l-1}
\tau(u) \|^2_{L^2} 
+C\left(1+\|\nabla u\|^{3n+12}_{H^{\left[\frac{n}{2}\right]+4}}\right)
\sum^{2l+4}_{m=2} \|\nabla u\|^m_{H^{l-1}} \nonumber \\
&& +
C\e\|\nabla^lu\|_{L^2}\left(1+\|\nabla
u\|^{\left[\frac{n}{2}\right]+2}_{H^{\left[\frac{n}{2}\right]+4}}\right)
\sum^{l+2}_{m=1}\|\nabla u\|^m_{H^{l-1}}\|\nabla^{\left[\frac{n}{2}\right]+3}\tau(u)\|_{L^2}
\nonumber \\
&\le & C\left(1+\|\nabla
u\|^{3n+12}_{H^{\left[\frac{n}{2}\right]+4}}\right)\|\nabla
u\|^2_{H^{l-1}}\left(1+\|\nabla u\|^{2l+2}_{H^{l-1}}\right) \nonumber
\end{eqnarray}
Thus (\ref{in-eqn62}) and (\ref{in-eqn61A}) conclude the
proof of Lemma \ref{in-lem1}.
\end{proof}


Since our ultimate goal is to estimate $\frac{d}{dt}\|\nabla
u\|^2_{H^{l-1}}$ for $l\ge 1$, we still need to analyze
$\frac{d}{dt}\|\nabla^lu\|^2_{L^2}$ for $1\le l\le
\left[\frac{n}{2}\right]+1$.

\begin{lemma}\label{in-lem2}
Let $u\in C([0,T], H^{\left[\frac{n}{2}\right]+4}(\RR^n,N))$ be a solution
of (\ref{in-eqn1}). Let $1\le l\le \left[\frac{n}{2}\right]+1$ then if
$s_0=\left[\frac{n}{2}\right]+2$ we have 
\begin{equation}\label{in-eqn64}
\frac{d}{dt}\|\nabla^lu\|^2_{L^2}\le c\|\nabla
u\|^2_{H^{s_0}}\left(1+\|\nabla u\|^{M_l}_{H^{s_0+2}}\right).
\end{equation}
where $M_l=3n+2l+12$.
\end{lemma}

\begin{proof}
Note that (\ref{in-eqn5.4})  and the Sobolev embedding theorem yields
\begin{equation}\label{in-eqn65}
\frac{d}{dt}\|\nabla u\|^2_{L^2} \le C\|\nabla
u\|^4_{H^{\left[\frac{n}{2}\right]+1}}\left(1+\|\nabla
u\|^{2n+8}_{H^{\left[\frac{n}{2}\right]+2}}\right).
\end{equation}
Note that for $l\ge 2$ computation (\ref{in-eqn5.19}) remains valid. In fact
we only used $l>\left[\frac{n}{2}\right]+1$ when we started to interpolate
as in Proposition \ref{an-prop1}. Let $s_0=\left[\frac{n}{2}\right]+2$.
Consider $3\le m\le l+2$ $1\le j_s\le l-1$ and $j_1+\cdots+j_m=l+2$ then by
Cauchy-Schwarz, H\"older's inequality applied with
$\frac{1}{p_1}+\cdots+\frac{1}{p_m}=\frac{1}{2}$ where
\begin{equation}\label{in-eqn66}
\frac{1}{p_i}  = \frac{j_i-1}{n} + \frac{1}{2}-\frac{s_0}{n}a_i 
\end{equation}
and
\begin{equation}\label{in-eqn67}
\frac{j_i-1}{s_0}\le a_i=\frac{j_i-1}{s_0}+\frac{n}{2(l-1)s_0}
\left(l-1-j_i+\frac{3}{m}\right)<1
\end{equation}
and (\ref{an-eqn5}) in the case $m>3$ or $m=3$ and $j_i\ge 2$ we obtain as
in (\ref{in-eqn5.26})
\begin{equation}\label{in-eqn68}
\int|\nabla^l\tau(u)||\nabla^{j_1}u|\cdots|\nabla^{j_m}u|\le
C\left(\int|\nabla^l\tau(u)|^2\right)^{\frac{1}{2}}\|\nabla u\|^m_{H^{s_0}}.
\end{equation}
In the case $m=3$ we proceed as in the proof of (\ref{in-eqn5.32}) (where $s_0$
now plays the role of $l_0$) and obtain
\begin{equation}\label{in-eqn69}
\int|\nabla^l\tau(u)|^2|\nabla^{j_1}u|\cdots|\nabla^{j_m}u| \le C\|\nabla
u\|_{L^\infty}\left(\int|\nabla^l\tau(u)|^2\right)^{\frac{1}{2}} \|\nabla
u\|^2_{H^{s_0}}.
\end{equation}
Thus for $2\le l\le \left[\frac{n}{2}\right]+1$ (\ref{in-eqn5.33}) becomes
\begin{eqnarray}\label{in-eqn70}
&&\kern-.7in\sum^{l+2}_{m=3}
\sum_{\mathop{\scriptscriptstyle{j_1+\cdots+j_m=l+2}}\limits_{\scriptscriptstyle{1\le
j_s\le l-1}}} \int|\nabla^l\tau(u)|\,|\nabla^{j_1}u|\cdots|\nabla^{j_m}u|\\
&\le &
C\left(\int|\nabla^l\tau(u)|^2\right)^{\frac{1}{2}}\left(1+\|\nabla
u\|^2_{L^\infty}\right) \sum^{l+2}_{m=1}\|\nabla u\|^m_{H^{s_0}}. \nonumber
\end{eqnarray}
The same type of argument as the one used to prove (\ref{in-eqn5.45}),
(\ref{in-eqn59}) and (\ref{in-eqn60}) yields
\begin{eqnarray}
&&\kern-.2in\sum^{l+3}_{m=5}\sum_{\mathop{\scriptscriptstyle{j_1+\cdots+j_m=l+4}}\limits_{\mathop{\scriptscriptstyle{1\le
j_s\le l-1}}}}
\int|\nabla^lu|\,|\nabla^{j_1}u|\cdots|\nabla^{j_m}u| \label{in-eqn71}\\
& \le &
C\left(\int|\nabla^lu|^2\right)^{\frac{1}{2}}(1+\|\nabla
u\|^3_{L^\infty}+\|\nabla\tau(u)\|_{L^\infty})\|\nabla
u\|^m_{H^{s_0}} \nonumber \\
&&\kern-.2in\sum^{l+2}_{m=3}\sum_{\mathop{\scriptscriptstyle{j_1+\cdots+j_m=l+2}}\limits_{\mathop{\scriptscriptstyle{1\le
j_s\le l-1\quad s\ge 2}}\limits_{\scriptscriptstyle{j_1\le
l-2}}}}\int|\nabla^lu|\,|\nabla^{j_1}\tau(u)|\cdots|\nabla^{j_m}u|\label{in-eqn72}\\
 & \le &
C\left(\int|\nabla^lu|^2\right)^{\frac{1}{2}} (1+\|\nabla u\|_{L^\infty})
\sum^{l+1}_{m=1} \|\nabla u\|^m_{H^{s_0}}
\|\nabla^{s_0}\tau(u)\|^{a_1}_{L^2} \|\tau(u)\|^{1-a_1}_{L^2} \nonumber \\
&&+C(\|\tau(u)\|_{L^\infty} +
\|\nabla\tau(u)\|_{L^\infty})\left(\int|\nabla^lu|^2\right)^{\frac{1}{2}}
\|\nabla u\|^2_{H^{s_0}} \nonumber \\
& \le & C\left(\int|\nabla^lu|^2\right)^{\frac{1}{2}} (1+\|\nabla
u\|_{L^\infty}) \sum^{l+2}_{m=2} \|\nabla u\|^m_{H^{s_0+1}} \nonumber \\
&&+C(\|\tau(u)\|_{L^\infty} + \|\nabla\tau(u)\|_{L^\infty})
\left(\int|\nabla^lu|^2\right)^{\frac{1}{2}} \|\nabla u\|_{H^{s_0}}
\nonumber \\
&&\kern-.2in\sum^{l+2}_{m=3} \sum_{j_1+\cdots+j_l=l+2}
\int|\nabla^lu|^2|\nabla^{j_1}u|\cdots|\nabla^{j_m}u| \label{in-eqn73} \\
& \le & C(1+\|\nabla
u\|^2_{L^\infty})\left(\int|\nabla^lu|^2\right)^{\frac{1}{2}}
\sum^{l+2}_{m-1} \|\nabla u\|^m_{H^{s_0}}. \nonumber
\end{eqnarray}

Combining (\ref{in-eqn5.19}), (\ref{in-eqn70}), (\ref{in-eqn71}),
(\ref{in-eqn72}) and (\ref{in-eqn73}); using (\ref{an-eqn60}),
(\ref{an-eqn61}) and (\ref{an-eqn61A}) we have for $\e\in(0,1)$, $l\le
\left[\frac{n}{2}\right]+1$, $s_0=\left[\frac{n}{2}\right]+2$
\begin{eqnarray}\label{in-eqn74}
&&\kern-.3in\frac{1}{2}\,\frac{d}{dt} \|\nabla^lu\|^2_{L^2}\\
 & \le & -\frac{\e}{2}
\int|\nabla^l\tau(u)|^2 + C\left(1+\|\nabla u\|^6_{L^\infty} +
\|\tau(u)\|^3_{L^\infty}\right) \int|\nabla^lu|^2 \nonumber\\
&&+C\left(1+\|\nabla u\|^4_{L^\infty}\right) \sum^{2l+4}_{m=2} \|\nabla
u\|^m_{H^{s_0}} \nonumber \\
&&+C\left(1+\|\nabla u\|_{L^\infty} + \|\tau(u)\|_{L^\infty} +
\|\nabla\tau(u)\|_{L^\infty}\right)\sum^{l+3}_{m=2}\|\nabla
u\|^m_{H^{s_0+1}} \nonumber \\
& \le & C\|\nabla u\|^2_{H^{s_0}}\left(1+\|\nabla
u\|^{3n+2l+12}_{H^{s_0+2}}\right) \nonumber
\end{eqnarray}
\end{proof}



\begin{corollary}[Uniform energy estimate]\label{in-cor2}
Let $u_\e(t)\in H^{s+1}(\RR^n,N)$, with $s\in\NN$ and $s\ge
\left[\frac{n}{2}\right]+4$ be a solution of (\ref{in-eqn1}). There exists
$T_0=T_0(\|\nabla u_0\|_{H^s})$ such that for $0\le t\le T_0$
\begin{equation}\label{in-eqn76}
\|\nabla u_\e(t)\|_{H^{s}} \le 3\|\nabla u_0\|_{H^{s}}.
\end{equation}
\end{corollary}

\begin{proof} 
Let $E(t)=\|\nabla u\|^2_{H^{s}}(t)$. Then (\ref{in-eqn5.1ast}),
(\ref{in-eqn63}) and (\ref{in-eqn64}) imply 
for 
$\left[\frac{n}{2}\right]+4\le s$
\begin{equation}\label{in-eqn77}
\frac{d}{dt}E\le C_0E(1+E^{2n+s+8}),
\end{equation}
which leads, after integrating from $0$ to $t$, to 
\begin{equation}\label{in-eqn78}
\ln\frac{E(t)}{E(0)} 
- \frac{1}{2n+s+8} \ln\frac{E(t)^{2n+s+8}{+1}}{E(0)^{2n+s+8}+1}\le C_0t
\end{equation}
which implies
\begin{eqnarray}\label{in-eqn79}
\frac{E(t)^{2n+s+8}}{1+E(t)^{2n+s+8}}& \le  & e^{C_0t(2n+s+8)}
\frac{E(0)^{2n+s+8}}{1+E(0)^{2n+s+8}}\\
& \le & (1+4C_0t(2n+s+8))
\frac{E(0)^{2n+s+8}}{1+E(0)^{2n+s+8}},\nonumber
\end{eqnarray}
for $t$ such that $C_0t(2n+s+8)<\frac{1}{8}$ for example.
A simple computation yields
\begin{eqnarray}\label{in-eqn79A}
E(t)^{2n+s+8}&\le &(1+4C_0t(2n+s+8))E(0)^{2n+s+8} \\
&&\qquad + 
4C_0t(2n+s+8)E(0)^{2n+s+8}E(t)^{2n+s+8}.\nonumber
\end{eqnarray}
For $t$ such that $4C_0t(2n+s+8) E(0)^{2n+s+8}<\frac{1}{2}$ we have
\begin{equation}\label{in-eqn79B}
E(t)^{2n+s+8}\le 2(1+4C_0t(2n+s+8))E(0)^{2n+s+8}.
\end{equation}
Thus for $s\ge\left[\frac{n}{2}\right]+4$, (\ref{in-eqn79B}) shows that if
\begin{equation}\label{in-eqn85}
0<t\le T_0=\min\{\frac{1}{8C_0(2n+s+8)}, \frac{1}{8C_0(2n+s+8)}E(0)^{2n+s+8}\}
\end{equation}
then


\begin{equation}\label{in-eqn86}
\|\nabla u_\e(t)\|_{H^{s}}\le 3\|\nabla u_0\|_{H^{s}}.
\end{equation}
\end{proof}

\begin{lemma}\label{in-lem3}
Let $u_\e(t)\in H^{s+1}(\RR^n,N)$ with $s\in\NN$ and
$s\geq \left[\frac{n}{2}\right]+4$ be a solution of (\ref{in-eqn5.1}). Let
$v=v_\e=w\circ u_\e$. For $T_0=T_0(\|\nabla u_0\|_{H^s})$ as in (\ref{in-eqn85}) we
have
\begin{equation}\label{in-eqn86A}
\sup_{0<t\le T_0}\|v(t)-v_0\|_{L^2}\le C\|\nabla
u_0\|_{H^{\left[\frac{n}{2}\right]+4}}\left(1+\|\nabla
u_0\|^{3\left[\frac{n}{2}\right]+6}_{H^{\left[\frac{n}{2}\right]+4}}\right)T_0.
\end{equation}
\end{lemma}

\begin{proof}
Our goal is to study how $\|v(t)-v_0\|_{L^2}$ evolves. 
Using (\ref{du-eqn17}) and (\ref{du-eqn18}) we have
\begin{eqnarray}\label{in-eqn87}
\frac{1}{2}\frac{d}{dt}\int|v-v_0|^2 
&= &\int\langle\partial_tv,v-v_0\rangle\\
 &
\le & \|v-v_0\|_{L^2}\left(\int(\partial_t
v)^2\right)^{\frac{1}{2}}\nonumber \\
& \le & C\left(\|\Delta^2v\|_{L^2}+ \|\partial^2v\|_{L^2} + \|\partial
v\|_{L^2}\|\partial v\|_{L^\infty}\right. \nonumber \\
&&\qquad + \|\partial^2v\|_{L^2} \|\partial^2v\|_{L^\infty} + \|\partial
v\|_{L^2} \|\partial v\|^3_{L^\infty} \nonumber \\
&&\left.\qquad + \|\partial^2v\|_{L^2} \|\partial v\|^2_{L^\infty}\right)
\|v-v_0\|_{L^2}. \nonumber
\end{eqnarray}
Recall that
\begin{equation}\label{in-eqn87a}
|\partial^2 v|\le |\nabla^2 u|+C|\nabla u|^2\mbox{ and }|\partial v|=|\nabla
u|.
\end{equation}
Moreover by (\ref{an-eqn36a}) we have
\begin{eqnarray}\label{in-eqn105A}
|\partial^4v| & \le & C|\nabla^4u|+C\sum^4_{l=2} \sum_{j_1+\cdots+j_l=4}
|\nabla^{j_1}u|\cdots|\nabla^{j_l}u| \\
& \le & C|\nabla^4u|+C|\nabla^2u|^2+C|\nabla u|^4+C|\nabla^3u|\,|\nabla
u|.\nonumber 
\end{eqnarray}
Using (\ref{an-eqn58}), (\ref{an-eqn59}) and (\ref{in-eqn105A}),
(\ref{in-eqn87}) yields
\begin{eqnarray}\label{in-eqn106}
&&\kern-.5in\frac{d}{dt}\int|v-v_0|^2 \\
 & \le & C\left\{\|\nabla^4u\|_{L^2} +
\|\nabla^2u\|_{L^\infty}\|\nabla^2u\|_{L^2}\right.\nonumber \\
&&\qquad \|\nabla u\|^3_{L^\infty} \|\nabla u\|_{L^2} + \|\nabla
u\|_{L^\infty} \|\nabla^3u\|_{L^2} + \|\nabla^3u\|_{L^2} \nonumber \\
&&\left.\qquad + \|\nabla u\|_{L^2}\|\nabla u\|_{L^\infty} + \|\nabla^2u\|_{L^2}
\|\nabla u\|_{L^\infty}^2\right\} \|v-v_0\|_{L^2} \nonumber \\
& \le & C\|\nabla u\|_{H^{\left[\frac{n}{2}\right]+3}} \left(1+\|\nabla
u\|^{3\left[\frac{n}{2}\right]+6}_{H^{\left[\frac{n}{2}\right]+3}}\right)
\|v-v_0\|_{L^2}.\nonumber
\end{eqnarray}
For $t\in[0,T_0]$ as in (\ref{in-eqn85}),  \eqref{in-eqn106} combined with
(\ref{in-eqn76}) yields
\begin{equation}\label{in-eqn107}
\frac{d}{dt}\|v-v_0\|^2_{L^2} \le C\|\nabla
u_0\|_{H^{\left[\frac{n}{2}\right]+4}} \left(1+\|\nabla
u_0\|^{3\left[\frac{n}{2}\right]+6}_{H^{\left[\frac{n}{2}\right]+4}}\right)
\|v-v_0\|_{L^2}.
\end{equation}
Integrating from $0$ to $T_0$ (as defined in (\ref{in-eqn85})) we deduce from
(\ref{in-eqn107}) that
\begin{equation}\label{in-eqn107A}
\|v(t)-v_0\|_{L^2} \le CT_0\|\nabla u_0\|_{H^{\left[\frac{n}{2}\right]+4}}
\left(1+\|\nabla
u_0\|^{3\left[\frac{n}{2}\right]+6}_{H^{\left[\frac{n}{2}\right]+4}}\right).
\end{equation}
\end{proof}

\begin{theorem}\label{in-thm1}
Let $s\ge \left[\frac{n}{2}\right]+4$. 
Given $u_0\in H^{s+1}(\RR^n,N)$ there exists
$T_0=T_0(\|\nabla u_0\|_{H^s}, N)>0$ and a solution $u_\e\in C([0,T_0],
H^{s+1}(\RR^n,N))$ of (\ref{in-eqn1}). Furthermore
\begin{equation}\label{in-eqn100}
\sup_{0\le t\le T_0} \|\nabla u_\e(t)\|_{H^s}\le 3\|\nabla u_0\|_{H^s}.
\end{equation}
\end{theorem}

\begin{proof}
Lemma \ref{fo-lem1}, Lemma \ref{an-lem8}, Theorem \ref{du-thm2} and Lemma 
\ref{an-lem8A} imply that there exist
$T_\e=T(\e,\|\nabla u_0\|_{H^s}, \|v_0-\gamma\|_{L^2},N)$ for some
$\gamma\in\RR^p$, and a solution of (\ref{in-eqn1}) given by
$u_\e\in C([0,T_\e],
H^{s+1}(\RR^n,N))$. Either $T_\e\ge T_0$ as
defined in (\ref{in-eqn85}) and we are done or $T_\e<T_0$. Using the fact
that 
\[
\|v(T_\e)-v_0\|_{L^2}\le CT_0\|\nabla
u_0\|_{H^{\left[\frac{n}{2}\right]+4}} \left(1+\|\nabla
u_0\|^{n+2}_{H^{\left[\frac{n}{2}\right]+4}}\right)
\] 
the same argument as
above ensures that there exists $T'_\e=T(\e,\|\nabla u_0\|_{H^s})$ and
$u_\e\in C([T_\e, T_\e+T'_\e], H^s(\RR^n,N))$ a solution of (\ref{in-eqn1}).
The uniqueness statement in Theorem \ref{du-thm2} ensures that we can extend $u_\e\in
C([0, T_\e+T'_\e], H^s(\RR^n,N))$ to be a solution of (\ref{in-eqn1}).

After a finite number of steps (namely $l$ where $T_\e+lT'_\e\le
T_0<T_\e+(l+1)T'_\e)$ we manage to extend $\forall\e\in(0,1)$, $u_\e$ to be a
solution of (\ref{in-eqn1}) in $C([0,T_0], H^s(\RR^n,N))$. 
Note that (\ref{in-eqn100})
is simply a restatement of (\ref{in-eqn86}).
\end{proof}

{\bf Proof of Theorem \ref{int-thm1}.}
For $s\ge\left[\frac{n}{2}\right]+4$, let $u_\e\in C([0,T_0],
H^{s+1}(\RR^n,N))$ be a solution of (\ref{in-eqn5.1}). 
Choosing a sequence $\e_i\to 0$ we conclude, by means of 
Theorem \ref{in-thm1} and Lemma \ref{an-lem8A}, \ref{fo-lem1} that 
there exist functions
$u\in C([0,T_0], H^{s+1}(\RR^n,N))$ 
and $v\in C([0,T_0], H^{s+1}(\RR^n,\RR^m))$ with $v=\omega\circ u$ 
satisfying the initial value problems
(\ref{du-eqn17}) and (\ref{fo-eqn13})
with $\e=0$ and $v_0=\omega\circ u_0$. 

To prove the well-posedness of the \sh flow 
(i.e. when $\beta=0$ in (\ref{int-eqn3}))
we refer to work of Ding and Wang \cite{DW1998} and McGahagan \cite{McG1}.
By adapting the argument of Ding and Wang \cite{DW1998}
one can show that if a solution, $u \in C([0,T_0],H^{s+1}(\RR^n,N))$ with
$s\ge\left[\frac{n}{2}\right]+4$,
to the initial value problem (\ref{int-eqn3}) (with $\beta=0$) exists then it 
is unique.  This argument makes explicit use of the fact that the target 
is compact and isometrically embedded into some Euclidean space. 
We present and extend here part of an argument that appears in the proof of
Theorem 4.1 in \cite{McG1}. These inequalities yield uniqueness and 
continuous dependence on the initial data for general $\beta \ge 0$. Let 
$u_1, u_2\in C([0,T_0], H^{s+1}(\RR^n,N))$ be solutions of 
(\ref{int-eqn3}) with initial data 
$u_1^0, u_2^0\in H^{s+1}(\RR^n,N)$ with
$s\ge\left[\frac{n}{2}\right]+4$. Following the notation in \cite{McG1}
let $V=\nabla u_1$ and $W=\nabla u_2$. Let $\widetilde V(x)$ represent 
the parallel transport of $V$ to the point $u_2(x)$ along the unique geodesic joining the
points. McGahagan proves
(see end of the proof of Theorem 4.1 in \cite{McG1}) that
whenever $\|u_1^0-u_2^0\|_{H^{\left[\frac{n}{2}\right]+4}}$ is small enough
(depending only on the geometry of $N$) and $\beta=0$, then
\begin{equation}\label{tt7}
\frac{d}{dt}\biggl(\|W-\widetilde V\|^2_{L^2} +\|u_1-u_2\|^2_{L^2}\biggr)
\le C
\biggl (\|W-\widetilde V\|^2_{L^2} +\|u_1-u_2\|^2_{L^2}\biggr ),
\end{equation}
where $C$ depends on the $H^{\left[\frac{n}{2}\right]+4}$ norms of 
$u_1$ and $u_2$. In the case that $u_1^0=u_2^0$ McGahagan concludes (using
Gronwall's) that $\|W-\widetilde V\|^2_{L^2} =\|u_1-u_2\|^2_{L^2}=0$,
and that therefore $u_1=u_2$ a.e.. In appendix $A$ we show that the inequality \eqref{tt7} (and therefore also the uniqueness result) remains true for all $\beta\ge 0$.
 
Since the unique solution is constructed as a limit 
of solutions of equation (\ref{in-eqn1}) letting $\e\rightarrow 0$, 
the estimate in Theorem \ref{in-thm1} yields that 
\begin{equation}\label{tt14}
\sup_{0\le t\le T_0} \|\nabla u(t)\|_{H^s}\le 3 \|\nabla u_0\|_{H^s}.
\end{equation}

To prove the continuous dependence on the initial data note 
that, in general, (\ref{tt7})
yields
\begin{equation}\label{tt8}
\|W-\widetilde V\|^2_{L^2}(t) +\|u_1-u_2\|^2_{L^2}(t)
\le e^{Ct}\biggl(
\|W^0-\widetilde V^0\|^2_{L^2} +\|u_1^0-u_2^0\|^2_{L^2}\biggr ),
\end{equation}
where $W^0=\nabla u_2^0$ and $\widetilde V^0(x)$ is the parallel transport of
$V^0=\nabla u_1^0$ to $u_2^0(x)$.
Since
\begin{equation}\label{tt9}
\|W-\widetilde V\|^2_{L^2}(t)\lesssim \|\partial u_1-\partial u_2\|^2_{L^2}(t)
 +\|u_1-u_2\|^2_{L^2}(t),
\end{equation}
and
\begin{equation}\label{tt10}
 \|\partial u_1-\partial u_2\|^2_{L^2}(t) \lesssim
\|W-\widetilde V\|^2_{L^2}(t)
 +\|u_1-u_2\|^2_{L^2}(t),
\end{equation}
(\ref{tt8}) yields
\begin{equation}\label{tt11}
 \|\partial u_1-\partial u_2\|^2_{L^2}(t) +\|u_1-u_2\|^2_{L^2}(t)
\lesssim e^{Ct}\biggl (
\|\partial u_1^0-\partial u_2^0\|^2_{L^2} +\|u_1^0-u_2^0\|^2_{L^2}\biggr ).
\end{equation}
Note that (\ref{tt11}) ensures that $C([0,T_0], H^{s+1}(\RR^n,N))$
solutions to (\ref{int-eqn3}) with $s\ge \left[\frac{n}{2}\right]+4$ depend continuously in $H^1$ 
on the initial data. To show continuous dependence in  
$H^{s'}$ for $s'<s$ we need to use a classic interpolation inequality in 
$\RR^n$. If $v_i=\omega\circ u_i$ for $i=1,2$, where $\omega$ 
denotes the embedding of 
$N$ into $\RR^p$ then combining (\ref{an-eqn34a}) and (\ref{tt11}) we
have
\begin{equation}\label{tt12}
 \|\partial v_1-\partial v_2\|^2_{L^2}(t) +\|v_1-v_2\|^2_{L^2}(t)
\lesssim e^{Ct}\biggl (
\|\partial v_1^0-\partial v_2^0\|^2_{L^2} +\|v_1^0-v_2^0\|^2_{L^2}\biggr ).
\end{equation}
Interpolation, Lemma \ref{an-lem8}, Lemma \ref{an-lem8A}, (\ref{tt14})
and (\ref{tt12}) yield for $s'<s$
\begin{eqnarray}\label{tt13}
&&\|\partial v_1-\partial v_2\|_{H^{s'}}(t) \lesssim
\|\partial v_1-\partial v_2\|^{\frac{s'}{s}}_{H^{s}}(t)
\|\partial v_1-\partial v_2\|^{1-\frac{s'}{s}}_{L^2}(t)\\
&& \qquad\lesssim  \biggl(\|\partial v_1\|^{\frac{s'}{s}}_{H^{s}}(t) 
+ \|\partial v_2\|^{\frac{s'}{s}}_{H^{s}}(t)\biggr)
\|\partial v_1-\partial v_2\|^{1-\frac{s'}{s}}_{L^2}(t)\nonumber\\
&&\qquad \lesssim  \biggl(\|u_1^0\|_{H^s}^m + \|u_2^0\|_{H^s}^m\biggr)
\|\partial v_1-\partial v_2\|^{1-\frac{s'}{s}}_{L^2}(t)\nonumber\\
&&\qquad \lesssim  \biggl(\|u_1^0\|_{H^s}^m + \|u_2^0\|_{H^s}^m\biggr)
e^{Ct}\biggl(
\|\partial v_1^0-\partial v_2^0\|_{L^2} +\|v_1^0-v_2^0\|_{L^2}\biggr)^{1-\frac{s'}{s}}
\nonumber
\end{eqnarray}
Inequalities (\ref{tt12}) and (\ref{tt13}) prove that if 
$u_1, u_2 \in C([0,T_0], H^{s+1}(\RR^n,N))$  are solutions to 
 (\ref{int-eqn3}) and 
$\|u_1^0-u_2^0\|_{H^{\left[\frac{n}{2}\right]+4}}$ is small enough
then the functions $v_1=\omega \circ u_1$ and $v_2=\omega \circ u_2$,
which are solutions to the ambient equation, depend continuously in the 
$H^{s'+1}(\RR^n,\RR^p)$-norm on the initial data for $s'<s$.
As mentioned in the introduction 
by means of the standard Bona-Smith regularization procedure
(\cite{B-S,I-I,Ke}) one can prove that the dependence on the initial data is 
continuous in $H^{s+1}(\RR^n,\RR^p)$. It is in this sense that we express the 
well-posedness of (\ref{int-eqn3}). This concludes the proof  
Theorem \ref{int-thm1}.\qed

\appendix
\numberwithin{equation}{section}
\section{Proof of \eqref{tt7} for $\beta\ge 0$}

Since the proof follows closely the one of Theorem $4.1$ in \cite{McG1} we only sketch the main ideas here. 

We let $u_1,u_2 \in C([0,T_0],H^{s+1}(\RR^n,N)$ be two solutions of \eqref{int-eqn3} with initial data $u_1^0$ respectively $u_2^0$ and we assume that $||u_1^0-u_2^0||_{H^{\left[\frac{n}{2}\right]+4}}$ is small. As in \cite{McG1} we let $\gamma(x;x,t)$ be the unique length minimizing geodesic (parametrized by arclength $s\in [0,l(x,t)]$) between $u_1(x,t)$ and $u_2(x,t)$, where $\gamma(0;x,t)=u_1(x,t)$ and $\gamma(l(x,t);x,t)=u_2(x,t)$ (the existence of the geodesic follows from the argument on page $392$ in \cite{McG1}; note that this argument is also applicable if $u_1^0$ and $u_2^0$ are only close to each other in $L^n$). Moreover we define $V_k=\partial_k u_1$, $W_k=\partial_k u_2$ and $\bar V_k=X(l,0)V_k$ as the parallel transport of $V_k$ to the point $u_2$. In the following we let $X(l,0)=:X$.

In \cite{McG1}, page $391$, the following commutator formulas are derived: $\forall F\in T_{u_1} N$ we have
\begin{align}
XJ(u_1) F &= J(u_2) X F, \nonumber \\
[D_k,X]F&=\int_0^l X(l,\tau)R(\partial_k \gamma,\partial_s \gamma) X(\tau,0)F d\tau. \label{com1}
\end{align}
Additionally the estimates 
\begin{align}
||[D_t,X]V||_{L^2}&+||[D,X]\partial_t u_1||_{L^2}+||D[D_k,X]V_k||_{L^2}\nonumber\\
\le& c(||W-\bar V||_{L^2}+||u_1-u_2||_{L^2}\label{com2}
\end{align}
have been derived in \cite{McG1} (see estimates ($42$), ($43$) and page $395$).
In the following we also need the fact that
\begin{align}
||[D_k,X]V_k||^2_{L^2}\le& c(||W-\bar V||^2_{L^2}+||u_1-u_2||^2_{L^2})\label{com3}.
\end{align}
In order to see this we note that 
\begin{align*}
||[D_k,X]V_k||^2_{L^2}\le& c|| l \nabla_k \gamma V_k||_{L^2}^2.
\end{align*}
For $n\ge 3$ we can use the Sobolev embedding theorem and H\"older's inequality to get (note that $|\nabla \gamma|\le c(|V|+|W|))$
\begin{align*} 
||[D_k,X]V_k||^2_{L^2}\le& c||\nabla l||^2_{L^2}|||V|(|V|+|W|)||_{L^n}^2\\
\le& c||W-\bar V||_{L^2}^2.
\end{align*}
In the case $n=2$ one argues with the help of the Brezis-Wainger theorem as in \cite{McG1}, page $395$.

Now we are finally able to prove \eqref{tt7}. Since $u_1$ and $u_2$ are both solutions of \eqref{int-eqn3} we get
\begin{align*}
\partial_t u_2-(J(u_2)+\beta)\tau(u_2)-X\Big(\partial_t u_1-(J(u_1)+\beta)\tau(u_1)\Big)=0.
\end{align*}
Using the previous definitions and the commutator formula \eqref{com1} we can rewrite this equation as follows
\begin{align*}
\partial_t u_2 -X\partial_t u_1 -J(u_2)(D_k W_k-X D_k V_k)=\beta(D_k W_k-X D_k V_k).
\end{align*}
Multiplying this equation with $(D_k W_k-X D_k V_k\in T_{u_2}N)$ and integrating we get
\begin{align*}
\int_{\RR^n} \langle \partial_t u_2 -X\partial_t u_1,D_k W_k-X D_k V_k\rangle=\beta \int_{\RR^n} |D_k W_k-X D_k V_k|^2.
\end{align*}
Next we calculate 
\begin{align*}
\int_{\RR^n} \langle &\partial_t u_2 -X\partial_t u_1,D_k W_k-X D_k V_k\rangle\\
=& -\int_{\RR^n} \langle D_k(\partial_t u_2 -X\partial_t u_1),W_k-XV_k\rangle-\int_{\RR^n} \langle \partial_t u_2 -X\partial_t u_1,[X,D_k]V_k\rangle\\
=& -\int_{\RR^n} \langle \partial_t W_k -X\partial_t V_k,W_k-XV_k\rangle-\int_{\RR^n} \langle [D_k,X]\partial_t u_1,W_k-XV_k\rangle -I\\
=&-\frac12 \partial_t \int_{\RR^n}|W_k-XV_k|^2 +\int_{\RR^n} \langle [X,\partial_t]V_k,W_k-XV_k\rangle -I-II\\
=&-\frac12 \partial_t \int_{\RR^n}|W_k-XV_k|^2-I-II+III.
\end{align*}
Combining this with the above equality we conclude
\begin{align*}
\frac12 \partial_t \int_{\RR^n}|W_k-XV_k|^2+\beta \int_{\RR^n} |D_k W_k-X D_k V_k|^2=-I-II+III.
\end{align*}
Next we estimate the three terms on the right hand side. We start with 
\begin{align*}
|II|\le& c||W_k-XV_k||_{L^2}||[D_k,X]\partial_t u_1||_{L^2}\\
\le& c(||W-\bar V||_{L^2}^2+||u_1-u_2||_{L^2}^2),
\end{align*}
where we used \eqref{com2} in the last line. Using the same arguments we also get
\begin{align*}
|III|\le& c||W_k-XV_k||_{L^2}||[X,\partial_t]V_k||_{L^2}\\
\le& c(||W-\bar V||_{L^2}^2+||u_1-u_2||_{L^2}^2).
\end{align*}
In order to estimate $I$ we use equation \eqref{int-eqn3}, the fact that $\nabla J=0$ and \eqref{com1} to rewrite
\begin{align*} 
-I=& -\int_{\RR^n} \langle \partial_t u_2 -X\partial_t u_1,[X,D_k]V_k\rangle\\
=&-\int_{\RR^n} \langle (J(u_2)+\beta_\tau(u_2)-X(J(u_1)+\beta)\tau(u_1),[X,D_k]V_k\rangle\\
=&-\int_{\RR^n} \langle (J(u_2)+\beta)\big(D_kW_k-D_k(XV_k)-[X,D_k]V_k\big),[X,D_k]V_k\rangle\\
=&\int_{\RR^n} \langle (J(u_2)+\beta)(W_k-XV_k),D_k([X,D_k]V_k)\rangle+\beta \int_{\RR^n} |[X,D_k]V_k|^2.
\end{align*}
Using H\"older's inequality, \eqref{com2} and \eqref{com3} we get
\begin{align*}
|I|\le& c(||W-\bar V||_{L^2}^2+||u_1-u_2||_{L^2}^2).
\end{align*}
Altogether this implies  
\begin{align}
\frac12 \partial_t \int_{\RR^n}|W-\bar V|^2+\beta \int_{\RR^n} |D_k W_k-X D_k V_k|^2\le c(||W-\bar V||_{L^2}^2+||u_1-u_2||_{L^2}^2).\label{der1}
\end{align}
Next we need to estimate $\frac12 \partial_t \int_{\RR^n}|u_1-u_2|^2$. In order to do this we argue as in \cite{McG1} and we consider $N$ to be isometrically embedded into $\RR^p$ and we extend $J$ as a continuous linear operator on $\RR^p$. In the following we denote the second fundamental form of the embedding $N\hookrightarrow \RR^p$ by $A$. With the help of these conventions we calculate
\begin{align*}
\frac12 &\partial_t \int_{\RR^n}|u_1-u_2|^2= \int_{\RR^n} \langle \partial_t(u_1-u_2), u_1-u_2 \rangle\\
=&\int_{\RR^n}\langle J(u_1)\Delta u_1 -J(u_2)\Delta u_2, u_1-u_2\rangle+\beta \int_{\RR^n}\langle \Delta u_1 -\Delta u_2, u_1-u_2\rangle\\
&+\int_{\RR^n}\langle (J(u_1)+\beta)A(u_1)(\nabla u_1,\nabla u_1)-(J(u_2)+\beta)A(u_2)(\nabla u_2, \nabla u_2), u_1-u_2\rangle.
\end{align*}
Arguing as in \cite{McG1}, page $296$, we get the estimate
\begin{align*}
\int_{\RR^n}\langle (J(u_1)+\beta)A(u_1)(\nabla u_1,\nabla u_1)&-(J(u_2)+\beta)A(u_2)(\nabla u_2, \nabla u_2), u_1-u_2\rangle\\
&+\int_{\RR^n}\langle J(u_1)\Delta u_1 -J(u_2)\Delta u_2, u_1-u_2\rangle\\
\le& c(||W-\bar V||_{L^2}^2+||u_1-u_2||_{L^2}^2).
\end{align*}
Moreover we note that
\begin{align*}
\beta \int_{\RR^n}\langle \Delta u_1 -\Delta u_2, u_1-u_2\rangle=-\beta ||\nabla u_1-\nabla u_2||_{L^2}^2
\end{align*}
and hence we conclude
\begin{align}
\frac12 \partial_t \int_{\RR^n}|u_1-u_2|^2+\beta ||\nabla u_1-\nabla u_2||_{L^2}^2\le& c(||W-\bar V||_{L^2}^2+||u_1-u_2||_{L^2}^2).\label{der2}
\end{align}
Combining \eqref{der1} and \eqref{der2} finishes the proof of of \eqref{tt7} for general $\beta \ge 0$.

\end{document}